\documentclass[12pt]{article}

\usepackage{amsmath}
\usepackage{amsthm}
\usepackage{amssymb}
\usepackage{amscd}
\usepackage{xspace}
\usepackage{verbatim}
\usepackage{geometry}
\usepackage{esint}
\usepackage{geometry}
\usepackage{fancyhdr}
\usepackage{xargs}
\usepackage{xifthen}
\usepackage{xcolor}

\newcommandx{\Set}[2][2=]{
    \ifthenelse{\isempty{#2}}
        {\left\{ {#1} \right\}}
        {\left\{  {#1}  \, \middle| \, {#2} \right\}}
}

\newtheorem{theorem}{Theorem}[section]
\newtheorem*{theorem*}{Theorem}
\newtheorem{corollary}[theorem]{Corollary}
\newtheorem{lemma}[theorem]{Lemma}
\newtheorem{proposition}[theorem]{Proposition}
\newtheorem*{proposition*}{Proposition}
\newtheorem{example}[theorem]{Example}

\theoremstyle{definition}
\newtheorem{definition}[theorem]{Definition}
\theoremstyle{remark}
\newtheorem{remark}{Remark}[section]
\numberwithin{equation}{section}
\newtheorem{question}{Question}[section]

\newcommand{\e}{\varepsilon}

\newcommand{\N}{\mathbb N}

\renewcommand{\vec}[1]{\mathbf{#1}}

\newcommand{\R}{\mathbb{R}}

\newcommand{\er}{\eqref}

\DeclareMathOperator{\Div}{div} \DeclareMathOperator{\dist}{dist}
\DeclareMathOperator{\supp}{supp}

\DeclareMathOperator{\diam}{diam}

\newcommandx{\norm}[1][1=\cdot]{\left|{#1}\right|}
\newcommandx{\supnorm}[2][1=\cdot,2=]{\nnorm[#1]_{\infty,#2}}
\newcommandx{\nnorm}[1][1=\cdot]{\left|\left|{#1}\right|\right|}

\newcommand{\Haus}{\mathcal H}

\renewcommand{\vec}[1]{\boldsymbol{#1}}
\date{\today}

\begin{document}
\title{BMO-Interpolation and Jump Detection for Functions in $BV\cap BMO$}

\author{Paz Hashash and Arkady Poliakovsky}
 
\maketitle
\date{}

\begin{abstract}
A generalization of the John--Nirenberg inequality is established. 
As a consequence, local and global $BMO$--interpolation inequalities in Lorentz spaces $L^{q,\gamma}$ are obtained for the full range $0<q<\infty$ and $0<\gamma\leq\infty$. 
These inequalities yield interpolation results in Besov spaces, fractional Sobolev spaces, and the space $BV$, including corresponding weak spaces. As a geometric application, consequences for the jump set of functions in $BV\cap BMO$ are derived. 
\end{abstract}
\tableofcontents

\section{Introduction}
Interpolation inequalities involving the space $BMO$ play an important role in harmonic analysis, partial differential equations, and geometric measure theory. The space $BMO(\mathbb{R}^N)$ was introduced by John and Nirenberg \cite{JN1961} and arises naturally as a limiting substitute for $L^\infty$ in several critical embedding problems.

A classical example is provided by the Sobolev embedding theorem. If $1\le p<N$, then
\[
W^{1,p}(\mathbb{R}^N)\subset L^{p^*}(\mathbb{R}^N),
\qquad
p^*=\frac{Np}{N-p},
\]
see \cite[Section 4.5.1]{EG}. In the critical case $p=N$, the embedding into $L^\infty(\mathbb{R}^N)$ fails, while
\[
W^{1,N}(\mathbb{R}^N)\subset BMO(\mathbb{R}^N),
\]
see \cite[Section 5.8.1]{Evans2022}. Another manifestation of the limiting character of $BMO$ is the John--Nirenberg inequality, which implies exponential integrability of $BMO$ functions. In particular, for every cube $Q\subset\mathbb{R}^N$ and every $p\in(0,\infty)$,
\[
L^\infty(Q)\subset BMO(Q)\subset L^p(Q),
\]
see Proposition \ref{prop:BMO-function is integrable of any order on a cube}.

A basic problem is to determine how bounded mean oscillation interacts with other regularity assumptions. This leads to interpolation inequalities of the form
\[
\|u\|_{X_1}\leq C\,\|u\|_{X_0}^{\theta}\|u\|_{BMO}^{1-\theta},
\]
where $X_0$ and $X_1$ are function spaces and $\theta\in(0,1)$. Such inequalities occur in connection with Sobolev embeddings, regularity theory, and harmonic analysis; see \cite{BrezisGallouet1980,FeffermanStein1972,KozonoTaniuchi2000,Stein1993}. Interpolation results involving $BMO$ and Lebesgue spaces were obtained in \cite{ChenZhu2005,Dao2018,KozonoWadade2008,McCormick2013}. More recently, Brezis and Mironescu \cite{BrezisMironescu2021} and Van Schaftingen \cite{VanSchaftingen2022} established interpolation inequalities for fractional Sobolev spaces.

The purpose of this paper is to study interpolation phenomena involving $BMO$ in weak Lebesgue spaces, Lorentz spaces, fractional Sobolev spaces, Besov spaces, and the space of functions of bounded variation. Our starting point is the John--Nirenberg inequality. We first establish a generalization of this, and then derive an interpolation inequality in Lorentz spaces. This result serves as the main tool from which the remaining interpolation inequalities are derived, yielding a unified treatment of several classes of function spaces and recovering a number of known results as special cases.

The method employed here differs from the classical interpolation-theoretic approach. It is based primarily on the exponential decay estimate supplied by the John--Nirenberg inequality rather than on the real or complex interpolation methods \cite{BerghLofstrom1976,Calderon1964}. For some results we also provide alternative proofs based on the real interpolation method and the theory of function spaces.

Weak variants of the classical function spaces arise naturally in terms of distribution-function estimates. 
For a measurable function $u$, the distribution function is defined by
\[
d_u(t) := \mathcal{L}^N(\{|u|>t\}), \qquad t \ge 0.
\] The weak Lebesgue space $L^q_w$ is characterized by the decay of the distribution function of $u$, while weak fractional Sobolev and Besov spaces may be formulated in terms of the distribution functions of the differences $u(x+h)-u(x)$. These formulations are particularly convenient for the study of interpolation inequalities and remain meaningful even when the corresponding strong quasi-norms are infinite.
\\
---
\\
We now turn to a more rigorous formulation. 

The following theorem was proved by John and Nirenberg in \cite{JN1961}:
\begin{theorem}[Lemma 1' in \cite{JN1961}]
\label{thm:BMO inequality}
Let $C$ be a cube and $u\in BMO(C)$, and
let us denote for $\sigma\in [0,\infty)$
\begin{equation}
S_\sigma:=\Set{x\in C}[|u(x)-u_C|>\sigma],\quad k:=\|u\|_{BMO(C)},\quad u_C:=\frac{1}{\mathcal{L}^N(C)}\int_C\,u(x)dx.
\end{equation}
Then,
\begin{equation}
\mathcal{L}^N\left(S_\sigma\right)
\leq \frac{A}{k}e^{-\alpha\sigma k^{-1}}\int_{C}|u(x)-u_{C}|dx\quad \text{for}\quad \frac{\sigma}{k}\geq a.
\end{equation}
Here $A\leq 1$, $\alpha,a$ are positive numbers depending only on the dimension $N$.
\end{theorem}

Our first results are the following two theorems, which are a useful generalization of the John--Nirenberg inequality:
\begin{theorem}[Local version]
\label{thm:relation between Lebesgue measures of level sets of BMO} 
Let $Q\subset\R^N$ be a cube and $u\in BMO(Q)$. 
Define for $\sigma\in [0,\infty)$
\begin{equation}
S_\sigma:=\Set{x\in Q}[|u(x)-u_Q|>\sigma].
\end{equation}
Then,
\begin{equation}\label{jkhhhhhjrkjkjhoiiofg}
\mathcal{L}^N\left(S_{\sigma}\right)
\leq \gamma \, e^{-\alpha\frac{\sigma-\eta}{\|u\|_{BMO(Q)}}}\mathcal{L}^N\left(S_{\eta}\right)\quad\quad\forall 0\leq\eta\leq \sigma\,.
\end{equation}
Here, $\gamma=\gamma(N)$ and $\alpha=\alpha(N)$ are positive constants depending only on the dimension $N$.
\end{theorem}

\begin{theorem}[Global version]
\label{thm:relation between Lebesgue measures of level sets of BMO,global} 
Let $u\in BMO(\R^N)$. 
Define for $\sigma\in [0,\infty)$
\begin{equation}
S_\sigma:=\Set{x\in \R^N}[|u(x)|>\sigma].
\end{equation}
Then,
\begin{equation}\label{jkhhhhhjrkjkjhoiiofg333}
\mathcal{L}^N\left(S_{\sigma}\right)
\leq \gamma \, e^{-\alpha\frac{\sigma-\eta}{\|u\|_{BMO(\R^N)}}}\mathcal{L}^N\left(S_{\eta}\right)\quad\quad\forall 0\leq\eta\leq \sigma\,.
\end{equation}
Here, $\gamma=\gamma(N)$ and $\alpha=\alpha(N)$ are positive constants depending only on the dimension $N$.
\end{theorem}

From Theorems \ref{thm:relation between Lebesgue measures of level sets of BMO}  and \ref{thm:relation between Lebesgue measures of level sets of BMO,global}  we derive the following two theorems, which will be used to prove interpolation theorems. They concern Lorentz spaces $L^{q,\gamma}$. We also use the weak Lebesgue space notation $L^q_w := L^{q,\infty}$, $[\cdot]_{L^q_w}:=\|\cdot\|_{L^{q,\infty}}$. 

\begin{theorem}(Local version)
\label{thm:the minimum theorem for Lorentz quasi-norms including characteristic function1} 
Let $Q\subset\R^N$ be a cube, $q\in (0,\infty)$, $\gamma\in (0,\infty]$. Let $u\in BMO(Q)$ and $\vec k \in [0,\infty)$ be such that
\begin{equation}
\label{eq:uniform bound on BMO norms137}
\left\|u\right\|_{BMO(Q)}\leq \vec k.
\end{equation}
Then,
\begin{equation}
\vec k \, \mathcal{L}^N\left(Q\cap \Set{\left|u-u_Q\right|>\vec k}\right)^{\frac{1}{q}}\approx \left\|\chi_{\{|u-u_Q|>\vec k\}}(u-u_Q)\right\|_{L^{q,\gamma}(Q)}.
\end{equation}
The constants in the equivalence depend only on $\gamma$, $q$, and $N$.
\end{theorem}

\begin{theorem}(Global version)
\label{thm:the minimum theorem for Lorentz quasi-norms including characteristic function,global,intro} 
Let $q\in (0,\infty)$, $\gamma\in (0,\infty]$. Let $u\in BMO(\R^N)$ and $\vec k \in [0,\infty)$ be such that
\begin{equation}
\label{eq:uniform bound on BMO norms1374}
\left\|u\right\|_{BMO(\R^N)}\leq \vec k.
\end{equation}
Then,
\begin{equation}
\vec k \, \mathcal{L}^N\left(\R^N\cap  \Set{\left|u\right|>\vec k}\right)^{\frac{1}{q}}\approx \left\|\chi_{\{|u|>\vec k\}}u\right\|_{L^{q,\gamma}(\R^N)}.
\end{equation}
The constants in the equivalence depend only on $\gamma$, $q$, and $N$.
\end{theorem}
The above theorems say the following. Let $u \in BMO$ with $\|u\|_{BMO} \le \vec{k}$. Then the Lorentz quasi-norm of the truncation $\chi_{\{|u|>\vec{k}\}}u$ is equivalent, up to constants depending only on $q,\gamma,N$, to $\vec{k}\,\mathcal{L}^N(\{|u|>\vec{k}\})^{1/q}$.
In particular, the $L^{q,\gamma}$-size of this truncation is completely determined by the measure of the super level set $\{|u|>\vec{k}\}$.
\\

We prove the following two  interpolation results using Theorems \ref{thm:the minimum theorem for Lorentz quasi-norms including characteristic function1} and \ref{thm:the minimum theorem for Lorentz quasi-norms including characteristic function,global,intro}:
\begin{theorem}[Local $BMO$-interpolation in Lorentz spaces]
\label{thm:interpolation for BMO in cube,intro}
Let $Q\subset\R^N$ be a cube, and let $0<p< q<\infty$, $\gamma\in (0,\infty]$. Let $u\in BMO(Q)$. Then,
\begin{multline}
\label{thm:interpolation for BMO in cube,intro1}
\|u-u_{Q}\|_{L^{q,\gamma}(Q)}\leq C
\Big\|\min\left\{\left|u-u_{Q}\right|\,,\,\|u\|_{BMO(Q)}\right\}\Big\|_{L^{q,\gamma}(Q)}
\\
\leq C\sup\limits_{0<s\leq \|u\|_{BMO(Q)}^p}\Big[s\mathcal{L}^N\left(\Set{x\in
Q}[\left|u(x)-u_{Q}\right|^p>s]\right)\Big]^{\frac{1}{q}}\|u\|_{BMO(Q)}^{1-\frac{p}{q}}
\\
\leq C[u-u_{Q}]^{\frac{p}{q}}_{L^p_w(Q)}\|u\|_{BMO(Q)}^{1-\frac{p}{q}}.
\end{multline}
Here, $C=C(\gamma,p,q,N)$ is a constant that depends only on $\gamma$, $p$, $q$, and $N$.
\end{theorem}

\begin{theorem}[Global $BMO$-interpolation in Lorentz spaces]
\label{thm:interpolation for BMO1222,intro} 
Let $0< p<q<\infty$, $\gamma\in (0,\infty]$. If $u\in BMO(\R^N)$ and it is not a nonzero constant function almost everywhere, then 
\begin{multline}
\label{eq:weak22,intro}
\|u\|_{L^{q,\gamma}(\R^N)}\leq C
\Big\|\min\left\{\left|u\right|\,,\,\|u\|_{BMO(\R^N)}\right\}\Big\|_{L^{q,\gamma}(\R^N)}
\\
\leq C\sup\limits_{0<s\leq \|u\|_{BMO(\R^N)}^p}\Big[s\mathcal{L}^N\left(\Set{x\in
\R^N}[\left|u(x)\right|^p>s]\right)\Big]^{\frac{1}{q}}\|u\|_{BMO(\R^N)}^{1-\frac{p}{q}}
\\
\leq C[u]^{\frac{p}{q}}_{L^p_w(\R^N)}\|u\|_{BMO(\R^N)}^{1-\frac{p}{q}}.
\end{multline}
Here, $C=C(\gamma,p,q,N)$ is a positive constant that depends only on
$p$, $q$,$\gamma$ and $N$.
\end{theorem}
Notice that Theorems \ref{thm:interpolation for BMO in cube,intro} and
\ref{thm:interpolation for BMO1222,intro} are sharp in the following sense: if $p=q$,
then they would yield 
$\|u\|_{L^{p}} \leq C [u]_{L^{p}_w}$,
which is false in general, since $L^{p}_w$ is not contained in
$L^{p}$. The same observation applies to the interpolation results
involving weak norms proved below.

Note that the first inequality of \eqref{eq:weak22,intro} fails for nonzero constant functions. For example, if $u=1$, then $\|u\|_{L^{q,\gamma}(\R^N)}=\infty$, while 
\[
\left\|\min\left\{\left|u\right|,\,\|u\|_{BMO(\R^N)}\right\}\right\|_{L^{q,\gamma}(\R^N)}=0,
\]
since $\|u\|_{BMO(\R^N)}=0$.

The inequality 
\begin{equation}
\label{eq:weak interpolation proved in two ways}
\|u\|_{L^{q,\gamma}(\R^N)}\leq C[u]^{\frac{p}{q}}_{L^p_w(\R^N)}\|u\|^{1-\frac{p}{q}}_{BMO(\R^N)},
\end{equation} 
which is part of \eqref{eq:weak22,intro}, was proved in~\cite{McCormick2013} for the case $\gamma=q$ and $p > 1$, but not for $0<p\leq 1$. In~\cite{Dao2018}, this interpolation inequality was proven for any $\gamma\in (0,\infty)$ and $p > 1$. See  also \cite{ChenZhu2005,KozonoWadade2008} for a similar $BMO$-interpolation result for Lebesgue spaces.

We prove Theorems \ref{thm:interpolation for BMO in cube,intro} and \ref{thm:interpolation for BMO1222,intro}  in Section \ref{sec:$BMO$-interpolation in Lorentz spaces}. We give an alternative proof of \eqref{eq:weak interpolation proved in two ways} based on results from \cite{JawerthTorchinsky1985}. The idea of this alternative proof was suggested by an anonymous referee and is not due to the authors. See Subsection \ref{subsec: An alternative proof for $BMO$-interpolation in Lorentz spaces via real interpolation and sharp maximal function}.

We immediately get the following known result from Theorem \ref{thm:interpolation for BMO1222,intro} by choosing $\gamma=q$: 
Let $0<p\leq q<\infty$. If $u \in BMO(\R^N)$ and it is not a nonzero constant almost everywhere, then $\|u\|_{L^q(\R^N)}
\leq C \|u\|^{\frac{p}{q}}_{L^p(\R^N)}\|u\|_{BMO(\R^N)}^{\,1-\frac{p}{q}}$, where $C$ depends on $p,q,N$ only.

Let $W^{s,p}$, $s\in(0,1)$, be the fractional Sobolev space. In the following definition we introduce the space $W^{s,p}_w$, where $w$ stands for "weak", and we call it the weak fractional Sobolev space. In Section \ref{sec:interpolation of weak Gagliardo spaces}, it will be shown that $W^{s,p} \subset W^{s,p}_w$.

\begin{definition}
\label{def:weak fractional Sobolev space,intro}
Let $0< q < \infty$ and $s\in (0,1)$. For Lebesgue measurable function $u: \R^N \to \R$, we define
\begin{equation}
[u]_{W^{s,q}_w(\R^N)}:=\sup_{0<t<\infty}\left(\int_{\R^N}t\mathcal{L}^N\Big(\Set{x\in
\R^N}[|u(x+y)-u(x)|^q>t]\Big)\frac{1}{|y|^{sq+N}} \, dy\right)^{\frac{1}{q}}.
\end{equation}
We say that $u \in W^{s,q}_w(\R^N)$ if and only if $u \in L^q_w(\R^N)$ and $[u]_{W^{s,q}_w(\R^N)} < \infty$.
\end{definition}
After investigating the space $W^{s,q}_w(\mathbb{R}^N)$, the authors learned, through an anonymous referee, that in case $q>1$ this space is a particular case of the spaces $\dot{\mathcal{B}}^{s}_{q}(\gamma, r)$ for $r=\infty$ and $\gamma=-sq$, see \cite[Theorem~1.3]{Dominguez2023}; meaning that
\begin{equation}
\label{eq:weak W and B coincide}
W^{s,q}_w(\mathbb{R}^N)=\dot{\mathcal{B}}^{s}_{q}(-sq, \infty).
\end{equation}

The following theorem was proved by Brezis and Mironescu in \cite[Lemma 15.7]{BrezisMironescu2021} in case $p,q\in [1,\infty)$, and in an alternative way by Jean Van Schaftingen for $p,q\in (1,\infty)$ in \cite{VanSchaftingen2022}:

\begin{theorem*}
\label{thm:BMO-interpolation for fractional Sobolev functions by Bresis and Mironescu}
Let $1\leq p < q < \infty$ and $s \in (0,1)$. If $u \in BMO(\mathbb{R}^N) \cap W^{s,p}(\mathbb{R}^N)$, then $u\in W^{\frac{ps}{q},q}(\R^N)$ and  
\begin{equation}
\|u\|_{W^{\frac{ps}{q},q}(\mathbb{R}^N)}
\leq C \|u\|^{\frac{p}{q}}_{W^{s,p}(\mathbb{R}^N)} \|u\|^{1-\frac{p}{q}}_{BMO(\mathbb{R}^N)},
\end{equation}
where $C = C(s,p,q,N)$ is a constant depending only on $p$, $q$, $N$, and $s$.
\end{theorem*}

We generalize this theorem to the weak fractional Sobolev space $W^{s,p}_w$ over the full range $0<p<q<\infty$:

\begin{theorem}[$BMO$-interpolation for weak fractional Sobolev functions $W_w^{s,p}$, $s\in (0,1)$]
\label{thm:BMO-interpolation for fractional Sobolev functions,intro}
Let $0< p < q < \infty$ and $s \in (0,1)$. If $u\in BMO(\R^N)$, then  
\begin{equation}
\label{eq:BMO-interpolation for fractional Sobolev functions,intro}
\|u\|_{W^{\frac{ps}{q},q}(\mathbb{R}^N)}
\leq C [u]^{\frac{p}{q}}_{W_w^{s,p}(\mathbb{R}^N)} \|u\|^{1-\frac{p}{q}}_{BMO(\mathbb{R}^N)}
\leq C \|u\|^{\frac{p}{q}}_{W^{s,p}(\mathbb{R}^N)} \|u\|^{1-\frac{p}{q}}_{BMO(\mathbb{R}^N)},
\end{equation}
where $C = C(p,q,N)$ is a constant depending only on $p$, $q$, and $N$ (and not on $s$). In particular, if $u \in W^{s,p}_w(\mathbb{R}^N) \cap BMO(\mathbb{R}^N)$, then $u \in W^{\frac{ps}{q},q}(\mathbb{R}^N)$.
\end{theorem}
Theorem~\ref{thm:BMO-interpolation for fractional Sobolev functions,intro} says that if $u \in BMO$ and $u$ belongs to the weak space $W^{s,p}_w$, then it has finite norm in the strong space $W^{\frac{ps}{q},q}$.

We give two proofs of Theorem~\ref{thm:BMO-interpolation for fractional Sobolev functions,intro}. The first follows our main approach,  see Section \ref{sec:$BMO$-interpolation in weak fractional Sobolev}, while the second provides an alternative argument. The idea of this alternative proof was suggested by the same anonymous referee and is not due to the authors. Moreover, this approach yields an inequality similar to \eqref{eq:BMO-interpolation for fractional Sobolev functions,intro}, with the norm of the homogeneous Besov space $\dot{B}^0_{\infty,\infty}$ replacing $BMO$ (the latter is embedded in the former), albeit under the restriction $p>1$, see Subsection \ref{An alternative proof for BMO-interpolation in Wsp via real interpolation theory and the spaces}.

The following theorem was proven by Jean Van Schaftigen in \cite{VanSchaftingen2022}:
\begin{theorem}
\label{thm: Van Schaftigen's theorem,intro}
For every $N\in \mathbb{N} \setminus \{0\}$ and every $p \in (1, \infty)$, there exists a constant $C(p) > 0$ such that for every $s \in (1/p, 1)$, every open convex set $\Omega \subset \mathbb{R}^N$ satisfying $\kappa(\Omega) < \infty$, and every function $f \in {W}^{1,sp}(\Omega) \cap BMO(\Omega)$, one has $f \in {W}^{s,p}(\Omega)$ and
\begin{equation}
\int_\Omega \int_\Omega \frac{|f(y) - f(x)|^p}{|y - x|^{N+sp}} \, dy \, dx 
\leq C(p) \frac{\kappa(\Omega)^{sp}}{(sp - 1)(1 - s)} \|f\|_{BMO(\Omega)}^{(1-s)p} \int_\Omega |\nabla f(x)|^{sp}dx.
\end{equation}
Here $f\in BMO(\Omega)$ if $f:\Omega\to \R$ is Lebesgue measurable and
\begin{equation}
\label{def:BMO norm according to Van Schaftigen,intro}
\|f\|_{BMO(\Omega)}:= \sup_{x \in \Omega, \, r > 0} 
 \fint_{\Omega \cap B_r(x)}\fint_{\Omega \cap B_r(x)} |f(y) - f(z)| \, dy \, dz<\infty;
\end{equation}
and 
\begin{equation}
\label{eq:definition of the geometric quantity G,intro}
\kappa(\Omega):= \sup\Set{ 
\frac{\mathcal{L}^N(B_r(x))}{\mathcal{L}^N(\Omega \cap B_r(x))}}[x \in \Omega, \, r \in \big(0, \operatorname{diam}(\Omega)\big)],
\end{equation}
where $\operatorname{diam}(\Omega)$ is the diameter of the set $\Omega$, which is equal to $\infty$ if $\Omega$ is unbounded.
\end{theorem}

Using our method, we prove a version of Theorem \ref{thm: Van Schaftigen's theorem,intro} without the convexity assumption on the domain \( \Omega \). However, we assume that the domain \( \Omega \) is an extension domain for the intersection of Sobolev space with \( BMO \). More precisely:
\begin{theorem}
\label{thm:BMO-interpolation for fractional Sobolev functions, the case s=1 in extension domain,intro}
Let $1<p< \infty$, $s\in (1/p,1)$, and $\Omega\subset \R^N$ an open set. Suppose that $u \in W^{1,sp}(\Omega) \cap BMO(\Omega)$. Assume that there exists a function $f\in W^{1,sp}(\R^N)\cap BMO(\R^N)$ such that $f=u$ almost everywhere in $\Omega$ and that there exists a constant $C=C(N,\Omega,p,s)$ satisfying  
\begin{equation}
\label{eq:extension for BMO and Sobolev at the same time}
\|f\|_{BMO(\R^N)}\leq C \|u\|_{BMO(\Omega)};\quad \|\nabla f\|_{L^{sp}(\R^N)}\leq C \|\nabla u\|_{L^{sp}(\Omega)}. 
\end{equation} 
Then,
\begin{equation}
\label{eq:interpolation for fractional Sobolev functions in BMO, extension domain,intro}
\|u\|_{W^{s,p}(\Omega)}
\leq C\, \|\nabla u\|^{s}_{L^{sp}(\Omega)} \|u\|^{1-s}_{BMO(\Omega)},
\end{equation}
where $C=C(s,p,N,\Omega)$ is a constant depending only on $s$, $p$, $N$, and the domain $\Omega$.
\end{theorem}
Theorem \ref{thm:BMO-interpolation for fractional Sobolev functions, the case s=1 in extension domain,intro} is proven in Section \ref{sec:BMO-interpolation of W1sp into Ws,p}.

\begin{remark}
The interpolation result \eqref{eq:interpolation for fractional Sobolev functions in BMO, extension domain,intro} for extension domains in Theorem \ref{thm:BMO-interpolation for fractional Sobolev functions, the case s=1 in extension domain,intro} can be derived from Theorem \ref{thm: Van Schaftigen's theorem,intro} simply by applying it with $\Omega = \mathbb{R}^N$ and then using the assumption \eqref{eq:extension for BMO and Sobolev at the same time}. However, we employ a completely different approach to prove \eqref{eq:interpolation for fractional Sobolev functions in BMO, extension domain,intro} than the one used in \cite{VanSchaftingen2022}.
\end{remark}

The following theorem concerns weak Besov functions \( (B^s_{p,q})_w \). A precise definition will be given below (see Definition \ref{def:definition of general Besov space1}):

\begin{theorem}[$BMO$-interpolation for Besov functions $B^s_{p,q}$]
\label{$BMO$-interpolation for Besov functions,intro}
Let $0<p<\infty$, $0< q\leq \infty$, $s\in (0,\infty)$ and $0<\lambda<p$. There exists a constant $C=C(N,s,p,q,\lambda)$ such that for every $u\in BMO(\R^N)\cap (B^{s}_{\lambda,q\frac{\lambda}{p}})_w(\R^N)$ it follows that
\begin{equation}
\label{eq:interpolations in Besov Bspq11333}
\begin{cases}
\|u\|_{B^{s\frac{\lambda}{p}}_{p,q}(\R^N)}
\leq C \|u\|^{1-\frac{\lambda}{p}}_{BMO(\R^N)}\left(\|u\|_{(B^{s}_{\lambda,q\frac{\lambda}{p}})_w(\R^N)}\right)^{\frac{\lambda}{p}}\leq C \|u\|^{1-\frac{\lambda}{p}}_{BMO(\R^N)}\left(\|u\|_{B^{s}_{\lambda,q\frac{\lambda}{p}}(\R^N)}\right)^{\frac{\lambda}{p}}\quad  &q<\infty
\\
\|u\|_{B^{s\frac{\lambda}{p}}_{p,\infty}(\R^N)}\leq C \|u\|^{1-\frac{\lambda}{p}}_{BMO(\R^N)}\left(\|u\|_{(B^{s}_{\lambda,\infty})_w(\R^N)}\right)^{\frac{\lambda}{p}}
\leq C \|u\|^{1-\frac{\lambda}{p}}_{BMO(\R^N)}\left(\|u\|_{B^{s}_{\lambda,\infty}(\R^N)}\right)^{\frac{\lambda}{p}}\quad  &q=\infty.
\end{cases}
\end{equation}

\end{theorem}
For closely related interpolation results concerning Besov functions, see \cite{Dominguez2023,DominguezTikhonov2022}.
We will prove a result similar to Theorem~\ref{$BMO$-interpolation for Besov functions,intro}, with the space $\mathrm{bmo}$ in place of $BMO$. The idea was suggested by the same anonymous referee; see Section \ref{sec: bmo(local) interpolation in Besov space}.

As a corollary from Theorem \ref{$BMO$-interpolation for Besov functions,intro} we get

\begin{corollary}
Let $m\in\mathbb{N}$, $0<s<1$, and $1\leq p<\infty$ satisfy
\[
\frac{2}{p}\leq s<\frac{1}{m}.
\]
Then there is a constant $C=C(N,m,p,s)>0$ such that, for every $u\in BMO(\R^N)\cap{W}^{m,sp}(\mathbb{R}^N)$
\[
\|u\|_{W^{ms,p}(\mathbb{R}^N)}
\leq
C\|u\|_{BMO(\mathbb{R}^N)}^{1-s}
\|u\|_{W^{m,sp}(\mathbb{R}^N)}^s.
\]
Here $\|u\|_{W^{ms,p}(\mathbb{R}^N)}$ denotes the Gagliardo seminorm, whereas
$\|u\|_{W^{m,sp}(\mathbb{R}^N)}$ denotes the homogeneous Sobolev seminorm
\[
\|u\|_{W^{m,sp}(\mathbb{R}^N)}
:=
\left(
\sum_{|\alpha|=m}
\|D^\alpha u\|_{L^{sp}(\mathbb{R}^N)}^{sp}
\right)^{1/(sp)}.
\]
\end{corollary}

This corollary may be viewed as a version of Theorem~\ref{thm:BMO-interpolation for fractional Sobolev functions, the case s=1 in extension domain,intro} for Sobolev spaces $W^{m,p}$ of functions with $m$ weak derivatives. However, the range of exponents obtained here is more restrictive: we require $s\geq \frac{2}{p}$, whereas in the case $m=1$ Theorem~\ref{thm:BMO-interpolation for fractional Sobolev functions, the case s=1 in extension domain,intro} holds for all $s>\frac{1}{p}$.
\\

Recall the following characterization of functions of bounded variation. For extensive references on $BV$-functions, see, for example, \cite{AFP,EG,Ziemer}. We also introduce a weak concept of it:
\begin{definition}
Let $u\in L^1(\R^N)$. We write $u\in BV(\R^N)$ if and only if
\begin{equation}
\|u\|_{BV(\R^N)}:=\sup_{h\in \R^N\setminus \{0\}}\frac{\|u(\cdot+h)-u\|_{L^1(\R^N)}}{|h|}<\infty.
\end{equation}
We define also the weak space $BV_w(\R^N)$ to be the set of functions $u\in L^1_w(\R^N)$ such that
\begin{equation}
[u]_{BV_w(\R^N)}:=\sup_{h\in \R^N\setminus \{0\}}\frac{[u(\cdot+h)-u]_{L^1_w(\R^N)}}{|h|}<\infty.
\end{equation}
\end{definition}
Note that $BV(\R^N)\subset BV_w(\R^N)$.
The following interpolation result is about  the strong $BV$ and the weak $BV_w$:
\begin{theorem}[$BMO$-interpolation for $BV_w$-functions]
\label{cor:BMO-interpolation for BV,intro}
Let $1 < q < \infty$. If $u \in BMO(\R^N)$, then 
\begin{equation}
\label{eq:interpolation for BV functions in BMO,intro}
\|u\|_{B^{1/q}_{q,\infty}(\R^N)}
\leq C [u]^{\frac{1}{q}}_{BV_w(\R^N)} \|u\|^{1-\frac{1}{q}}_{BMO(\R^N)}\leq C \|u\|^{\frac{1}{q}}_{BV(\R^N)} \|u\|^{1-\frac{1}{q}}_{BMO(\R^N)},
\end{equation}
where $C=C(q,N)$ is a constant depending only on $q$ and $N$. In particular, if $u \in BMO(\R^N)\cap BV_w(\R^N)$, then $u\in B^{1/q}_{q,\infty}(\R^N)$ for every $q\in (1,\infty)$.
\end{theorem}

Recall the definition of the space $VMO$, the space of functions with vanishing mean oscillation:
\begin{definition}[$VMO$-functions]
\label{def:definition of VMO}
Let $u\in BMO(\R^N)$. We say that $u\in VMO(\R^N)$ if and only if
\begin{equation}
\label{eq:vanishing property,intro}
\lim_{R\to 0^+}\left(\sup_{C,\diam(C)\leq R}\fint_C\fint_C|u(x)-u(z)|dxdz\right)=0,
\end{equation}
where the sumpremum is taken over all cubes $C$ with diameter, $\diam(C)$, at most $R$.
\end{definition}

The following theorem is about functions in $\left(B^s_{p,\infty}\right)_w\cap VMO$:
\begin{theorem}
\label{thm: Besov constants vanish for VMO,intro}
Let $0<p<q<\infty$ and $s\in(0,1)$. If  $u\in \left(B^{s}_{p,\infty}\right)_w(\R^N)\cap VMO(\R^N)$, then
\begin{equation}
\label{eq: Besov constants vanish for VMO ,intro}
\lim_{h\to 0}\frac{\|u_h-u\|_{L^q(\R^N)}^q}{|h|^{sp}}=\lim_{h\to 0}\int_{\R^N}\frac{|u(x+h)-u(x)|^q}{|h|^{sp}}dx=0,\quad u_h(x):=u(x+h).
\end{equation}
\end{theorem}
The result of Theorem \ref{thm: Besov constants vanish for VMO,intro} \emph{is sharp} in the following sense: it is no longer true when $p = q$; see Remark \ref{rem: the limit of Besov constant in case q=p is not necessarily zero}.

The similar theorem deals with functions in $BV_w\cap VMO$:
\begin{theorem}
\label{thm: Besov constants vanish for VMO,intro111}
Let $1<q<\infty$. If  $u\in BV_w(\R^N)\cap VMO(\R^N)$, then
\begin{equation}
\label{eq: Besov constants vanish for VMO ,intro111}
\lim_{h\to 0}\frac{\|u_h-u\|_{L^q(\R^N)}^q}{|h|}=\lim_{h\to 0}\int_{\R^N}\frac{|u(x+h)-u(x)|^q}{|h|}dx=0,\quad u_h(x):=u(x+h).
\end{equation}
\end{theorem}
Again note that the result of Theorem \ref{thm: Besov constants vanish for VMO,intro111} \emph{is sharp}: it is no longer true when $q=1$, since \eqref{eq: Besov constants vanish for VMO ,intro111} in that case implies that $u$ is a constant.

\subsection*{Application to jump detection for functions $BV\cap BMO$}
Recall the notions of approximate jump point and jump set:

\begin{definition}[Approximate Jump Points]
\label{def:approximate jump points,intro}
Let
$\Omega\subset\mathbb{R}^N$ be an open set, $u\in
L^1_{\text{loc}}(\Omega,\mathbb{R}^d)$, and $x\in \Omega$. We say that $x$
is an {\it approximate jump point} of $u$ if and only if there
exist $a,b\in \mathbb{R}^d,\,a\neq b,$ and $\nu\in S^{N-1}:=\Set{z\in\R^N}[|z|=1]$ such
that
\begin{equation}
\label{eq:one-sided approximate limit,intro}
\lim_{\rho\to
0^+}\frac{1}{\rho^N}\left(\int_{B^+_\rho(x,\nu)}|u(z)-a|dz+\int_{B^-_\rho(x,\nu)}|u(z)-b|dz\right)=0,
\end{equation}
where
\begin{equation}
B^+_\rho(x,\nu):=\big\{y\in B_\rho(x):(y-x)\cdot
\nu>0\big\},\quad B^-_\rho(x,\nu):=\big\{y\in
B_\rho(x):(y-x)\cdot \nu<0\big\}.\footnote{The dot symbol here stands, as usual, for the scalar product between vectors $v = (v_1, \dots, v_N)$ and $w = (w_1, \dots, w_N)$ in $\mathbb{R}^N$:
$v \cdot w := \sum_{i=1}^{N} v_i w_i$.}
\end{equation}
The triple $(a,b,\nu)$, uniquely determined by \eqref{eq:one-sided
approximate limit,intro} up to a permutation of $(a,b)$ and the change
of sign of $\nu$, is denoted by $(u^+(x),u^-(x),\nu_u(x))$. The set
of approximate jump points is denoted by $\mathcal{J}_u$.
\end{definition}

The following definition introduces the notion of a kernel. We will use this concept to present our next result, which is Theorem \ref{thm: jumps for BMO and BV functions,intro}.

\begin{definition}(Kernel)
\label{def:kernel,intro}
Let $a\in(0,\infty]$. Let $\rho_\e:(0,\infty)\to
[0,\infty),\e\in (0,a),$ be a family of $\mathcal{L}^1$-measurable
functions. We say that the family $\{\rho_\e\}_{\e\in (0,a)}$ is a
{\it kernel} if it has the following two properties:
\begin{enumerate}
\item $\int_{\R^N}\rho_{\e}(|z|)dz=1$ for every $\e\in (0,a)$,
and
\item for every 
$\delta\in (0,\infty)$
\begin{equation}
\label{eq:defining property of decreasing support property,intro}
\lim_{\e\to 0^+}\int_{\delta}^\infty\rho_\e(r)r^{N-1}dr=0.
\end{equation}
\end{enumerate}
\end{definition}

The following theorem was proved in \cite{HashashPoliakovsky2024}:
\begin{theorem}[Jump detection in $BV\cap B^{1/p}_{p,\infty}$]
\label{thm: limit for Besov integral in terms of jumps,intro}
Let $1<p<\infty$, $u\in BV(\R^N,\R^d)\cap B^{1/p}_{p,\infty}(\R^N,\R^d)$, and $1<q<p$. Then, for every $n\in \R^N$ and every Borel set $B\subset\R^N$ such that $\mathcal{H}^{N-1}(\partial B\cap \mathcal{J}_u)=0$, we have
\begin{equation}
\label{eq:equation2779,intro} 
\lim_{\e\to 0^+}\int_B\chi_{B}(x+\e
n)\frac{|u(x+\e n)-u(x)|^q}{\e}dx =\int_{B\cap
\mathcal{J}_u}\left|u^+(x)-u^-(x)\right|^q|\nu_u(x)\cdot
n|d\mathcal{H}^{N-1}(x),
\end{equation}
and for every kernel $\rho_\e$ (see Definition \ref{def:kernel,intro}), we have
\begin{multline}
\label{eq:BBM formula for BV and BMO}
\lim_{\e\to 0^+}\int_{B}\int_{B}\rho_{\e}(|x-y|)\frac{|u(x)-u(y)|^q}{|x-y|}dydx
\\
=\lim_{\e\to 0^+}\fint_{S^{N-1}}\int_B\chi_{B}(x+\e n)\frac{|u(x+\e n)-u(x)|^q}{\e}dxd\mathcal{H}^{N-1}(n)
\\
=\left(\fint_{S^{N-1}}|z_1|~d\Haus^{N-1}(z)\right)\int_{B\cap \mathcal{J}_u}
\Big|u^+(x)-u^-(x)\Big|^q d\mathcal{H}^{N-1}(x).
\end{multline}
Here, $u^+,u^-$ are the one-sided approximate limits of $u$, $\nu_u$ is the unit normal, and $\mathcal{J}_u$ is the jump set of $u$ (see Definition \ref{def:approximate jump points,intro}).
\end{theorem}
The following geometric result is the same as in Theorem \ref{thm: limit for Besov integral in terms of jumps,intro}, but stated explicitly for functions that have both bounded mean oscillation and bounded variation, i.e., $BMO \cap BV$. It is a straightforward consequence of Theorem \ref{thm: limit for Besov integral in terms of jumps,intro} and Theorem \ref{cor:BMO-interpolation for BV,intro}.

\begin{theorem}[Jump detection in $BV\cap BMO$]
\label{thm: jumps for BMO and BV functions,intro}
Let $1<q<\infty$, and $u\in BV(\R^N)\cap BMO(\R^N)$. Then $u\in B^{1/q}_{q,\infty}(\R^N)$, and for every $n\in \R^N$ and every Borel set $B\subset\R^N$ such that $\mathcal{H}^{N-1}(\partial B\cap \mathcal{J}_u)=0$, we have  \eqref{eq:equation2779,intro},  and for every kernel $\rho_\e$ we have \eqref{eq:BBM formula for BV and BMO}.
\end{theorem}

For extensive motivation regarding the energy-type results \eqref{eq:equation2779,intro} and \eqref{eq:BBM formula for BV and BMO}, see \cite{PAjumps,HashashPoliakovsky2024,Poliakovsky2018}.

\section{Functions of Bounded Mean Oscillation(BMO)}
Throughout the paper, $N$ denotes a natural number. We use the term "cube" to refer to a set that is a Cartesian product of intervals of the form $\Pi_{i=1}^N J_i$, where each $J_i$, for $1 \leq i \leq N$, $i\in\N$, is a non-empty, closed, and bounded interval in $\R$ of the same length.  

For a Lebesgue measurable set $A \subset \R^N$ with finite and positive $N$-dimensional Lebesgue measure, denoted by $\mathcal{L}^N(A)$, and a function $u \in L^1(A)$, we use the notation $u_A$ for the average of $u$ over the set $A$, meaning that  
\[
u_A:=\fint_A u(z)\,dz:=\frac{1}{\mathcal{L}^N(A)}\int_A u(z)\,dz.
\]
Recall the definition of $BMO$-functions:
\begin{definition}[Definition of $BMO$-functions in a cube and $\R^N$]
Let \(C_0\) be a cube and \(u \in L^1(C_0)\). We say that \(u \in BMO(C_0)\) if and only if
\begin{equation}
\sup \left\{ \fint_{C} |u(z) - u_C| \, dz \mid C \subset C_0 \text{ is a cube} \right\}< \infty.
\end{equation}
Similarly, if \(u \in L^1_{\text{loc}}(\mathbb{R}^N)\), we say that \(u \in BMO(\mathbb{R}^N)\) if and only if
\begin{equation}
\sup \left\{ \fint_{C} |u(z) - u_C| \, dz \mid C \subset \mathbb{R}^N \text{ is a cube} \right\} < \infty.
\end{equation}
\end{definition}

\begin{proposition}[Characterization of $BMO$-functions via double averages]
\label{prop:Characterization of $BMO$-functions via double averages}
Let \( u \in L^1(C_0) \). Then, \( u \in BMO(C_0) \) if and only if
\begin{equation}
\sup \left\{ \fint_{C} \fint_{C} |u(z) - u(x)| \, dz \, dx \mid C \subset C_0 \text{ is a cube} \right\} < \infty.
\end{equation}
Similarly, for \( u \in L^1_{\text{loc}}(\mathbb{R}^N) \), \( u \in BMO(\mathbb{R}^N) \) if and only if
\begin{equation}
\sup \left\{ \fint_{C} \fint_{C} |u(z) - u(x)| \, dz \, dx \mid C \subset \mathbb{R}^N \text{ is a cube} \right\} < \infty.
\end{equation}
\end{proposition}

\begin{proof}
By the triangle inequality we get for every cube $C$
\begin{multline}
\label{eq:double integral for proving the characterisation of BMO}
\fint_{C}|u(z)-u_C|dz\leq \fint_{C}\fint_{C}|u(z)-u(x)|dxdz\leq \fint_{C}\fint_{C}|u(z)-u_C|+|u(x)-u_C|dxdz
\\
=\fint_{C}|u(x)-u_C|dx+\fint_{C}|u(z)-u_C|dz=2\fint_{C}|u(z)-u_C|dz.
\end{multline}
Taking the supremum in \eqref{eq:double integral for proving the characterisation of BMO} over all cubes $C \subset C_0$, we obtain the first equivalence. Similarly, taking the supremum in \eqref{eq:double integral for proving the characterisation of BMO} over all cubes $C \subset \mathbb{R}^N$, we obtain the second equivalence.
\end{proof}
\begin{definition}[$BMO$-Semi-norms]
For \(u \in BMO(C_0)\), we define the \(BMO\)-semi-norm by
\begin{equation}
\|u\|_{BMO(C_0)} := \sup \left\{ \fint_{C} \fint_{C} |u(z) - u(x)| \, dz \, dx \mid C \subset C_0 \text{ is a cube} \right\}.
\end{equation}
Similarly, for \(u \in BMO(\mathbb{R}^N)\), we define the \(BMO\)-semi-norm by
\begin{equation}
\|u\|_{BMO(\mathbb{R}^N)} := \sup \left\{ \fint_{C} \fint_{C} |u(z) - u(x)| \, dz \, dx \mid C \subset \mathbb{R}^N \text{ is a cube} \right\}.
\end{equation}
\end{definition}

The following corollary is a slight modification of Theorem \ref{thm:BMO inequality}:
\begin{corollary}
\label{cor:BMO inequality with eta}
 Let $C\subset\R^N$ be a cube, $u\in BMO(C)$. Define for $\sigma\in [0,\infty)$
\begin{equation}
S_\sigma:=\Set{x\in C}[|u(x)-u_C|>\sigma].
\end{equation}
Then,
\begin{equation}
\mathcal{L}^N\left(S_\sigma\right)
\leq \frac{1}{\eta}e^{-\alpha\sigma \eta^{-1}}\int_{C}|u(x)-u_{C}|dx,
\end{equation}
for every $\eta,\sigma$ such that
\begin{equation}
\eta\geq \|u\|_{BMO(C)}\quad \text{and}\quad \frac{\sigma}{\eta}\geq a+\alpha^{-1}.
\end{equation}
Here, $\alpha$ and $a$ are positive constants depending only on the dimension $N$.
\end{corollary}

\begin{proof}
Let $\alpha,a\in (0,\infty)$ be the positive numbers from Theorem \ref{thm:BMO inequality}. Define a function $g(t):=te^{-\alpha t}$. Then,
$g'(t)=e^{-\alpha t}-\alpha te^{-\alpha t}=(1-\alpha
t)e^{-\alpha t}<0,\,\forall t>\alpha^{-1}$. Therefore, the function $g$ is monotonically decreasing in the interval $(\alpha^{-1},\infty)$.
If $\eta\geq
\|u\|_{BMO(C)}$ and $\frac{\sigma}{\eta}\geq a+\frac{1}{\alpha}$,
then $\frac{\sigma}{\|u\|_{BMO(C)}}\geq \frac{\sigma}{\eta}\geq
a+\frac{1}{\alpha}>\frac{1}{\alpha}$, and thus
$g\left(\frac{\sigma}{\|u\|_{BMO(C)}}\right)\leq
g\left(\frac{\sigma}{\eta}\right)$.
Let us denote $I:=\int_{C}|u(x)-u_{C}|dx$. By Theorem \ref{thm:BMO inequality}, since $\frac{\sigma}{\|u\|_{BMO(C)}}\geq a$ and the parameter $A$ satisfies $A\leq 1$, we obtain
\begin{equation}
\label{eq:the BMO level set inequality with eta}
\mathcal{L}^N\left(S_\sigma\right)\leq 
 \frac{1}{\|u\|_{BMO(C)}}e^{-\alpha\sigma \|u\|_{BMO(C)}^{-1}}I=\frac{1}{\sigma}g\left(\frac{\sigma}{\|u\|_{BMO(C)}}\right)I
\leq\frac{1}{\sigma}g\left(\frac{\sigma}{\eta}\right)I
=\frac{1}{\eta}e^{-\alpha\sigma \eta^{-1}}I.
\end{equation}
\end{proof}

The following key lemma is a further generalization of John-Nirenberg inequality: 
\begin{lemma}
\label{lem:BMO inequality with eta refined}
 Let $C\subset\R^N$ be a cube, $u\in BMO(C)$. Define for $\sigma\in [0,\infty)$
\begin{equation}
S_\sigma:=\Set{x\in C}[|u(x)-u_C|>\sigma].
\end{equation}
Then,
\begin{equation}
\label{eq:JN-with indicator}
\mathcal{L}^N\left(S_\sigma\right)
\leq \frac{\beta}{\eta}e^{-\alpha\sigma \eta^{-1}}\int_{\Set{x\in C}[|u(x)-u_C|>\eta]}|u(x)-u_{C}|dx,
\end{equation}
for every $\eta,\sigma$ such that
\begin{equation}
\label{eq:additional assumption on eta and sigma}
\sigma \geq \eta\geq \vec k:=\|u\|_{BMO(C)}.
\end{equation}
Here, $\beta=\beta(N)$ and $\alpha=\alpha(N)$ are positive constants depending only on the dimension $N$.
\end{lemma}

\begin{proof}
Denote $g:=u\chi_{\{|u|>\eta\}}$. The idea of the proof is to use the function $g$, which is also in $BMO$, in place of $u$ in the John--Nirenberg inequality.

Without loss of generality, assume that $u_C=0$; otherwise, replace $u$ by $u-u_C$. Notice that 
\begin{equation}
|g_C|\leq \fint_{C}|g(x)|dx\leq \fint_{C}|u(x)|dx=\fint_{C}|u(x)-u_C|dx\leq \vec k,
\end{equation}
Let $a,\alpha$ be the dimensional constants from Corollary \ref{cor:BMO inequality with eta}. Assume that $\sigma>\vec k$; otherwise, if $\sigma=\vec k$, we obtain \eqref{eq:JN-with indicator} from Chebyshev's inequality with $\beta:=e^{\alpha}$. We get
\begin{multline}
\label{eq:Ssigma subsets of g-level set}
S_\sigma=\Set{x\in C}[|u(x)|>\sigma]=\Set{x\in C}[|u(x)\chi_{\{|u|>\eta\}}(x)|>\sigma]=\Set{x\in C}[|g(x)|>\sigma]
\\
\subset \Set{x\in C}[|g(x)-g_C|>\sigma-|g_C|]
\subset \Set{x\in C}[|g(x)-g_C|>\sigma-\vec k].
\end{multline}
We have $g\in BMO(C)$ because $g=u-u\chi_{\{|u|\leq \eta\}}$, $u\in BMO(C)$, and $u\chi_{\{|u|\leq \eta\}}\in L^\infty(C)$. Note that
\begin{equation}
\|g\|_{BMO(C)}\leq \|u\|_{BMO(C)}+\|u\chi_{|u|\leq \eta}\|_{BMO(C)}\leq \eta+2\eta = 3\eta.
\end{equation}
Assume first that
\begin{equation}
\label{eq:restrication1}
\frac{\sigma-\vec k}{3\eta}\geq a+\alpha^{-1}.
\end{equation}
Therefore, by \eqref{eq:Ssigma subsets of g-level set} and Corollary \ref{cor:BMO inequality with eta},
\begin{multline}
\label{eq:Leb of singma}
\mathcal{L}^N\left(S_\sigma\right)\leq \mathcal{L}^N\left(\Set{x\in C}[|g(x)-g_C|>\sigma-\vec k]\right)
\leq \frac{1}{3\eta}e^{-\alpha(\sigma-\vec k) (3\eta)^{-1}}\int_{C}|g(x)-g_{C}|dx
\\
\leq \frac{2}{3\eta}e^{-\alpha(\sigma-\vec k) (3\eta)^{-1}}\int_{C}|g(x)|dx
=\frac{2}{3\eta}e^{-\alpha(\sigma-\vec k) (3\eta)^{-1}}\int_{\Set{x\in C}[|u(x)|>\eta]}|u(x)|dx.
\end{multline}
Notice that
\begin{equation}
\label{eq:Leb of singma1}
e^{-\alpha(\sigma-\vec k) (3\eta)^{-1}}=e^{-\alpha\sigma(3\eta)^{-1}}e^{\alpha \vec k (3\eta)^{-1}}\leq e^{-\frac{\alpha}{3}\sigma\eta^{-1}}e^{\alpha/3}\quad (\eta\geq \vec k).
\end{equation}
Thus, \eqref{eq:Leb of singma} and \eqref{eq:Leb of singma1} give
\begin{equation}
\label{eq:JN-with indicator1}
\mathcal{L}^N\left(S_\sigma\right)
\leq \frac{2e^{\alpha/3}}{3\eta}e^{-\frac{\alpha}{3}\sigma \eta^{-1}}\int_{\Set{x\in C}[|u(x)|>\eta]}|u(x)|dx.
\end{equation}
Inequality \eqref{eq:JN-with indicator1} is the desired result \eqref{eq:JN-with indicator} under the additional restriction \eqref{eq:restrication1}.

In the case
\begin{equation}
\label{eq:case where sigma-k/eta less that dimensional constant}
\frac{\sigma-\vec k}{3\eta}< a+\alpha^{-1},
\end{equation}
we have, by Chebyshev's inequality,
\begin{multline}
\label{eq:ljljl}
\mathcal{L}^N\left(S_\sigma\right)
=\mathcal{L}^N\left(\Set{x\in C}[|u(x)|>\sigma]\right)=\mathcal{L}^N\left(\Set{x\in C}[|u(x)\chi_{S_\sigma}(x)|>\sigma]\right)
\\
\leq \frac{1}{\sigma}\int_{C}\chi_{S_\sigma}(x)|u(x)|dx
=
e^{\alpha(\sigma-\vec k) (3\eta)^{-1}}e^{-\alpha (\sigma-\vec k) (3\eta)^{-1}}\frac{1}{\sigma}\int_{\Set{x\in C}[|u(x)|>\sigma]}|u(x)|dx.
\end{multline}
Using \eqref{eq:case where sigma-k/eta less that dimensional constant} and \eqref{eq:Leb of singma1}, we have
\begin{equation}
\label{eq:ljljl1}
e^{\alpha(\sigma-\vec k) (3\eta)^{-1}}e^{-\alpha (\sigma-\vec k) (3\eta)^{-1}}\leq e^{\alpha a+1}e^{-\alpha (\sigma-\vec k) (3\eta)^{-1}}
\leq e^{\alpha a+1}\left(e^{-\frac{\alpha}{3}\sigma\eta^{-1}}e^{\alpha/3}\right)=\left(e^{\alpha (a+1/3)+1}\right)e^{-\frac{\alpha}{3}\sigma\eta^{-1}}.
\end{equation}
Therefore, from \eqref{eq:ljljl} and \eqref{eq:ljljl1} we get
\begin{multline}
\mathcal{L}^N\left(S_\sigma\right)
\leq \left(e^{\alpha (a+1/3)+1}\right)e^{-\frac{\alpha}{3}\sigma\eta^{-1}}\frac{1}{\sigma}\int_{\Set{x\in C}[|u(x)|>\sigma]}|u(x)|dx
\\
\leq \left(e^{\alpha (a+1/3)+1}\right)e^{-\frac{\alpha}{3}\sigma\eta^{-1}}\frac{1}{\eta}\int_{\Set{x\in C}[|u(x)|>\eta]}|u(x)|dx.
\end{multline}
Denote
\begin{equation}
\label{eq:constant}
\beta':=\frac{2}{3}e^{\alpha/3}+e^{\alpha (a+1/3)+1},\quad \alpha':=\alpha/3.
\end{equation}
Note that $\beta',\alpha'$ depend only on the dimension $N$. This completes the proof.
\end{proof}

\begin{corollary}
\label{thm:L1 norm of u times indicator} 
Let $Q\subset\R^N$ be a cube, $u\in BMO(Q)$ and $\vec k \in [0,\infty)$ be such that
\begin{equation}
\label{eq:uniform bound on BMO norms137jjkjkj}
\left\|u\right\|_{BMO(Q)}\leq \vec k.
\end{equation}
Then,
\begin{equation}
\label{eq:weak137gjdjj}
\left\|\chi_{\{|u-u_Q|>\vec k\}}(u-u_Q)\right\|_{L^{1}(Q)}\leq C\,\vec k \,\mathcal{L}^N\left({\Set{x\in
Q}[\left|u(x)-u_Q\right|>\vec k]}\right).
\end{equation}
Here, $C:=C(N)$ is a positive constant that depends only on
$N$.
\end{corollary}

\begin{proof}
Assume that $u_Q=0$ and $\vec k=1$. Then, using Fubini's theorem, for every $\gamma\in (0,1)$ we get

\begin{multline}
\label{eq:estimating of u>1nnn333331059}
\left\|u\chi_{\{|u|>1\}}\right\|_{L^{1}(Q)}
=\int_0^\infty\mathcal{L}^N\left(Q\cap \left\{\left|u\chi_{\{|u|>1\}}\right|>t\right\}\right)dt
\\
=\int_0^\infty\mathcal{L}^N\left(Q\cap \left\{\left|u\chi_{\{|u|>1\}}\right|>t\right\}\right)^{1-\gamma}
\mathcal{L}^N\left(Q\cap \left\{\left|u\chi_{\{|u|>1\}}\right|>t\right\}\right)^\gamma dt
\\
=\int_0^\infty\mathcal{L}^N\left(Q\cap \left\{|u|>\max\{t,1\}\right\}\right)^{1-\gamma}\mathcal{L}^N\left(Q\cap \left\{\left|u\chi_{\{|u|>1\}}\right|>t\right\}\right)^\gamma dt
\\
\leq
\mathcal{L}^N\left(Q\cap \{|u|>1\}\right)^{1-\gamma}
\int_0^\infty\mathcal{L}^N\left(Q\cap \left\{\left|u\chi_{\{|u|>1\}}\right|>t\right\}\right)^\gamma dt.
\end{multline}

Note that by Chebyshev's inequality
\begin{equation}
\label{eq:estimating of u>1nnn33333105988899999}
\int_0^1\mathcal{L}^N\left(Q\cap \left\{\left|u\chi_{\{|u|>1\}}\right|>t\right\}\right)^\gamma dt
=
\mathcal{L}^N\left(Q\cap \left\{\left|u\chi_{\{|u|>1\}}\right|>1\right\}\right)^\gamma
\leq \left\|u\chi_{\{|u|>1\}}\right\|^\gamma_{L^{1}(Q)}.
\end{equation}

So far we didn't use the $BMO$-assumption. Now, by Lemma \ref{lem:BMO inequality with eta refined}, we get
\begin{multline}
\label{eq:step where using BMO assumption}
\int_1^\infty\left(\mathcal{L}^N\left(Q\cap \{|u|>t\}\right)\right)^\gamma dt
\leq \int_1^\infty\left(\beta e^{-\alpha t }\int_{Q\cap \{|u|>1\}}|u|dx\right)^\gamma dt
\\
=\left(\beta \int_{Q\cap \{|u|>1\}}|u|dx\right)^\gamma
\int_1^\infty  e^{-\alpha \gamma t } dt
=\left(\beta \left\|u\chi_{\{|u|>1\}}\right\|_{L^{1}(Q)}\right)^\gamma
\frac{e^{-\alpha\gamma}}{\alpha\gamma}.
\end{multline}

Thus, by \er{eq:estimating of u>1nnn33333105988899999} and \er{eq:step where using BMO assumption} together we get
\begin{multline}
\label{eq:estimating of u>1nnn333331059kjhjjh}
\int_0^\infty\mathcal{L}^N\left(Q\cap \left\{\left|u\chi_{\{|u|>1\}}\right|>t\right\}\right)^\gamma dt
=\int_0^1\mathcal{L}^N\left(Q\cap \left\{\left|u\chi_{\{|u|>1\}}\right|>t\right\}\right)^\gamma dt
\\
+\int_1^\infty\mathcal{L}^N\left(Q\cap \left\{\left|u\chi_{\{|u|>1\}}\right|>t\right\}\right)^\gamma dt
\leq \left\|u\chi_{\{|u|>1\}}\right\|_{L^{1}(Q)}^\gamma
\left(1+\beta^\gamma\frac{e^{-\alpha\gamma}}{\alpha\gamma}\right).
\end{multline}

Therefore, using \eqref{eq:estimating of u>1nnn333331059} and \eqref{eq:estimating of u>1nnn333331059kjhjjh}, we obtain
\begin{equation}
\label{eq:step where using BMO assumption1}
\left\|u\chi_{\{|u|>1\}}\right\|_{L^{1}(Q)}
\leq \mathcal{L}^N\left(Q\cap \{|u|>1\}\right)^{1-\gamma}
\left\|u\chi_{\{|u|>1\}}\right\|_{L^{1}(Q)}^\gamma
\left(1+\beta^\gamma\frac{e^{-\alpha\gamma}}{\alpha\gamma}\right).
\end{equation}

From \eqref{eq:step where using BMO assumption1} we finally get
\begin{equation}
\label{eq:estimating of u>1nnn33333100}
\left\|u\chi_{\{|u|>1\}}\right\|_{L^{1}(Q)}
\leq
\left(1+\beta^\gamma\frac{e^{-\alpha\gamma}}{\alpha\gamma}\right)^{\frac{1}{1-\gamma}}
\,\mathcal{L}^N\left(Q\cap \{|u|>1\}\right).
\end{equation}
\end{proof}

The following corollary proves Theorem \ref{thm:relation between Lebesgue measures of level sets of BMO}:
\begin{corollary}[Refined John--Nirenberg]
\label{cor:relation between Lebesgue measures of level sets of BMO}
Let $Q\subset\R^N$ be a cube and $u\in BMO(Q)$. 
Define for $\sigma\in [0,\infty)$
\begin{equation}
S_\sigma:=\Set{x\in Q}[|u(x)-u_Q|>\sigma].
\end{equation}
Then,
\begin{equation}
\mathcal{L}^N\left(S_\sigma\right)
\leq \beta e^{-\alpha\sigma \eta^{-1}}\mathcal{L}^N\left(S_\eta\right)
\end{equation}
for every $\eta,\sigma>0$ such that
\begin{equation}
\label{eq:choice of sigma and eta}
 \|u\|_{BMO(Q)}\leq\eta\leq\sigma.
\end{equation}
Here, $\beta=\beta(N)$ and $\alpha=\alpha(N)$ are positive constants depending only on the dimension $N$.
\end{corollary}

\begin{proof}
Let $\sigma,\eta$ be as in \eqref{eq:choice of sigma and eta}. By Lemma \ref{lem:BMO inequality with eta refined}, we get
\begin{multline}
\mathcal{L}^N\left(S_\sigma\right)
\leq \frac{\beta}{\eta}e^{-\alpha\sigma \eta^{-1}}
\int_{\Set{x\in Q}[|u(x)-u_Q|>\eta]}|u(x)-u_Q|\,dx
\\
= \frac{\beta}{\eta}e^{-\alpha\sigma \eta^{-1}}
\left\|\chi_{\{|u-u_Q|>\eta\}}(u-u_Q)\right\|_{L^{1}(Q)}.
\end{multline}
By Corollary \ref{thm:L1 norm of u times indicator},
\begin{equation}
\mathcal{L}^N\left(S_\sigma\right)
\leq \frac{\beta}{\eta}e^{-\alpha\sigma \eta^{-1}}
\Big[C(N)\eta\,\mathcal{L}^N\left(\Set{x\in Q}[|u(x)-u_Q|>\eta]\right)\Big]
=\beta C(N)e^{-\alpha\sigma \eta^{-1}}
\,\mathcal{L}^N\left(S_\eta\right).
\end{equation}
\end{proof}
\begin{proof}[Proof of Theorem \ref{thm:relation between Lebesgue measures of level sets of BMO}]
Let $u\in BMO(Q)$. Without loss of generality we may assume that $u_Q=0$. Denote $\vec k:=\|u\|_{BMO(Q)}$. Given $\eta\geq 0$
let us define a function $\psi:\R\to \R$ by
\begin{equation}
\psi(t)=0 \quad \text{whenever } |t|\leq \eta,\quad \text{and}\quad 
\psi(t)=|t|-\eta \quad \text{for } |t|\geq \eta.
\end{equation}
The function $\psi(u)$ satisfies $\|\psi(u)\|_{BMO(Q)}\leq \|u\|_{BMO(Q)}=\vec k$ because $\psi$ is Lipschitz with Lipschitz constant $1$. Moreover, $0\leq\psi(u)\leq |u|$ and thus,
\begin{equation}
\label{eq:average is bounded by k}
|\psi(u)_Q|\leq |u|_Q=|u-u_Q|_Q\leq\vec k.
\end{equation}
For a function $g$ and number $\delta$ we define:
\begin{equation}
S_{\delta}(g):=\Set{x\in Q}[|g(x)-g_Q|>\delta].
\end{equation}
Then, given $\tau\geq 2\vec k$, by Corollary \ref{cor:relation between Lebesgue measures of level sets of BMO} we have
\begin{equation}\label{jkhhhhhj}
\mathcal{L}^N\left(S_\tau(\psi(u))\right)
\leq \beta e^{-\alpha\tau {(2\vec k)}^{-1}}\mathcal{L}^N\left(S_{2\vec k}(\psi(u))\right)\,.
\end{equation}
The non-negativity of $\psi(u)$ and \eqref{eq:average is bounded by k} give
\begin{multline}
\label{eq:set-bound for S2k}
S_{2\vec k}(\psi(u)):=\Set{x\in Q}[|(\psi(u))(x)-(\psi(u))_Q|>2\vec k]\subset \Set{x\in Q}[|(\psi(u))(x)|+|(\psi(u))_Q|>2\vec k]\\
\subset \Set{x\in Q}[(\psi(u))(x)>\vec k]= \Set{x\in Q}[|u(x)|-\eta>\vec k]=\Set{x\in Q}[|u(x)-u_Q|>\eta+\vec k]\\=S_{\eta+\vec k}(u)\subset S_{\eta}(u)\,.
\end{multline}
Thus, by \er{jkhhhhhj} and \eqref{eq:set-bound for S2k}
\begin{equation}\label{jkhhhhhjr}
\mathcal{L}^N\left(S_\tau(\psi(u))\right)
\leq \beta e^{-\alpha\tau {(2\vec k)}^{-1}}\mathcal{L}^N\left(S_{\eta}(u)\right)\,.
\end{equation}
In addition,
\begin{multline}
S_{\tau}(\psi(u)):=\Set{x\in Q}[|(\psi(u))(x)-(\psi(u))_Q|>\tau]\supset \Set{x\in Q}[|(\psi(u))(x)|-|(\psi(u))_Q|>\tau]\\
\supset \Set{x\in Q}[(\psi(u))(x)>\tau+\vec k]= \Set{x\in Q}[|u(x)|-\eta>\tau+\vec k]\\=\Set{x\in Q}[|u(x)-u_Q|>\tau+\eta+\vec k]=S_{\tau+\eta+\vec k}(u).
\end{multline}
Thus,
\begin{equation}\label{jkhhhhhjrrf}
\mathcal{L}^N\left(S_{\tau+\eta+\vec k}(u)\right)
\leq \beta e^{-\alpha\tau {(2\vec k)}^{-1}}\mathcal{L}^N\left(S_{\eta}(u)\right)\quad\quad\forall \tau\geq 2\vec k\,.
\end{equation}
Then, by substituting $\sigma=\tau+\eta+\vec k$, it follows that 
\begin{equation}\label{jkhhhhhjrkjkjh}
\mathcal{L}^N\left(S_{\sigma}(u)\right)
\leq \beta e^{\alpha/2} e^{-\alpha(\sigma-\eta) {(2\vec k)}^{-1}}\mathcal{L}^N\left(S_{\eta}(u)\right)\quad\quad\forall \sigma\geq \eta+3\vec k\,.
\end{equation}
On the other hand if $\eta\leq\sigma\leq\eta+3\vec k$ then 
\begin{multline}\label{jkhhhhhjrkjkjhhhhg}
\mathcal{L}^N\left(S_{\sigma}(u)\right)
\leq\mathcal{L}^N\left(S_{\eta}(u)\right)
=e^{\alpha(\sigma-\eta){(2\vec k)}^{-1}}e^{-\alpha(\sigma-\eta){(2\vec k)}^{-1}}\mathcal{L}^N\left(S_{\eta}(u)\right)
\\
\leq e^{3\alpha/2}e^{-\alpha(\sigma-\eta){(2\vec k)}^{-1}}\mathcal{L}^N\left(S_{\eta}(u)\right)\quad\quad\forall \eta\leq\sigma\leq\eta+3\vec k\,.
\end{multline}
Thus taking $\gamma:=\max\Big\{\beta e^{\alpha/2},e^{3\alpha/2}\Big\}$
we finally deduce
\begin{equation}\label{jkhhhhhjrkjkjhoiio}
\mathcal{L}^N\left(S_{\sigma}(u)\right)
\leq \gamma\,e^{-\frac{\alpha}{2}\frac{\sigma-\eta}{\vec k}}\mathcal{L}^N\left(S_{\eta}(u)\right)\quad\quad\forall 0\leq\eta\leq \sigma\,.
\end{equation}
Thus, \er{jkhhhhhjrkjkjhoiiofg} is proven with constant $\frac{\alpha}{2}$ in the place of $\alpha$.
\end{proof}
\begin{lemma}
\label{lem:convergence to zero of averages of bounded functions in Lqw}
Let $q\in (0,\infty)$, 
$u\in L^\infty(\R^N)$, and let $\{Q_n\}_{n\in\N}$ be a sequence of cubes monotonically increasing to $\R^N$. 
Assume that there exists $\delta\in (0,\infty)$ such that 
\begin{equation}
\label{eq:finiteness of Lqw on finite interval}
\sup_{0 < s < \delta} s\,\mathcal{L}^N\!\left(\Set{x \in
\R^N}[|u(x)|^q > s]\right)<\infty.
\end{equation}
Then,
\begin{equation}
\label{eq:uniform limit of averages5}
\lim_{n \to \infty} |u|_{Q_n} = 0.
\end{equation}
\end{lemma}

\begin{proof}
By a change of variable in the supremum, the finiteness of \eqref{eq:finiteness of Lqw on finite interval} is equivalent to the finiteness of the same quantity when $\delta$ is replaced by any positive number. 
Indeed, if $a,b\in (0,\infty)$ with $a\leq b$, then
\begin{equation}
\sup_{0 < s < a} \Big[s\,\mathcal{L}^N\big(\Set{x \in
\R^N}[|u(x)|^q > s]\big)\Big]\leq \sup_{0 < s < b} \Big[s\,\mathcal{L}^N\big(\Set{x \in
\R^N}[|u(x)|^q > s]\big)\Big]
\end{equation}
and
\begin{equation}
\sup_{0 < s < b} \Big[s\,\mathcal{L}^N\big(\Set{x \in
\R^N}[|u(x)|^q > s]\big)\Big]\leq \frac{b}{a}\sup_{0 < t < a} \Big[t\,\mathcal{L}^N\big(\Set{x \in
\R^N}[|u(x)|^q > t]\big)\Big].
\end{equation}

Let $\alpha := \|u\|^q_{L^\infty(\R^N)}$ and fix $\xi$ such that $0 < \xi < \min\{1/q,1\}$.  
By Fubini's theorem and a change of variables, we obtain
\begin{multline}
\label{eq:infinitesimality of the average of u on the cube in case p>15}
|u|_{Q_n} 
= \frac{1}{\mathcal{L}^N(Q_n)} \int_{Q_n} |u(x)| \, dx 
= \frac{1}{\mathcal{L}^N(Q_n)} \int_0^{\infty} \mathcal{L}^N\!\big(\Set{x \in
Q_n}[|u(x)|^q > t^q]\big) \, dt
\\
= \frac{1}{q \mathcal{L}^N(Q_n)} \int_0^{\infty} s^{\frac{1}{q}-1} \mathcal{L}^N\!\big(\Set{x \in
Q_n}[|u(x)|^q > s]\big) \, ds
\\
= \frac{1}{q \mathcal{L}^N(Q_n)} \int_0^{\alpha} s^{\frac{1}{q}-1} \mathcal{L}^N\!\big(\Set{x \in
Q_n}[|u(x)|^q > s]\big) \, ds
\\
= \frac{1}{q \mathcal{L}^N(Q_n)} \int_0^{\alpha} s^{\frac{1}{q}-1} \big(\mathcal{L}^N(\Set{x \in
Q_n}[|u(x)|^q > s])\big)^{1-\xi} \big(\mathcal{L}^N(\Set{x \in
Q_n}[|u(x)|^q > s])\big)^{\xi} \, ds
\\
\leq \frac{1}{q \mathcal{L}^N(Q_n)} \int_0^{\alpha} s^{\frac{1}{q}-1} \big(\mathcal{L}^N(Q_n)\big)^{1-\xi} \big(\mathcal{L}^N(\Set{x \in
Q_n}[|u(x)|^q > s])\big)^{\xi} \, ds
\\
=\frac{1}{q \mathcal{L}^N(Q_n)^{\xi}} \int_0^{\alpha} s^{\frac{1}{q}-1-\xi} \big(s\,\mathcal{L}^N(\Set{x \in
Q_n}[|u(x)|^q > s])\big)^{\xi} \, ds
\\
\leq \sup_{0 < s < \alpha} \big(s\,\mathcal{L}^N(\Set{x \in
\R^N}[|u(x)|^q > s])\big)^{\xi}  \frac{\alpha^{\frac{1}{q}-\xi}}{q \mathcal{L}^N(Q_n)^{\xi} \left( \frac{1}{q} - \xi \right)} < \infty.
\end{multline}  
Therefore, \eqref{eq:uniform limit of averages5} follows by taking $n \to \infty$ in \eqref{eq:infinitesimality of the average of u on the cube in case p>15}.
\end{proof}

\begin{lemma}
\label{lem:BMO inequality with eta refined,global}
Let $u\in BMO(\R^N)$. Define for $\sigma\in [0,\infty)$
\begin{equation}
S_\sigma:=\Set{x\in \R^N}[|u(x)|>\sigma].
\end{equation}
Then,
\begin{equation}
\label{eq:refined version of JN - global}
\mathcal{L}^N\left(S_\sigma\right)
\leq \frac{\beta}{\eta}e^{-\alpha \,\sigma \eta^{-1}}\int_{\Set{|u|>\eta}}|u(x)|dx,
\end{equation}
for every $\eta,\sigma>0$ such that
\begin{equation}
\label{eq:additional assumption on eta and sigma,global}
\sigma \geq \eta\geq \vec k:=\|u\|_{BMO(\R^N)}.
\end{equation}
Here, $\beta=\beta(N)$ and $\alpha=\alpha(N)$ are positive constants depending only on the dimension $N$.
\end{lemma}

\begin{proof}
Let $\{Q_n\}$ be a sequence of cubes monotonically increasing to $\R^N$. By Lemma \ref{lem:BMO inequality with eta refined}, we get for every $n\in\N$ and $\sigma \geq \eta\geq \vec k$
\begin{equation}
\label{eq:jn-local-cubes3}
\mathcal{L}^N\left(\Set{x\in Q_n}[|u(x)-u_{Q_n}|>\sigma]\right)
\leq \frac{\beta}{\eta}e^{-\alpha\sigma \eta^{-1}}\int_{\Set{x\in Q_n}[|u(x)-u_{Q_n}|>\eta]}|u(x)-u_{Q_n}|dx.
\end{equation}

\textbf{Step 1:}
\\
We first prove \eqref{eq:refined version of JN - global} under the assumptions that $u\in L^\infty(\R^N)$ and there exists $\delta\in (0,\infty)$ such that 
\begin{equation}
\label{eq:lqw-finite}
\sup_{0 < s < \delta} s\,\mathcal{L}^N\!\left(\Set{x \in
\R^N}[|u(x)| > s]\right)<\infty.
\end{equation}
By Lemma \ref{lem:convergence to zero of averages of bounded functions in Lqw}, we get 
\begin{equation}
\label{eq:averages-converge-zero}
\lim_{n\to \infty}|u_{Q_n}|=0.
\end{equation}

By Fatou's lemma and \eqref{eq:averages-converge-zero}, we get
\begin{multline}
\label{eq:fatou-level-sets}
\liminf_{n\to \infty}\mathcal{L}^N\left(\Set{x\in Q_n}[|u(x)-u_{Q_n}|>\sigma]\right)
\geq \int_{\R^N}\liminf_{n\to \infty}\chi_{\Set{x\in Q_n}[|u(x)-u_{Q_n}|>\sigma]}(x)dx
\\
\geq  \int_{\R^N}\chi_{\Set{x\in \R^N}[|u(x)|>\sigma]}(x)dx
= \mathcal{L}^N\left(\Set{x\in \R^N}[|u(x)|>\sigma]\right).
\end{multline}

Assume that 
\begin{equation}
\label{eq:integrable-tail}
\int_{\Set{|u|>\eta}}|u(x)|dx<\infty;
\end{equation}
otherwise, inequality \eqref{eq:refined version of JN - global} trivially holds.

Let $\e>0$. Since \eqref{eq:averages-converge-zero} holds, there exists $n_0\in\N$ such that $|u_{Q_n}|<\min\{\e,\eta\}$ for every $n\geq n_0$. Therefore, for every $n\geq n_0$, since $|u_{Q_n}|<\e$
\begin{equation}
\Set{x\in Q_n}[|u(x)-u_{Q_n}|>\eta+\e]\subset \{|u|>\eta\}.
\end{equation}

Thus, for every $n\geq n_0$,
\begin{equation}
\label{eq:dct-preparation}
\int_{\Set{x\in Q_n}[|u(x)-u_{Q_n}|>\eta +\e]}|u(x)-u_{Q_n}|dx
\leq \int_{\{|u|>\eta\}}|u(x)-u_{Q_n}|dx.
\end{equation}

For almost every $x\in \{|u|>\eta\}$, we have 
\begin{equation}
\label{eq:pointwise-convergence}
\lim_{n\to \infty}|u(x)-u_{Q_n}|=|u(x)|.
\end{equation}
and 
\begin{equation}
\label{eq:dct-dominating-bound}
|u(x)-u_{Q_n}|\leq |u(x)|+|u_{Q_n}|<|u(x)|+\eta<2|u(x)|.
\end{equation}

Therefore, taking into account \eqref{eq:integrable-tail}, \eqref{eq:pointwise-convergence}, \eqref{eq:dct-dominating-bound}, and the dominated convergence theorem, we obtain by taking the lower limit as $n\to \infty$ on both sides of \eqref{eq:dct-preparation}
\begin{multline}
\label{eq:liminf-estimate-eps}
\liminf_{n\to \infty}\int_{\Set{x\in Q_n}[|u(x)-u_{Q_n}|>\eta+\e]}|u(x)-u_{Q_n}|dx
\\
\leq \lim_{n\to \infty}\int_{\{|u|>\eta\}}|u(x)-u_{Q_n}|dx
= \int_{\{|u|>\eta\}}|u(x)|dx.
\end{multline}
Assume that $\sigma>\eta$; if $\sigma=\eta$, then \eqref{eq:refined version of JN - global} follows from Chebyshev's inequality with $\beta:=e^{\alpha}$. Replacing $\eta$ with $\eta+\e$, provided $\eta+\e<\sigma$, in \eqref{eq:jn-local-cubes3}, we obtain
\begin{equation}
\label{eq:jn-local-cubes1}
\mathcal{L}^N\left(\Set{x\in Q_n}[|u(x)-u_{Q_n}|>\sigma]\right)
\leq \frac{\beta}{\eta+\e}e^{-\alpha\sigma (\eta+\e)^{-1}}\int_{\Set{x\in Q_n}[|u(x)-u_{Q_n}|>\eta+\e]}|u(x)-u_{Q_n}|dx.
\end{equation}
Taking the lower limit as $n\to \infty$ on both sides of \eqref{eq:jn-local-cubes1} and taking into account \eqref{eq:fatou-level-sets} and \eqref{eq:liminf-estimate-eps}, we obtain
\begin{equation}
\label{eq:jn-local-cubes}
\mathcal{L}^N\left(\Set{x\in \R^N}[|u(x)|>\sigma]\right)
\leq \frac{\beta}{\eta+\e}e^{-\alpha\sigma (\eta+\e)^{-1}}\int_{\Set{x\in \R^N}[|u(x)|>\eta]}|u(x)|dx.
\end{equation}
Taking the limit as $\e\to 0^+$ in \eqref{eq:jn-local-cubes}, we get \eqref{eq:refined version of JN - global}.

\textbf{Step 2:}
\\
We now generalize \eqref{eq:refined version of JN - global} to any $u\in BMO(\R^N)$, while assuming that \eqref{eq:lqw-finite} still holds. Let $\{u_h\}_{h \in \N}$ be the truncations of $u$ defined by 
\begin{equation}
\label{eq:truncation-definition}
u_h(x) = \max(-h,\, \min(u(x),h)), \quad h \in \N.
\end{equation}

Then $u_h \in L^\infty(\R^N)$, $|u_h(x)| \uparrow |u(x)|$ for a.e $x\in\R^N$, and \eqref{eq:lqw-finite} holds for each $u_h$. Hence, for every $h\in\N$ and every $\sigma \geq \eta\geq \vec k\geq \|u_h\|_{BMO(\R^N)}$, we get 
\begin{equation}
\label{eq:jn-truncated}
\mathcal{L}^N\left(\Set{|u_h|>\sigma}\right)
\leq \frac{\beta}{\eta}e^{-\alpha\sigma \eta^{-1}}\int_{\Set{|u_h|>\eta}}|u_h(x)|dx
\leq \frac{\beta}{\eta}e^{-\alpha\sigma \eta^{-1}}\int_{\Set{|u|>\eta}}|u(x)|dx.
\end{equation}

Therefore, from \eqref{eq:jn-truncated}, we get
\begin{equation}
\label{eq:jn-general-lqw}
\mathcal{L}^N\left(\Set{|u|>\sigma}\right)
=\mathcal{L}^N\left(\bigcup_{h\in\N}\Set{|u_h|>\sigma}\right)
=\lim_{h\to \infty}\mathcal{L}^N\left(\Set{|u_h|>\sigma}\right)
\leq \frac{\beta}{\eta}e^{-\alpha\sigma \eta^{-1}}\int_{\Set{|u|>\eta}}|u(x)|dx
\end{equation}
for every $u\in BMO(\R^N)$ satisfying \eqref{eq:lqw-finite}. 

\textbf{Step 3:}
\\
We prove \eqref{eq:refined version of JN - global}. First of all note that as before, we may assume \eqref{eq:integrable-tail}, since otherwise \eqref{eq:refined version of JN - global} is trivial. In the case where \eqref{eq:lqw-finite} does not hold, let us define a function $f_\eta:\R\to \R$ by
\begin{equation}
f_\eta(t)=0 \quad \text{whenever } t\leq \eta,\quad \text{and}\quad 
f_\eta(t)=t-\eta \quad \text{for } t\geq \eta.
\end{equation}

The function $f_\eta(u)$ satisfies $\|f_\eta(u)\|_{BMO(\R^N)}\leq \|u\|_{BMO(\R^N)}$ because $f_\eta$ is Lipschitz with Lipschitz constant $1$. Notice that 
$\Set{|f_\eta(u)| > s}\subset \Set{|u| > \eta}$ for every $s\in (0,\infty)$. Therefore,
\begin{multline}
\label{eq:f-lqw-finite}
\sup_{0 < s < \eta} s\,\mathcal{L}^N\!\left(\Set{|f_\eta(u)| > s}\right)
\leq\sup_{0 < s < \eta} s\,\mathcal{L}^N\!\left(\Set{|u| > \eta}\right)
\\
=\eta\mathcal{L}^N\!\left(\Set{|u| > \eta}\right)
=\eta\mathcal{L}^N\!\left(\Set{|u\chi_{\{|u|>\eta\}}| > \eta}\right)
\leq \int_{\Set{|u|>\eta}}|u(x)|dx<\infty.
\end{multline}

Note that since $f_\eta$ is Lipschitz with Lipschitz constant $1$ and satisfies $f_\eta(0)=0$, we have $|f_\eta(u)|\leq |u|$. Therefore, by the previous case,
\begin{equation}
\label{eq:jn-applied-to-f}
\mathcal{L}^N\left(\Set{|f_\eta(u)|>\sigma}\right)
\leq \frac{\beta}{\eta}e^{-\alpha\sigma \eta^{-1}}\int_{\Set{|f_\eta(u)|>\eta}}|f_\eta(u(x))|dx
\leq \frac{\beta}{\eta}e^{-\alpha\sigma \eta^{-1}}\int_{\Set{|u|>\eta}}|u(x)|dx.
\end{equation}

Inequality \eqref{eq:jn-applied-to-f} holds for every $u\in BMO(\R^N)$ and every $\sigma\geq \eta\geq \|u\|_{BMO(\R^N)}$. 

Fix $u\in BMO(\R^N)$ and $\eta\geq \|u\|_{BMO(\R^N)}$. Consider first the case $\sigma>3\eta$. Note that 
\begin{equation}
\sigma - \eta > 2\eta>\eta\geq \|u\|_{BMO(\R^N)}\geq \||u|\|_{BMO(\R^N)}=\||u|+\eta\|_{BMO(\R^N)}.
\end{equation}
From \eqref{eq:jn-applied-to-f} 
\begin{multline}
\label{eq:shifted-jn}
\mathcal{L}^N\left(\Set{|u|>\sigma}\right)=\mathcal{L}^N\left(\Set{(|u|+\eta)-2\eta>\sigma-\eta}\right)=\mathcal{L}^N\left(\Set{|f_{2\eta}(|u|+\eta)|>\sigma-\eta}\right)
\\ \leq \frac{\beta}{2\eta}e^{-\alpha(\sigma-\eta) (2\eta)^{-1}}\int_{\Set{|u|+\eta>2\eta}}(|u(x)|+\eta)dx
=\frac{\beta}{2\eta}e^{\frac{\alpha}{2}}e^{-\alpha\sigma (2\eta)^{-1}}\int_{\Set{|u|>\eta}}(|u(x)|+\eta)dx
\\
\leq e^{\frac{\alpha}{2}}\frac{\beta}{\eta}e^{-\alpha\sigma (2\eta)^{-1}}\int_{\Set{|u|>\eta}}|u(x)|dx.
\end{multline}
Now, in the case $\eta\leq \sigma\leq 3\eta$, and by Chebyshev's inequality,
\begin{multline}
\label{eq:small-sigma-case}
\mathcal{L}^N\left(\Set{|u|>\sigma}\right)
\leq \frac{1}{\sigma}\int_{\Set{|u|>\sigma}}|u(x)|dx
\\
=\frac{e^{\alpha\sigma (2\eta)^{-1}}e^{-\alpha\sigma (2\eta)^{-1}}}{\sigma}\int_{\Set{|u|>\sigma}}|u(x)|dx
\leq \left(e^{\frac{3\alpha}{2}}\right)\frac{e^{-\alpha\sigma (2\eta)^{-1}}}{\sigma}\int_{\Set{|u|>\sigma}}|u(x)|dx
\\
\leq \left(e^{\frac{3\alpha}{2}}\right)\frac{e^{-\alpha\sigma (2\eta)^{-1}}}{\eta}\int_{\Set{|u|>\eta}}|u(x)|dx.
\end{multline}

Finally, \eqref{eq:shifted-jn} and \eqref{eq:small-sigma-case} together complete the proof, with $\frac{\alpha}{2}$ instead of $\alpha$ and $\max\{e^{\frac{\alpha}{2}}\beta, e^{\frac{3\alpha}{2}}\}$ instead of $\beta$.
\end{proof}

\begin{corollary}
\label{thm:L1 norm of u times indicator,global}
Let $u\in BMO(\R^N)$ and $\vec k \in (0,\infty)$ be such that
\begin{equation}
\left\|u\right\|_{BMO(\R^N)}\leq \vec k.
\end{equation}
Then,
\begin{equation}
\left\|\chi_{\{|u|>\vec k\}}u\right\|_{L^{1}(\R^N)}\leq C(N)\vec k \,\mathcal{L}^N\left({\Set{x\in
\R^N}[\left|u(x)\right|>\vec k]}\right).
\end{equation}
Here, $C(N)$ is a positive constant that depends only on
$N$.
\end{corollary}

\begin{proof}
The proof is the same as the proof of Corollary \ref{thm:L1 norm of u times indicator}, replacing the cube $Q$ by $\R^N$ and using Lemma \ref{lem:BMO inequality with eta refined,global} in place of Lemma \ref{lem:BMO inequality with eta refined} in step \eqref{eq:step where using BMO assumption}.
\end{proof}

The following corollary proves Theorem \ref{thm:relation between Lebesgue measures of level sets of BMO,global}:
\begin{corollary}[Refined John-Nirenberg -- global version]
\label{cor:exponential growth for level sets of BMO}
Let $u\in BMO(\R^N)$ . 
Define for $\sigma\in [0,\infty)$
\begin{equation}
S_\sigma:=\Set{x\in \R^N}[|u(x)|>\sigma].
\end{equation}
Then, 
\begin{equation}
\mathcal{L}^N\left(S_\sigma\right)
\leq \beta e^{-\alpha\sigma \eta^{-1}}\mathcal{L}^N\left(S_\eta\right)
\end{equation}
for every $\eta,\sigma>0$ such that
\begin{equation}
\label{eq:choice of sigma and eta1}
 \|u\|_{BMO(\R^N)}\leq\eta\leq\sigma.
\end{equation}
Here, $\beta=\beta(N)$ and $\alpha=\alpha(N)$ are positive constants depending only on the dimension $N$.
\end{corollary}

\begin{proof}
Let $\sigma,\eta$ as in \eqref{eq:choice of sigma and eta1}. By Lemma \ref{lem:BMO inequality with eta refined,global}, we get
\begin{equation}
\mathcal{L}^N\left(S_\sigma\right)
\leq \frac{\beta}{\eta}e^{-\alpha\sigma \eta^{-1}}\int_{\Set{|u|>\eta}}|u(x)|dx= \frac{\beta}{\eta}e^{-\alpha\sigma \eta^{-1}}\left\|\chi_{\{|u|>\eta\}}u\right\|_{L^{1}(\R^N)}.
\end{equation}
By Corollary \ref{thm:L1 norm of u times indicator,global}
\begin{equation}
\mathcal{L}^N\left(S_\sigma\right)
\leq \frac{\beta}{\eta}e^{-\alpha\sigma \eta^{-1}}\Big[C(N)\eta \,\mathcal{L}^N\left({\Set{x\in
\R^N}[\left|u(x)\right|>\eta]}\right)\Big]
=\beta C(N)e^{-\alpha\sigma \eta^{-1}}\,\mathcal{L}^N\left(S_\eta\right).
\end{equation}
\end{proof}

\begin{proof}[Proof of Theorem \ref{thm:relation between Lebesgue measures of level sets of BMO,global}]
Let $u\in BMO(\R^N)$. Denote $\vec k:=\|u\|_{BMO(\R^N)}$. Given $\eta\geq 0$, as before,
let us define a function $\psi:\R\to \R$ by
\begin{equation}
\psi(t)=0 \quad \text{whenever } |t|\leq \eta,\quad \text{and}\quad 
\psi(t)=|t|-\eta \quad \text{for } |t|\geq \eta.
\end{equation}
The function $\psi(u)$ satisfies $\|\psi(u)\|_{BMO(\R^N)}\leq \|u\|_{BMO(\R^N)}=\vec k$ because $\psi$ is Lipschitz with Lipschitz constant $1$. Moreover, $0\leq\psi(u)\leq |u|$. Define for a function $g$ and $\delta>0$ 
\begin{equation}
S_{\delta}(g):=\Set{x\in \R^N}[|g(x)|> \delta].
\end{equation}
Then, given $\tau\geq \vec k$, by Corollary \ref{cor:exponential growth for level sets of BMO} we have
\begin{equation}\label{jkhhhhhj333}
\mathcal{L}^N\left(S_\tau(\psi(u))\right)
\leq \beta e^{-\alpha\tau {\vec k}^{-1}}\mathcal{L}^N\left(S_{\vec k}(\psi(u))\right)\,.
\end{equation}
Note that
\begin{multline}
\label{eq:identity for Stau}
S_{\vec k}(\psi(u)):=\Set{x\in \R^N}[|(\psi(u))(x)|>\vec k]=
\Set{x\in \R^N}[(\psi(u))(x)>\vec k]=\\ \Set{x\in \R^N}[|u(x)|-\eta>\vec k]=\Set{x\in \R^N}[|u(x)|>\eta+\vec k]=S_{\eta+\vec k}(u)\subset S_{\eta}(u)\,.
\end{multline}
Thus, by \er{jkhhhhhj333} and \eqref{eq:identity for Stau}
\begin{equation}\label{jkhhhhhjr333}
\mathcal{L}^N\left(S_\tau(\psi(u))\right)
\leq \beta e^{-\alpha\tau {\vec k}^{-1}}\mathcal{L}^N\left(S_{\eta}(u)\right)\,.
\end{equation}
Repeating the computation in \eqref{eq:identity for Stau} with $\vec k$ replaced by $\eta$ yields
\begin{equation}
S_{\tau}(\psi(u))=S_{\tau+\eta}(u).
\end{equation}
Thus,
\begin{equation}\label{jkhhhhhjrrf333}
\mathcal{L}^N\left(S_{\tau+\eta}(u)\right)
\leq \beta e^{-\alpha\tau {\vec k}^{-1}}\mathcal{L}^N\left(S_{\eta}(u)\right)\quad\quad\forall \tau\geq \vec k\,.
\end{equation}
Substituting $\sigma=\tau+\eta$ gives
\begin{equation}\label{jkhhhhhjrkjkjh333}
\mathcal{L}^N\left(S_{\sigma}(u)\right)
\leq \beta e^{-\alpha(\sigma-\eta) {\vec k}^{-1}}\mathcal{L}^N\left(S_{\eta}(u)\right)\quad\quad\forall \sigma\geq \eta+\vec k\,.
\end{equation}
If $\eta\leq\sigma\leq\eta+\vec k$, then 
\begin{equation}\label{jkhhhhhjrkjkjhhhhg333}
\mathcal{L}^N\left(S_{\sigma}(u)\right)
\leq\mathcal{L}^N\left(S_{\eta}(u)\right)\leq e^{\alpha}e^{-\alpha(\sigma-\eta){\vec k}^{-1}}\mathcal{L}^N\left(S_{\eta}(u)\right)\quad\quad\forall \eta\leq\sigma\leq\eta+\vec k\,.
\end{equation}
Thus taking $\gamma:=\max\Big\{\beta ,e^{\alpha}\Big\}$
we finally deduce
\begin{equation}\label{jkhhhhhjrkjkjhoiio333}
\mathcal{L}^N\left(S_{\sigma}(u)\right)
\leq \gamma\,e^{-\alpha\frac{\sigma-\eta}{\vec k}}\mathcal{L}^N\left(S_{\eta}(u)\right)\quad\quad\forall 0\leq\eta\leq \sigma\,.
\end{equation}
\end{proof}

\section{$BMO$-interpolation in Lorentz spaces $L^{q,\gamma}$}
\label{sec:$BMO$-interpolation in Lorentz spaces}
 Recall the definition of Lorentz spaces:
\begin{definition}
\label{def:definition of Lorents space}
Let $X$ be a set and let $\mu$ be a measure\footnote{By a measure $\mu$ we always mean a set function $\mu:2^X \to [0,\infty]$ that is $\sigma$-subadditive and satisfies $\mu(\emptyset)=0$.} on it. Let $0 < q < \infty$ and $0 < \gamma \leq \infty$. The \emph{Lorentz space} $L^{q,\gamma}(X)$ consists of all $\mu$-measurable functions $u : X \to \mathbb{R}$ such that
\begin{equation}
\label{eq:definition of Lorentz quasi-norms}
\|u\|_{L^{q,\gamma}(X)} = 
\begin{cases}
\left[ \displaystyle\int_0^\infty  \Big( \,\mu\left( \Set{x \in X}[|u(x)|^q > t ] \right)\Big)^{\frac{\gamma}{q}} \,t^{\frac{\gamma}{q}-1}\, dt \right]^{1/\gamma}, & \text{if } \gamma < \infty, 
\\[10pt]
\left[\sup\limits_{0<t<\infty} \, t \, \mu\left( \Set{x \in X}[|u(x)|^q > t ] \right)\right]^{1/q}, & \text{if } \gamma = \infty,
\end{cases}
\end{equation}

is finite. In case $\gamma=\infty$ we use the symbols $L^q_w(X):=L^{q,\infty}(X)$ and $[u]_{L^q_w(X)}:=\|u\|_{L^{q,\infty}(X)}$. If $u\notin L^{q,\gamma}(X)$, then we set $\|u\|_{L^{q,\gamma}(X)}=\infty$.
\end{definition}

\begin{remark}
Note that \eqref{eq:definition of Lorentz quasi-norms} is well defined even if the function $u$ is not $\mu$-measurable, since the function $t \mapsto \mu\left( \Set{x \in X}[|u(x)|^q > t ] \right)$ is monotone non-increasing and therefore $\mathcal{L}^1$-measurable. 
\end{remark}

\begin{remark}
The space $L^q_w$ is called the \emph{weak Lebesgue space}, as it generalizes the Lebesgue space $L^q$ in the sense that, if $u$ is a $\mu$-measurable function, then we get by Chebyshev's inequality
\begin{equation}
\label{eq:Lqw is less than Lq}
[u]_{L^q_w(X)} \leq \|u\|_{L^q(X)},
\end{equation}
and therefore $L^q(X) \subset L^q_w(X)$.
Moreover, we also have for every $q,\gamma\in (0,\infty)$
\begin{multline}
\label{bbbbb}
[u]_{L^q_w(X)}=\left[\sup\limits_{0<s<\infty} \, s \, \mu\left( \Set{x \in X}[|u(x)|^q > s ] \right)\right]^{1/q}
\\
=\sup\limits_{0<s<\infty} \left[ \displaystyle\int_0^s  \Big( \,\mu\left( \Set{x \in X}[|u(x)|^q > s ] \right)\Big)^{\frac{\gamma}{q}} \frac{\gamma}{q}\,t^{\frac{\gamma}{q}-1}\, dt \right]^{1/\gamma}
\\
\leq\left(\frac{\gamma}{q}\right)^{\frac{1}{\gamma}}\sup\limits_{0<s<\infty} \left[ \displaystyle\int_0^s  \Big( \,\mu\left( \Set{x \in X}[|u(x)|^q > t ] \right)\Big)^{\frac{\gamma}{q}} \,t^{\frac{\gamma}{q}-1}\, dt \right]^{1/\gamma}
\\
=\left(\frac{\gamma}{q}\right)^{\frac{1}{\gamma}}\left[ \displaystyle\int_0^\infty  \Big( \,\mu\left( \Set{x \in X}[|u(x)|^q > t ] \right)\Big)^{\frac{\gamma}{q}} \,t^{\frac{\gamma}{q}-1}\, dt \right]^{1/\gamma}=\left(\frac{\gamma}{q}\right)^{\frac{1}{\gamma}}\|u\|_{L^{q,\gamma}(X)}. 
\end{multline}

Note also that the Lorentz space $L^{q,q}(X)$ coincides with the Lebesgue space $L^q(X)$ by Fubini's theorem: for $q<\infty$
\begin{equation}
\|u\|_{L^{q,q}(X)}=\left[ \displaystyle\int_0^\infty  \mu\left( \Set{x \in X}[|u(x)|^q > t ] \right) \, dt \right]^{1/q}=
\left[\int_{X}|u(x)|^q\,d\mu(x) \right]^{1/q}
=\|u\|_{L^{q}(X)}.
\end{equation}
\end{remark}
For an exposition of weak Lebesgue spaces and Lorentz spaces, see, for instance, \cite{Grafakos2008}.

Hereafter, we use the measure-theoretic conventions: $\sup \emptyset := -\infty$, $\inf \emptyset := \infty$, $\infty \cdot 0 := 0$, and $-\infty \cdot 0 := 0$. 

\subsection{Equivalence between super-level set measures and Lorentz norms of truncated BMO functions}
\begin{lemma}
\label{lem:the minimum theorem for Lorentz quasi-norms including characteristic function1} 
Let $X$ be a set and $\mu$ a measure on it. 
Let $E \subset X$, and let $\gamma, q \in (0,\infty)$. 
Suppose $u:E \to \R$ is a $\mu$-measurable function and let $\vec k \in [0,\infty)$. Then,
\begin{equation}
\label{eq:weak13lemma1} 
\left(\frac{q}{\gamma}\right)^{\frac{1}{\gamma}}\vec k\mu\left({\Set{x\in
E}[\left|u(x)\right|>\vec k]}\right)^{\frac{1}{q}}\leq\left(\frac{q}{\gamma}\right)^{\frac{1}{\gamma}}\left[\chi_{\{|u|>\vec k\}}u\right]_{L^{q}_{w}(E)}\leq
\left\|\chi_{\{|u|>\vec k\}}u\right\|_{L^{q,\gamma}(E)}.
\end{equation}
\end{lemma}

\begin{proof}
By \eqref{bbbbb} we have the second inequality in \eqref{eq:weak13lemma1}. Using the definition of the quasi-norm $L^q_w$, we get 
\begin{multline}
\label{eq:weak5'} 
\left[\chi_{\{|u|>\vec k\}}u\right]_{L^{q}_{w}(E)}
=\sup_{0<t<\infty}\left[t\,\mu\left({\Set{x\in
E}[\left|\chi_{\{|u|>\vec k\}}(x)u(x)\right|^q>t]}\right)\right]^{\frac{1}{q}}
\\
\geq \left[\vec k^q\,\mu\left({\Set{x\in
E}[\left|\chi_{\{|u|>\vec k\}}(x)u(x)\right|^q>\vec k^q]}\right)\right]^{\frac{1}{q}}
=\vec k\mu\left({\Set{x\in
E}[\left|u(x)\right|>\vec k]}\right)^{\frac{1}{q}}.
\end{multline}
\end{proof}

\begin{theorem}
\label{thm:the minimum theorem for Lorentz quasi-norms including characteristic function} 
Let $q,\gamma\in (0,\infty)$.

\begin{enumerate}
\item Let $C\subset\R^N$ be a cube. Let $u\in BMO(C)$ and $\textbf{k}\in [0,\infty)$ be such that
\begin{equation}
\label{eq:uniform bound on BMO norms13}
\left\|u\right\|_{BMO(C)}\leq
\vec k.
\end{equation}
Then,
\begin{multline}
\label{eq:weak13} 
\left(\frac{q}{\gamma}\right)^{\frac{1}{\gamma}}\vec k\mathcal{L}^N\left({\Set{x\in
C}[\left|u(x)-u_C\right|>\vec k]}\right)^{\frac{1}{q}}
\leq
\left(\frac{q}{\gamma}\right)^{\frac{1}{\gamma}}\left[\chi_{\{|u-u_C|>\vec k\}}(u-u_C)\right]_{L^{q}_{w}(C)}
\\
\leq
\left\|\chi_{\{|u-u_C|>\vec k\}}(u-u_C)\right\|_{L^{q,\gamma}(C)}
\leq
C(\gamma,q,N)\vec k\mathcal{L}^N\left({\Set{x\in
C}[\left|u(x)-u_C\right|>\vec k]}\right)^{\frac{1}{q}}.
\end{multline}

\item Let $u\in BMO(\R^N)$ and $\textbf{k}\in [0,\infty)$ be such that
\begin{equation}
\left\|u\right\|_{BMO(\R^N)}\leq
\vec k.
\end{equation}
Then,
\begin{multline}
\left(\frac{q}{\gamma}\right)^{\frac{1}{\gamma}}\vec k\mathcal{L}^N\left({\Set{x\in
\R^N}[\left|u(x)\right|>\vec k]}\right)^{\frac{1}{q}}
\leq
\left(\frac{q}{\gamma}\right)^{\frac{1}{\gamma}}\left[\chi_{\{|u|>\vec k\}}u\right]_{L^{q}_{w}(\R^N)}
\\
\leq
\left\|\chi_{\{|u|>\vec k\}}u\right\|_{L^{q,\gamma}(\R^N)}
\leq
C(\gamma,q,N)\vec k\mathcal{L}^N\left({\Set{x\in
\R^N}[\left|u(x)\right|>\vec k]}\right)^{\frac{1}{q}}.
\end{multline}
\end{enumerate}

Here, $C(\gamma,q,N)$ is a positive constant that depends only on
$\gamma$, $q$, and $N$.
\end{theorem}

\begin{proof}
We prove both cases simultaneously. Let $E$ be a cube or $\R^N$, in case $E$ is a cube we assume that $u_E=0$. Assume $\gamma<\infty$. By the definition of the quasi-Lorentz norm and additivity of integral we have
\begin{multline}
\label{eq:decomposition of Lqgamma}
\left\|\chi_{\{|u|>\vec k\}}u\right\|^{\gamma}_{L^{q,\gamma}(E)}=\int_0^\infty  \Big( \mathcal{L}^N\left( \Set{x \in E}[|\chi_{\{|u|>\vec k\}}(x)u(x)|^q > t ] \right)\Big)^{\frac{\gamma}{q}} \,t^{\frac{\gamma}{q}-1}dt
\\
=\int_0^{\vec k^q}  \Big( \mathcal{L}^N\left( \Set{x \in E}[|\chi_{\{|u|>\vec k\}}(x)u(x)|^q > t ] \right)\Big)^{\frac{\gamma}{q}} \,t^{\frac{\gamma}{q}-1}dt
\\
+\int_{\vec k^q} ^\infty  \Big( \mathcal{L}^N\left( \Set{x \in E}[|\chi_{\{|u|>\vec k\}}(x)u(x)|^q > t ] \right)\Big)^{\frac{\gamma}{q}} \,t^{\frac{\gamma}{q}-1}dt.
\end{multline}
Notice that
\begin{equation}
\label{eq:decomposition of Lqgamma1}
\int_0^{\vec k^q}  \Big( \mathcal{L}^N\left( \Set{x \in E}[|\chi_{\{|u|>\vec k\}}(x)u(x)|^q > t ] \right)\Big)^{\frac{\gamma}{q}} \,t^{\frac{\gamma}{q}-1}dt
=\frac{q}{\gamma}\vec k^\gamma\mathcal{L}^N\left( \Set{x \in E}[|u(x)| > \vec k ] \right)^{\frac{\gamma}{q}}.
\end{equation}
Assume that $\vec k>0$. Using Corollary \ref{cor:relation between Lebesgue measures of level sets of BMO} in the case where $E$ is a cube, and Corollary \ref{cor:exponential growth for level sets of BMO} in the case where $E=\R^N$, we get
\begin{multline}
\label{eq:decomposition of Lqgamma2}
\int_{\vec k^q} ^\infty  \Big( \mathcal{L}^N\left( \Set{x \in E}[|\chi_{\{|u|>\vec k\}}(x)u(x)|^q > t ] \right)\Big)^{\frac{\gamma}{q}} \,t^{\frac{\gamma}{q}-1}dt
\\
=\int_{\vec k^q} ^\infty  \Big( \mathcal{L}^N\left( \Set{x \in E}[|u(x)|^q > t ] \right)\Big)^{\frac{\gamma}{q}} \,t^{\frac{\gamma}{q}-1}dt
\leq \int_{\vec k^q} ^\infty  \Big( \beta e^{-\alpha\frac{t^{1/q}}{\vec k}}\mathcal{L}^N\left( \Set{x \in E}[|u(x)| > \vec k ] \right)\Big)^{\frac{\gamma}{q}} \,t^{\frac{\gamma}{q}-1}dt
\\
=\Big( \beta \mathcal{L}^N\left( \Set{x \in E}[|u(x)| > \vec k ] \right)\Big)^{\frac{\gamma}{q}} \int_{\vec k^q} ^\infty  e^{-\alpha\frac{t^{1/q}}{\vec k}}t^{\frac{\gamma}{q}-1}dt
\\
= \Big( \beta \mathcal{L}^N\left( \Set{x \in E}[|u(x)| > \vec k ] \right)\Big)^{\frac{\gamma}{q}} q\vec k^{\gamma}\int_{1} ^\infty  e^{-\alpha s}s^{\gamma-1}ds.
\end{multline}
From \eqref{eq:decomposition of Lqgamma}, \eqref{eq:decomposition of Lqgamma1}, and \eqref{eq:decomposition of Lqgamma2}, we get
\begin{multline}
\left\|\chi_{\{|u|>\vec k\}}u\right\|^{\gamma}_{L^{q,\gamma}(E)}
\leq \frac{q}{\gamma}\vec k^\gamma\mathcal{L}^N\left( \Set{x \in E}[|u(x)| > \vec k ] \right)^{\frac{\gamma}{q}}
\\
+\Big( \beta \mathcal{L}^N\left( \Set{x \in E}[|u(x)|> \vec k ] \right)\Big)^{\frac{\gamma}{q}} q\vec k^{\gamma}\int_{1} ^\infty  e^{-\alpha s}s^{\gamma-1}ds
\\
=\vec k^\gamma\mathcal{L}^N\left( \Set{x \in E}[|u(x)| > \vec k ] \right)^{\frac{\gamma}{q}}\left(\frac{q}{\gamma}+\beta^{\frac{\gamma}{q}}q\int_{1} ^\infty  e^{-\alpha s}s^{\gamma-1}ds\right). 
\end{multline}
This completes the proof in case $\gamma<\infty$. The case $\gamma=\infty$ is proved by similar considerations:
\begin{multline}
\left\|\chi_{\{|u|>\vec k\}}u\right\|^{q}_{L^{q,\infty}(E)}
=\sup\limits_{0<t<\infty}t\mathcal{L}^N\left( \Set{x \in E}[|\chi_{\{|u|>\vec k\}}(x)u(x)|^q > t ] \right)
\\
\leq \sup\limits_{0<t\leq \vec k^q}t\mathcal{L}^N\left( \Set{x \in E}[|\chi_{\{|u|>\vec k\}}(x)u(x)|^q > t ] \right)
\\
+\sup\limits_{\vec k^q<t<\infty}t\mathcal{L}^N\left( \Set{x \in E}[|\chi_{\{|u|>\vec k\}}(x)u(x)|^q > t ] \right).
\end{multline}
We have
\begin{equation}
\sup\limits_{0<t\leq \vec k^q}t\mathcal{L}^N\left( \Set{x \in E}[|\chi_{\{|u|>\vec k\}}(x)u(x)|^q > t ] \right)=\vec k^q \mathcal{L}^N\left( \Set{x \in E}[|u(x)| > \vec k ] \right),
\end{equation}
and 
\begin{multline}
\sup\limits_{\vec k^q<t<\infty}t\mathcal{L}^N\left( \Set{x \in E}[|\chi_{\{|u|>\vec k\}}(x)u(x)|^q > t ] \right)=\sup\limits_{\vec k^q<t<\infty}t\mathcal{L}^N\left( \Set{x \in E}[|u(x)|^q > t ] \right)
\\
\leq \sup\limits_{\vec k^q<t<\infty}t \beta e^{-\alpha\frac{t^{1/q}}{\vec k}}\mathcal{L}^N\left( \Set{x \in E}[|u(x)| > \vec k ] \right)
\\
=\vec k^q\mathcal{L}^N\left( \Set{x \in E}[|u(x)| > \vec k ] \right)\left(\beta\sup\limits_{1<s<\infty}s^q e^{-\alpha s}\right)\quad (t=(s\vec k)^q)
\end{multline}
\end{proof}

\subsection{BMO-interpolation in Lorentz spaces}

\begin{lemma}
\label{lem:the minimum theorem for Lorentz quasi-norms including characteristic function195} 
Let $X$ be a set and $\mu$ a measure on it. 
Let $E \subset X$, and let $0<p<q<\infty$, $\gamma \in (0,\infty]$. 
Suppose $u:E \to \R$ is a function and let $\vec k \in [0,\infty)$. Then,
\begin{equation}
\label{eq:weak7'} 
\Big\|\min\left\{\left|u\right|\,,\,\vec k\right\}\Big\|_{L^{q,\gamma}(E)}
\leq C\,\sup\limits_{0<s\leq \vec k^p}\Big[s\mu\left(\Set{x\in
E}[\left|u(x)\right|^p>s]\right)\Big]^{\frac{1}{q}}\vec k^{1-\frac{p}{q}}
\leq C\,[u]^{\frac{p}{q}}_{L^p_w(E)}\vec k^{1-\frac{p}{q}}.
\end{equation}
Here $C:=\left[\frac{\gamma}{q}\left(1-\frac{p}{q}\right)\right]^{-\frac{1}{\gamma}}$ in case $\gamma<\infty$; and $C=1$ in case $\gamma=\infty$.
\end{lemma}

\begin{proof}
Assume that $\gamma<\infty$.
\begin{multline}
\label{eq:ineterpolation inequality involving BMO and Lorentz7}
\Big\|\min\left\{\left|u\right|\,,\,\vec k\right\}\Big\|^{\gamma}_{L^{q,\gamma}(E)}
=\int_{0}^{\infty}\Big(\mu\left(\Set{x\in
E}[\min\left\{\left|u(x)\right|^q\,,\,\vec k^q\right\}>\sigma]\right)\Big)^{\frac{\gamma}{q}}\sigma^{\frac{\gamma}{q}-1} d\sigma
\\
=\int_{0}^{\vec k^q}\Big(\mu\left(\Set{x\in
E}[\min\left\{\left|u(x)\right|^q\,,\,\vec k^q\right\}>\sigma]\right)\Big)^{\frac{\gamma}{q}}\sigma^{\frac{\gamma}{q}-1} d\sigma
=\int_{0}^{\vec k^q}\Big(\mu\left(\Set{x\in
E}[\left|u(x)\right|^q>\sigma]\right)\Big)^{\frac{\gamma}{q}}\sigma^{\frac{\gamma}{q}-1} d\sigma
\\
=\int_{0}^{\vec k^q}\Big(\mu\left(\Set{x\in
E}[\left|u(x)\right|^p>\sigma^{\frac{p}{q}}]\right)\Big)^{\frac{\gamma}{q}}\sigma^{\frac{\gamma}{q}-1} d\sigma
=\int_{0}^{\vec k^q}\Big(\sigma^{\frac{p}{q}}\mu\left(\Set{x\in
E}[\left|u(x)\right|^p>\sigma^{\frac{p}{q}}]\right)\Big)^{\frac{\gamma}{q}}\sigma^{-\frac{\gamma}{q}\frac{p}{q}}\sigma^{\frac{\gamma}{q}-1} d\sigma
\\
\leq \sup\limits_{0<\sigma\leq \vec k^q}\Big(\sigma^{\frac{p}{q}}\mu\left(\Set{x\in
E}[\left|u(x)\right|^p>\sigma^{\frac{p}{q}}]\right)\Big)^{\frac{\gamma}{q}}
\int_{0}^{\vec k^q}\sigma^{\frac{\gamma}{q}\left(1-\frac{p}{q}\right)-1}d\sigma
\\
=\sup\limits_{0<t\leq \vec k^p}\Big(t\mu\left(\Set{x\in
E}[\left|u(x)\right|^p>t]\right)\Big)^{\frac{\gamma}{q}}\frac{1}{\frac{\gamma}{q}\left(1-\frac{p}{q}\right)}
\vec k^{\gamma\left(1-\frac{p}{q}\right)}
\leq \frac{1}{\frac{\gamma}{q}\left(1-\frac{p}{q}\right)}[u]_{L^p_w(E)}^{\gamma \frac{p}{q}} \vec k^{\gamma\left(1-\frac{p}{q}\right)}.
\end{multline}

For $\gamma=\infty$,
\begin{multline}
\label{eq:ineterpolation inequality involving BMO and Lorentz71}
\Big\|\min\left\{\left|u\right|\,,\,\vec k\right\}\Big\|^q_{L^{q,\infty}(E)}
=\sup_{0<\sigma<\infty}\sigma\,\mu\left(\Set{x\in
E}[\min\left\{\left|u(x)\right|^q\,,\,\vec k^q\right\}>\sigma]\right)
\\
=\sup_{0<\sigma<\vec k^q}\sigma\,\mu\left(\Set{x\in
E}[\min\left\{\left|u(x)\right|^q\,,\,\vec k^q\right\}>\sigma]\right)
=\sup_{0<\sigma<\vec k^q}\sigma\,\mu\left(\Set{x\in
E}[\left|u(x)\right|^q>\sigma]\right)
\\
=\sup_{0<\sigma<\vec k^q}\sigma\,\mu\left(\Set{x\in
E}[\left|u(x)\right|^p>\sigma^{\frac{p}{q}}]\right)
=\sup_{0<t<\vec k^p}t^{\frac{q}{p}}\,\mu\left(\Set{x\in
E}[\left|u(x)\right|^p>t]\right)
\\
=\sup_{0<t<\vec k^p}t^{\frac{q}{p}-1}\,t\,\mu\left(\Set{x\in
E}[\left|u(x)\right|^p>t]\right)\leq \vec k^{q-p}\sup_{0<t<\vec k^p}t\,\mu\left(\Set{x\in
E}[\left|u(x)\right|^p>t]\right)
\\
\leq \vec k^{q-p}[u]_{L^p_w(E)}^{p}.
\end{multline}
\end{proof}

\begin{theorem}
\label{thm:the minimum theorem for Lorentz quasi-norms} 
Let $0<p<q<\infty$, $\gamma\in (0,\infty]$.

\begin{enumerate}
\item Let $C\subset\R^N$ be a cube. Let $u\in BMO(C)$, and $\vec k\in [0,\infty)$ such that
\begin{equation}
\label{eq:uniform bound on BMO norms}
\left\|u\right\|_{BMO(C)}\leq
\vec k.
\end{equation}
Then,
\begin{multline}
\label{eq:weak} 
\|u-u_C\|_{L^{q,\gamma}(C)}\leq C(\gamma,q,N)
\Big\|\min\left\{\left|u-u_C\right|\,,\,\textbf{k}\right\}\Big\|_{L^{q,\gamma}(C)}
\\
\leq C(\gamma,q,N)\,\sup\limits_{0<s\leq \vec k^p}\Big[s\mathcal{L}^N\left(\Set{x\in
C}[\left|u(x)-u_C\right|^p>s]\right)\Big]^{\frac{1}{q}}\vec k^{1-\frac{p}{q}}
\\
\leq C(\gamma,q,N)\,[u-u_C]^{\frac{p}{q}}_{L^p_w(C)}\vec k^{1-\frac{p}{q}}.
\end{multline}

\item Let $u\in BMO(\R^N)$, and $\textbf{k}\in [0,\infty)$ such that
\begin{equation}
\left\|u\right\|_{BMO(\R^N)}\leq
\textbf{k}.
\end{equation}
Then,
\begin{multline}
\|u\|_{L^{q,\gamma}(\R^N)}\leq C(\gamma,q,N)
\Big\|\min\left\{\left|u\right|\,,\,\textbf{k}\right\}\Big\|_{L^{q,\gamma}(\R^N)}
\\
\leq C(\gamma,q,N)\,\sup\limits_{0<s\leq \vec k^p}\Big[s\mathcal{L}^N\left(\Set{x\in
\R^N}[\left|u(x)\right|^p>s]\right)\Big]^{\frac{1}{q}}\vec k^{1-\frac{p}{q}}
\\
\leq C(\gamma,q,N)\,[u]^{\frac{p}{q}}_{L^p_w(\R^N)}\vec k^{1-\frac{p}{q}}.
\end{multline}
\end{enumerate}

Here, $C(\gamma,q,N)$ is a positive constant that depends only on
$\gamma$, $q$, and $N$.
\end{theorem}

\begin{proof}
We prove both cases simultaneously. Let $E$ be a cube or $\R^N$. Assume that $\vec k=1$, otherwise we replace $u$ by $\vec k^{-1}u$. Note that the case $\vec k=0$ is trivial. Assume first $\gamma<\infty$. By Definition \ref{def:definition of Lorents space} and additivity of the integral, we get
\begin{multline}
\label{eq:estimate for integral of u in q through Fubini theorem}
\|u\|^\gamma_{L^{q,\gamma}(E)}
=\int_{0}^{\infty}\Big(\mathcal{L}^N\left(\Set{x\in
E}[|u(x)|^q>\sigma]\right)\Big)^{\frac{\gamma}{q}}\sigma^{\frac{\gamma}{q}-1} d\sigma
\\
=\int_{0}^{1}\Big(\mathcal{L}^N\left(\Set{x\in
E}[|u(x)|^q>\sigma]\right)\Big)^{\frac{\gamma}{q}}\sigma^{\frac{\gamma}{q}-1} d\sigma
\\
+\int_{1}^\infty\Big(\mathcal{L}^N\left(\Set{x\in
E}[|u(x)|^q>\sigma]\right)\Big)^{\frac{\gamma}{q}}\sigma^{\frac{\gamma}{q}-1} d\sigma
:=S_1+S_2.
\end{multline}

We intend to analyse the quantities $S_1$, $S_2$. We first treat $S_1$. We get
\begin{multline}
\label{eq:treatment of S1}
S_1=\int_{0}^{1}\Big(\mathcal{L}^N\left(\Set{x\in
E}[|u(x)|^q>\sigma]\right)\Big)^{\frac{\gamma}{q}}\sigma^{\frac{\gamma}{q}-1} d\sigma
\\
=\int_{0}^{1}\Big(\mathcal{L}^N\left(\Set{x\in
E}[\min\{|u(x)|^q,1\}>\sigma]\right)\Big)^{\frac{\gamma}{q}}\sigma^{\frac{\gamma}{q}-1} d\sigma
\\
=\int_{0}^{\infty}\Big(\mathcal{L}^N\left(\Set{x\in
E}[\min\{|u(x)|^q,1\}>\sigma]\right)\Big)^{\frac{\gamma}{q}}\sigma^{\frac{\gamma}{q}-1} d\sigma
=\|\min\left\{\left|u\right|\,,\,1\right\}\|^\gamma_{L^{q,\gamma}(E)}.
\end{multline}

We now treat the expression $S_2$. In case $E$ is a cube we assume $u_E=0$. By Theorem \ref{thm:the minimum theorem for Lorentz quasi-norms including characteristic function} we have
\begin{multline}
\label{eq:estimate for S2}
S_2=
\int_{1}^\infty\Big(\mathcal{L}^N\left(\Set{x\in
E}[|u(x)|^q>\sigma]\right)\Big)^{\frac{\gamma}{q}}\sigma^{\frac{\gamma}{q}-1} d\sigma
\\
=\int_{1}^\infty\Big(\mathcal{L}^N\left(\Set{x\in
E}[|u(x)\chi_{\{|u|>1\}}(x)|^q>\sigma]\right)\Big)^{\frac{\gamma}{q}}\sigma^{\frac{\gamma}{q}-1} d\sigma
\\
\leq \int_{0}^\infty\Big(\mathcal{L}^N\left(\Set{x\in
E}[|u(x)\chi_{\{|u|>1\}}(x)|^q>\sigma]\right)\Big)^{\frac{\gamma}{q}}\sigma^{\frac{\gamma}{q}-1} d\sigma
\\
=\|\chi_{\{|u|>1\}}u\|^{\gamma}_{L^{q,\gamma}(E)}
\leq C(\gamma,q,N)
\mathcal{L}^N\left({\Set{x\in
E}[\left|u(x)\right|>1]}\right)^{\frac{\gamma}{q}}.
\end{multline}

Notice that
\begin{multline}
\label{eq:estimate for S2 before starting cases}
\|\min\left\{\left|u\right|\,,\,1\right\}\|^\gamma_{L^{q,\gamma}(E)}
=\int_{0}^{\infty}\Big(\mathcal{L}^N\left(\Set{x\in
E}[\min\{|u(x)|^q,1\}>\sigma]\right)\Big)^{\frac{\gamma}{q}}\sigma^{\frac{\gamma}{q}-1} d\sigma
\\
= \int_{0}^{1}\Big(\mathcal{L}^N\left(\Set{x\in
E}[\min\{|u(x)|^q,1\}>\sigma]\right)\Big)^{\frac{\gamma}{q}}\sigma^{\frac{\gamma}{q}-1} d\sigma 
\\
\geq \Big(\mathcal{L}^N\left(\Set{x\in
E}[|u(x)|>1]\right)\Big)^{\frac{\gamma}{q}}\int_{0}^{1}\sigma^{\frac{\gamma}{q}-1} d\sigma
=\Big(\mathcal{L}^N\left(\Set{x\in
E}[|u(x)|>1]\right)\Big)^{\frac{\gamma}{q}}\frac{q}{\gamma}.
\end{multline}

Thus, from \eqref{eq:estimate for S2} and \eqref{eq:estimate for S2 before starting cases}, we obtain
\begin{equation}
\label{eq:estimate for S2 final}
S_2\leq \|\chi_{\{|u|>1\}}u\|^{\gamma}_{L^{q,\gamma}(E)}
\leq\frac{\gamma}{q} C(\gamma,q,N)\|\min\left\{\left|u\right|\,,\,1\right\}\|^\gamma_{L^{q,\gamma}(E)}.
\end{equation}

This together with \eqref{eq:estimate for integral of u in q through Fubini theorem}
and
\eqref{eq:treatment of S1}
completes the proof in the case $\gamma<\infty$.

Finally, for $\gamma=\infty$ (recall that $L^{q,\infty}=L^q_w$, $\|\cdot\|_{L^{q,\infty}}=[\cdot]_{L^q_w}$),
\begin{multline}
\label{ihjhkhkhkhk}
[u]_{L^q_w(E)}=\sup\limits_{s\in(0,\infty)}\Big[s\mathcal{L}^N\left(\Set{x\in
E}[\left|u(x)\right|^q>s]\right)\Big]^{\frac{1}{q}}
\\
\leq\sup\limits_{s\in(0,1)}\Big[s\mathcal{L}^N\left(\Set{x\in
E}[\left|u(x)\right|^q>s]\right)\Big]^{\frac{1}{q}}
+\sup\limits_{s\in[1,\infty)}\Big[s\mathcal{L}^N\left(\Set{x\in
E}[\left|u(x)\right|^q>s]\right)\Big]^{\frac{1}{q}}
\\
=\sup\limits_{s\in(0,1)}\Big[s\mathcal{L}^N\left(\Set{x\in
E}[\min\left\{\left|u(x)\right|,1\right\}^q>s]\right)\Big]^{\frac{1}{q}}
+\sup\limits_{s\in[1,\infty)}\Big[s\mathcal{L}^N\left(\Set{x\in
E}[\left|\chi_{\{|u|>1\}}(x)u(x)\right|^q>s]\right)\Big]^{\frac{1}{q}}
\\
\leq\sup\limits_{s\in (0,\infty)}\Big[s\mathcal{L}^N\left(\Set{x\in
E}[\min\left\{\left|u(x)\right|,1\right\}^q>s]\right)\Big]^{\frac{1}{q}}
+\sup\limits_{s\in (0,\infty)}\Big[s\mathcal{L}^N\left(\Set{x\in
E}[\left|\chi_{\{|u|>1\}}(x)u(x)\right|^q>s]\right)\Big]^{\frac{1}{q}}
\\
=[\min\left\{\left|u\right|,1\right\}]_{L^q_w(E)}+[\chi_{\{|u|>1\}}u]_{L^q_w(E)}.
\end{multline}

However, by Theorem \ref{thm:the minimum theorem for Lorentz quasi-norms including characteristic function} we have
\begin{multline}
\label{eq:weak13jkjkljkkj} 
\left[\chi_{\{|u|>1\}}u\right]_{L^{q}_{w}(E)}\leq C'\mathcal{L}^N\left({\Set{x\in
E}[\left|u(x)\right|^q>1]}\right)^{\frac{1}{q}}
\leq C'\sup_{s\in (0,1)}\left[s\,\mathcal{L}^N\left({\Set{x\in
E}[\left|u(x)\right|^q>s]}\right)\right]^{\frac{1}{q}}
\\
=C'\sup_{s\in (0,1)}\left[s\,\mathcal{L}^N\left({\Set{x\in
E}[\min\{\left|u(x)\right|^q,1\}>s]}\right)\right]^{\frac{1}{q}}
\leq C'[\min\left\{\left|u\right|,1\right\}]_{L^q_w(E)},
\end{multline}
where $C'$ depends on $\gamma,q,N$ only. Thus, from \eqref{ihjhkhkhkhk} and \eqref{eq:weak13jkjkljkkj}, we obtain
\begin{equation}
\label{ihjhkhkhkhkklklkl}
[u]_{L^q_w(E)}\leq\left(1+C'\right)[\min\left\{\left|u\right|,1\right\}]_{L^q_w(E)}.
\end{equation}
This completes the proof.
\end{proof}

The following proposition is a consequence of Theorem \ref{thm:the minimum theorem for Lorentz quasi-norms}. We will use it in the proofs of Theorem \ref{thm:BMO-interpolation for Besov functions} and Theorem \ref{thm: Besov constants vanish for VMO}.  
\begin{proposition}
\label{prop:interpolation for translations in BMO - weak case}
Let $0<p<q<\infty$. If $u\in BMO(\R^N)$, then for every $h\in \R^N$ and for every $\vec k\in [0,\infty)$ such that $\|u\|_{BMO(\R^N)}\leq \vec k$, we have
\begin{multline}
\label{eq:BMO and Lp weak interpolation on RN with Lp weak norm1}
\int_{\R^N} |u(x+h) - u(x)|^q \, dx
\leq C
\int_{\R^N}\Big(\min\left\{\left|u(x+h)-u(x)\right|\,,\,\textbf{k}\right\}\Big)^q\,dx
\\
\leq C\textbf{k}^{q-p} 
\left(\sup\limits_{0<s\leq \textbf{k}^p}s\mathcal{L}^N\big(\Set{x\in
\R^N}[|u(x+h) - u(x)|^p>s]\big)\right)
\\
\leq C\textbf{k}^{q-p}[u(\cdot+h) - u]_{L^p_w(\R^N)}^{p},
\end{multline}
where $C=C(p,q,N)$ is a constant that depends only on $p$, $q$, and $N$.
\end{proposition}

\begin{proof}
Let us define for $h \in \R^N$ the function $v_h(x) := u(x+h) - u(x)$. If $u \in BMO(\R^N)$, then $v_h\in BMO(\R^N)$ because
\begin{equation}
\label{eq:if u in BMO so is vh1}
\|v_h\|_{BMO(\R^N)} \leq 2\|u\|_{BMO(\R^N)}.
\end{equation}
Let us prove \eqref{eq:if u in BMO so is vh1}. Note that for any cube $C$, by the definition of $v_h$, the triangle inequality, and the change of variable formula, we have
\begin{multline}
\label{eq:BMO norm of translation1}
\fint_{C} \fint_{C} |v_h(z) - v_h(x)| \, dz \, dx = \fint_{C} \fint_{C} |u(z+h) - u(z) - (u(x+h) - u(x))| \, dz \, dx \\
\leq \fint_{C} \fint_{C} |u(z) - u(x)| \, dz \, dx + \fint_{C} \fint_{C} |u(z+h) - u(x+h)| \, dz \, dx \\
= \fint_{C} \fint_{C} |u(z) - u(x)| \, dz \, dx + \fint_{C+h} \fint_{C+h} |u(z) - u(x)| \, dz \, dx \leq 2\|u\|_{BMO(\R^N)}.
\end{multline}
In \eqref{eq:BMO norm of translation1}, we used that if $C$ is a cube, then $C+h$ is also a cube and $\mathcal{L}^N(C) = \mathcal{L}^N(C+h)$. Taking the supremum in \eqref{eq:BMO norm of translation1} over all cubes $C$, we get $\|v_h\|_{BMO(\R^N)} \leq 2\|u\|_{BMO(\R^N)}$. Applying Theorem \ref{thm:the minimum theorem for Lorentz quasi-norms}
with $\gamma=q$, $v_h \in BMO(\R^N)$, and $\vec k:=2\|u\|_{BMO(\R^N)}$, we obtain inequalities \eqref{eq:BMO and Lp weak interpolation on RN with Lp weak norm1}.
\end{proof}

\section{$BMO$-interpolation in weak fractional Sobolev spaces $W_w^{s,p}$}
\label{sec:$BMO$-interpolation in weak fractional Sobolev}
\label{sec:interpolation of weak Gagliardo spaces}
Recall the definition of the space of fractional Sobolev functions $W^{s,p}, s\in (0,1)$:
\begin{definition}[Fractional Sobolev functions]
\label{def:Gagliardo seminorm,intro}
Let $0< q < \infty$ and $s\in (0,1)$. For a Lebesgue measurable set $E$ and Lebesgue measurable function $u: E \to \R$, the Gagliardo quasi-semi-norm is defined as
\begin{equation}
\|u\|_{W^{s,q}(E)} :=
\left(\int_{E}\int_{E}\frac{|u(x)-u(z)|^q}{|x-z|^{sq+N}} \, dxdz\right)^{\frac{1}{q}}.
\end{equation}
We say that $u \in W^{s,q}(E)$ if and only if $u \in L^q(E)$ and $\|u\|_{W^{s,q}(E)} < \infty$.
\end{definition}

Notice that:
\begin{multline}
\label{eq:strong fractional Sobolev included in weak}
\|u\|_{W^{s,q}(\R^N)} =
\left(\int_{\R^N}\int_{\R^N}\frac{|u(x)-u(z)|^q}{|x-z|^{sq+N}} \, dxdz\right)^{\frac{1}{q}}=\left(\int_{\R^N}\int_{\R^N}\frac{|u(y+z)-u(z)|^q}{|y|^{sq+N}} \, dydz\right)^{\frac{1}{q}}
\\
=\left(\int_{\R^N}\frac{\|u(y+\cdot)-u\|^q_{L^q(\R^N)}}{|y|^{sq+N}} \, dy\right)^{\frac{1}{q}}\geq \left(\int_{\R^N}\frac{[u(y+\cdot)-u]^q_{L^q_w(\R^N)}}{|y|^{sq+N}} \, dy\right)^{\frac{1}{q}}
\\
=\left(\int_{\R^N}\sup_{0<t<\infty}t\mathcal{L}^N\Big(\Set{z\in
\R^N}[|u(y+z)-u(z)|^q>t]\Big)\frac{1}{|y|^{sq+N}} \, dy\right)^{\frac{1}{q}}
\\
\geq \sup_{0<t<\infty}\left(\int_{\R^N}t\mathcal{L}^N\Big(\Set{z\in
\R^N}[|u(y+z)-u(z)|^q>t]\Big)\frac{1}{|y|^{sq+N}} \, dy\right)^{\frac{1}{q}}:=[u]_{W^{s,q}_w(\R^N)}.
\end{multline}
Recall the definition of the weak space $W^{s,q}_w$, see Definition~\ref{def:weak fractional Sobolev space,intro}. Therefore, \eqref{eq:strong fractional Sobolev included in weak} shows that $W^{s,q}(\mathbb{R}^N)$ is continuously embedded into $W^{s,q}_w(\mathbb{R}^N)$. 
Note that if we define $\bar{u}(x,y):=u(x+y)-u(x)$, then  
\begin{equation}
[u]_{W^{s,q}_w(\R^N)}=[\bar{u}]_{L^{q}_w\left(\R^N\times\R^N,\mathcal{L}^N\times \frac{1}{|y|^{sp+N}}\mathcal{L}^N\right)}.
\end{equation} 
Indeed, let us denote $\mu(E):=\int_{E}\frac{1}{|y|^{sq+N}}dy$ for $\mathcal{L}^N$ measurable sets $E\subset \R^N$. Then,
\begin{multline}
[u]^q_{W^{s,q}_w(\R^N)}
=\sup_{0<t<\infty}\int_{\R^N}t\mathcal{L}^N\Big(\Set{x\in
\R^N}[|u(x+y)-u(x)|^q>t]\Big)d\mu(y)
\\
=\sup_{0<t<\infty}t\int_{\R^N}\int_{\R^N}\chi_{\Set{(x,y)\in
 \R^N\times \R^N}[|u(x+y)-u(x)|^q>t]}(x,y)dx\,d\mu(y)
\\
=\sup_{0<t<\infty}t\left(\mathcal{L}^N\times \mu\right)\Big(\Set{(x,y)\in
\R^N\times \R^N}[|\bar{u}(x,y)|^q>t]\Big)
=[\bar{u}]^q_{L^{q}_w\left(\R^N\times\R^N,\mathcal{L}^N\times \mu \right)}.
\end{multline} 

\begin{lemma}
\label{lem:interpolation for BMO weak}
Let $0<p<q<\infty$ and $r\in(0,\infty)$. Let $\mu$ be a $\sigma$-finite measure on a set $Y$. Given a non-empty at most countable set $\mathcal{A}$, for every $\alpha\in\mathcal{A}$ let $u_\alpha(x,y):\R^N\times Y\to\R$ be an $\mathcal{L}^N\times\mu$-measurable function such that $u_\alpha(\cdot,y)\in BMO(\R^N)$ for $\mu$-a.e. $y\in Y$ and for every $\alpha\in\mathcal{A}$. Assume that there exists $\vec k\in (0,\infty)$ such that
\begin{equation}
\label{eq:definition of k}
\operatorname*{ess-sup}_{\substack{y \in Y}}\|u_\alpha(\cdot,y)\|_{BMO(\R^N)}\leq \vec k\quad\quad \forall\alpha\in\mathcal{A}.
\end{equation}
Then
\begin{multline}
\label{eq:BMO and Lp weak interpolation on RN}
\left(\int_Y\left(\sup\limits_{\alpha\in\mathcal{A}}\int_{\R^N}\left|u_\alpha(x,y)\right|^qdx\right)^r d\mu(y)\right)^{\frac{1}{rq}}
\\
\leq C\,\vec k^{1-\frac{p}{q}}\left(\sup_{0<\tau\leq \vec k^p}\int_Y\left(\sup\limits_{\alpha\in\mathcal{A}}\,\tau\mathcal{L}^N\left(\R^N\cap\{|u_\alpha(\cdot,y)|^p>\tau\}\right)\right)^rd\mu(y)\right)^{\frac{1}{rq}}.
\end{multline}
Here $C:=C(p,q,r,N)$ is a constant that depends only on $p$, $q$, $r$ and
$N$.
\end{lemma}
\begin{proof}
By Theorem \ref{thm:the minimum theorem for Lorentz quasi-norms} and Fubini's theorem we get for every $\alpha\in\mathcal{A}$ and for $\mu$-almost every $y\in Y$:
\begin{multline}
\label{eq:estimate for double integral containing general space Y}
\int_{\R^N}\left|u_\alpha(x,y)\right|^qdx\leq C
\int_{\R^N}\min\left\{\left|u_\alpha(x,y)\right|^q\,,\,\vec k^q\right\}dx
\\
=C\int_0^\infty\mathcal{L}^N\left(\R^N\cap\left\{\min\left\{\left|u_\alpha(\cdot,y)\right|^q\,,\,\vec k^q\right\}>\sigma\right\}\right) d\sigma
\\
=C\int_0^{\vec k^q}\mathcal{L}^N\left(\R^N\cap\left\{\min\left\{\left|u_\alpha(\cdot,y)\right|^q\,,\,\vec k^q\right\}>\sigma\right\}\right) d\sigma
\\
=C\int_0^{\vec k^q}\mathcal{L}^N\left(\R^N\cap\{|u_\alpha(\cdot,y)|^q>\sigma\}\right) d\sigma,
\end{multline}
where $C=C(q,N)$ is a constant depending only on $q,N$.
Therefore, in the case $r\geq 1$ by Minkowski's integral inequality, we get from \eqref{eq:estimate for double integral containing general space Y}:

\begin{multline}
\label{eq:estimate for double integral containing general space Y1}
\left(\int_Y\left(\sup\limits_{\alpha\in\mathcal{A}}\int_{\R^N}\left|u_\alpha(x,y)\right|^qdx\right)^rd\mu(y)\right)^{\frac{1}{r}}\leq\\
 C
\left(\int_Y\left(\sup\limits_{\alpha\in\mathcal{A}}\int_0^{\vec k^q}\mathcal{L}^N\left(\R^N\cap\{|u_\alpha(\cdot,y)|^q>\sigma\}\right) d\sigma\right)^rd\mu(y)\right)^{\frac{1}{r}}\\
\leq
C
\left(\int_Y\left(\int_0^{\vec k^q}\sup\limits_{\alpha\in\mathcal{A}}\mathcal{L}^N\left(\R^N\cap\{|u_\alpha(\cdot,y)|^q>\sigma\}\right) d\sigma\right)^rd\mu(y)\right)^{\frac{1}{r}}
\\
\leq
C\int_0^{\vec k^q}\left(\int_Y\left(\sup\limits_{\alpha\in\mathcal{A}}\mathcal{L}^N\left(\R^N\cap\{|u_\alpha(\cdot,y)|^q>\sigma\}\right)\right)^rd\mu(y) \right)^{\frac{1}{r}}d\sigma
\\
=C\int_0^{\vec k^q}\sigma^{-p/q}\left(\int_Y\left(\sup\limits_{\alpha\in\mathcal{A}}\,\sigma^{p/q}\mathcal{L}^N\left(\R^N\cap\{|u_\alpha(\cdot,y)|^p>\sigma^{p/q}\}\right)\,\right)^rd\mu(y)\right)^{\frac{1}{r}}d\sigma
\\
\leq C\left(\int_0^{\vec k^q}\sigma^{-p/q}d\sigma
\right)\left(\sup\limits_{0<\sigma \leq \vec k^q}
\int_Y\left(\sup\limits_{\alpha\in\mathcal{A}}\,\sigma^{p/q}\mathcal{L}^N\left(\R^N\cap\{|u_\alpha(\cdot,y)|^p>\sigma^{p/q}\}\right)\,\right)^rd\mu(y)\right)^{\frac{1}{r}}
\\
=C\frac{q}{q-p}\vec k^{q-p}\left(\sup\limits_{0<s\leq \vec k^p}
\int_Y\left(\sup\limits_{\alpha\in\mathcal{A}}\,s\mathcal{L}^N\left(\R^N\cap\{|u_\alpha(\cdot,y)|^p>s\}\right)\,\right)^rd\mu(y)\right)^{\frac{1}{r}}.
\end{multline}
So we get 
\er{eq:BMO and Lp weak interpolation on RN} in the case $r\geq 1$.

If $r\in(0,1)$, then by \er{eq:estimate for double integral containing general space Y} for every $\alpha\in\mathcal{A}$ we get
\begin{multline}
\label{eq:estimate for double integral containing general space Y333}
\left(\int_{\R^N}\left|u_\alpha(x,y)\right|^qdx\right)^r
\leq C^r\left(\int_0^{\vec k^q}\mathcal{L}^N\left(\R^N\cap\{|u_\alpha(\cdot,y)|^q>\sigma\}\right) d\sigma\right)^r\\=
C^r
\left(\sum\limits_{j=0}^{\infty}\int_{\vec k^q 2^{-j-1}}^{\vec k^q2^{-j}}\mathcal{L}^N\left(\R^N\cap\{|u_\alpha(\cdot,y)|^q>\sigma\}\right) d\sigma\right)^r\\
\leq C^r
\left(\sum\limits_{j=0}^{\infty}\vec k^q 2^{-j-1}\mathcal{L}^N\left(\R^N\cap\{|u_\alpha(\cdot,y)|^q>\vec k^q 2^{-j-1}\}\right) \right)^r.
\end{multline}
For $r\in(0,1)$ and $A,B\geq 0$, it follows that $(A+B)^r\leq A^r+B^r$. Therefore,
\begin{multline}
\label{eq:estimate for double integral containing general space Y555}
\left(\sup\limits_{\alpha\in\mathcal{A}}\int_{\R^N}\left|u_\alpha(x,y)\right|^qdx\right)^r
\leq C^r
\left(\sum\limits_{j=0}^{\infty}\sup\limits_{\alpha\in\mathcal{A}}\vec k^q 2^{-j-1}\mathcal{L}^N\left(\R^N\cap\left\{|u_\alpha(\cdot,y)|^q>\vec k^q 2^{-j-1}\right\}\right)\right)^r
\\
\leq
\sum\limits_{j=0}^{\infty}C^r\left(\sup\limits_{\alpha\in\mathcal{A}}\vec k^q 2^{-j-1}\mathcal{L}^N\left(\R^N\cap\left\{|u_\alpha(\cdot,y)|^q>\vec k^q 2^{-j-1}\right\}\right)\right)^r=
\\
\sum\limits_{j=0}^{\infty}C^r\vec k^{(q-p)r} 2^{-\frac{(q-p)r(j+1)}{q}}\left(\sup\limits_{\alpha\in\mathcal{A}}\,\vec k^p 2^{-\frac{p(j+1)}{q}}\mathcal{L}^N\left(\R^N\cap\left\{|u_\alpha(\cdot,y)|^p>\vec k^p 2^{-\frac{p(j+1)}{q}}\right\}\right)\,\right)^r.
\end{multline}
Therefore,
\begin{multline}
\label{eq:estimate for double integral containing general space Y777}
\int_Y\left(\sup\limits_{\alpha\in\mathcal{A}}\int_{\R^N}\left|u_\alpha (x,y)\right|^qdx\right)^rd\mu(y)
\leq
 \\
\sum\limits_{j=0}^{\infty}C^r\vec k^{(q-p)r} 2^{-\frac{(q-p)r(j+1)}{q}}\int_Y\left(\sup\limits_{\alpha\in\mathcal{A}}\,\vec k^p 2^{-\frac{p(j+1)}{q}}\mathcal{L}^N\left(\R^N\cap\left\{|u_\alpha(\cdot,y)|^p>\vec k^p 2^{-\frac{p(j+1)}{q}}\right\}\right)\,\right)^rd\mu(y)
\\
\leq
\left(\sup\limits_{0<s\leq \vec k^p}
\int_Y\left(\sup\limits_{\alpha\in\mathcal{A}}\,s\mathcal{L}^N\left(\R^N\cap\{|u_\alpha(\cdot,y)|^p>s\}\right)\,\right)^rd\mu(y)\right)
\sum \limits_{j=0}^{\infty}C^r\vec k^{(q-p)r} 2^{-\frac{(q-p)r(j+1)}{q}}.
\end{multline}
Thus,
\begin{multline}
\label{eq:estimate for double integral containing general space Y999}
\left(\int_Y\left(\sup\limits_{\alpha\in\mathcal{A}}\int_{\R^N}\left|u_\alpha(x,y)\right|^qdx\right)^rd\mu(y)\right)^{\frac{1}{r}}
\leq 
\\
C\vec k^{q-p}\left(\sup\limits_{0<s\leq \vec k^p}
\int_Y\left(\sup\limits_{\alpha\in\mathcal{A}}\,s\mathcal{L}^N\left(\R^N\cap\{|u_\alpha(\cdot,y)|^p>s\}\right)\right)^rd\mu(y)\right)^{\frac{1}{r}}\left(\sum\limits_{j=0}^{\infty} 2^{-\frac{(q-p)r(j+1)}{q}}\right)^{\frac{1}{r}}
\\=
C\vec k^{q-p}\left(\sup\limits_{0<s\leq \vec k^p}
\int_Y\left(\sup\limits_{\alpha\in\mathcal{A}}s\mathcal{L}^N\left(\R^N\cap\{|u_\alpha(\cdot,y)|^p>s\}\right)\right)^rd\mu(y)\right)^{\frac{1}{r}}\frac{1}{\left(2^{\frac{(q-p)r}{q}}-1\right)^{\frac{1}{r}}}.
\end{multline}
 So we get 
\er{eq:BMO and Lp weak interpolation on RN} also in the case $r\in(0,1)$.
\end{proof}

\begin{proof}[Proof of Theorem \ref{thm:BMO-interpolation for fractional Sobolev functions,intro}]
We intend to use Lemma \ref{lem:interpolation for BMO weak}. We choose $Y=\R^N$, $\mu=\frac{1}{|z|^{sp+N}}\mathcal{L}^N$, $u(x,z):=u(z+x)-u(x)$. Note that since $u\in BMO(\R^N)$, we get that $u(\cdot,z)\in BMO(\R^N)$ for every $z\in\R^N$. Note also that
\begin{equation}
\label{eq:definition of k23}
\sup_{z\in\R^N}\|u(\cdot+z)-u\|_{BMO(\R^N)}\leq 2\|u\|_{BMO(\R^N)}<\infty.
\end{equation}
Thus, the conditions of Lemma \ref{lem:interpolation for BMO weak} are satisfied. Therefore,
\begin{multline}
\label{eq:BMO and Lp weak interpolation on RN11}
\int_{\R^N}\int_{\R^N}\frac{\left|u(z+x)-u(x)\right|^q}{|z|^{sp+N}} dxdz
\\
\leq C\|u\|_{BMO(\R^N)}^{q-p}\left(\sup_{0<\tau\leq \|u\|_{BMO(\R^N)}^p}\tau\int_{\R^N}\mathcal{L}^N\left({\Set{x\in
\R^N}[\left|u(z+x)-u(x)\right|^p>\tau]}\right)\frac{1}{|z|^{sp+N}} dz\right)
\\
\leq C\|u\|_{BMO(\R^N)}^{q-p}[u]^p_{W^{s,p}_w(\R^N)}\leq C\|u\|_{BMO(\R^N)}^{q-p}\|u\|^p_{W^{s,p}(\R^N)},
\end{multline}
where $C=C(p,q,N)$ is a constant depending only on $p,q,N$. 
\end{proof}

\begin{corollary}
\label{cor:BMO-interpolation for fractional Sobolev functions}
Let $1 <p \leq q < \infty$. If $u \in W^{\frac{1}{p},p}(\R^N) \cap BMO(\R^N)$, then $u \in W^{\frac{1}{q},q}(\R^N)$ and
\begin{equation}
\label{eq:interpolation for fractional Sobolev functions in BMO in case 1/q-differentiability}
\|u\|_{W^{1/q,q}(\R^N)}
\leq C(p,q,N) \|u\|^{\frac{p}{q}}_{W^{1/p,p}(\R^N)} \|u\|^{1-\frac{p}{q}}_{BMO(\R^N)},
\end{equation}
where $C(p,q,N)$ is a constant depending only on $p$, $q$, and $N$.
\end{corollary}

\begin{proof}
Assume $p<q$. Choose $s=\frac{1}{p}$ in Theorem \ref{thm:BMO-interpolation for fractional Sobolev functions,intro}.
\end{proof}

\section{$BMO$-interpolation of $W^{1,sp}$ into $W^{s,p}$}
\label{sec:BMO-interpolation of W1sp into Ws,p}
In this section, our goal is to prove in an alternative way the interpolation theorem established by Jean Van Schaftingen (see Theorem \ref{thm: Van Schaftigen's theorem,intro}), which concerns $BMO$-interpolation involving functions with fractional derivatives and functions with first weak derivatives.   
\\

The following theorem was proved by Hoai-Minh Nguyen in \cite{Nguyen2006}:
\begin{theorem}
\label{thm:Nguyen inequality}
Let \(1 < p < \infty\). Then
there exists a constant \(C_{N,p}\), depending only on \(N\) and \(p\), such that
\begin{equation}
\sup\limits_{0<\delta<\infty}\delta^p\int_{\R^N}\mathcal{L}^N\left({\Set{x\in \R^N}[|g(x+y)-g(x)|>\delta]}\right)\frac{1}{|y|^{N+p}}dy\leq C_{N,p} \int_{\mathbb{R}^N} |\nabla g(x)|^p \, dx,
\end{equation}
for all \(g \in W^{1,p}(\mathbb{R}^N)\).
\end{theorem}

We now give a definition of $BMO$-functions in open sets that is weaker than condition \eqref{def:BMO norm according to Van Schaftigen,intro}: 
\begin{definition}[Definition of $BMO$-functions in open sets in $\R^N$]
\label{def:definition of BMO with cubes inside open sets,intro}
Let \(\Omega\) be an open set and \(u \in L^1_{\text{loc}}(\Omega)\). We say that \(u \in BMO(\Omega)\) if and only if
\begin{equation}
\label{eq:definition of BMO with cubes inside open sets,intro}
\sup \left\{ \fint_{C}\fint_{C} |u(z) - u(y)|dy \, dz \mid C \subset \Omega \text{ is a cube} \right\} < \infty.
\end{equation}
\end{definition}

\begin{proof}[Proof of Theorem \ref{thm:BMO-interpolation for fractional Sobolev functions, the case s=1 in extension domain,intro}]
Let us denote
\begin{equation}
\vec k:=\|f\|_{BMO(\R^N)},\quad \mu(E):=\int_{E}\frac{1}{|y|^{N+sp}}dy,
\end{equation}
for $\mathcal{L}^N$-measurable sets $E\subset \R^N$. Assume that $\vec k>0$, otherwise $\|u\|_{BMO(\Omega)}\leq \|f\|_{BMO(\R^N)}=0$, so $u$ is a constant function almost everywhere, so the interpolation inequality \eqref{eq:interpolation for fractional Sobolev functions in BMO, extension domain,intro} holds.
By Lemma \ref{lem:interpolation for BMO weak}, we have
\begin{multline}
\label{eq:BMO and Lp weak interpolation on RN,translation case,usage in proof,p}
\int_{\R^N}\int_{\R^N}\left|f(x+y)-f(x)\right|^pdxd\mu(y)
\\
\leq C\,\vec k^{p-sp}\sup_{0<\tau\leq \vec k^{sp}}\tau\int_{\R^N}\mathcal{L}^N\left({\Set{z\in
\R^N}[\left|f(z+y)-f(z)\right|^{sp}>\tau]}\right) d\mu(y).
\end{multline}
Here $C=C(s,p,N)$ is a constant that depends only on $p$, $s$, and
$N$. From \eqref{eq:extension for BMO and Sobolev at the same time}, Theorem \ref{thm:Nguyen inequality} and \eqref{eq:BMO and Lp weak interpolation on RN,translation case,usage in proof,p}, we obtain
\begin{multline}
\|u\|_{BMO(\Omega)}^{p-sp}\left(\int_{\Omega}|\nabla u(x)|^{sp}dx\right)
\geq C\vec k^{p-sp}\left(\int_{\R^N}|\nabla f(x)|^{sp}dx\right)
\\
\geq C\vec k^{p-sp} \left(\sup\limits_{0<\delta<\infty}\delta^{sp}\int_{\R^N}\mathcal{L}^N\left({\Set{z\in \R^N}[|f(z+y)-f(z)|>\delta]}\right)\frac{1}{|y|^{N+sp}}dy\right)
\\
=C\vec k^{p-sp}\sup\limits_{0<\delta<\infty}\delta^{sp}\int_{\R^N}\mathcal{L}^N\left({\Set{z\in \R^N}[|f(z+y)-f(z)|>\delta]}\right)d\mu(y)
\\
=C\vec k^{p-sp}\sup\limits_{0<\tau<\infty}\tau\int_{\R^N}\mathcal{L}^N\left({\Set{z\in \R^N}[|f(z+y)-f(z)|^{sp}>\tau]}\right)d\mu(y)
\\
\geq C\int_{\R^N}\int_{\R^N}\left|f(x+y)-f(x)\right|^p\,dxd\mu(y)
\\
=C\int_{\R^N}\int_{\R^N}\frac{\left|f(x+y)-f(x)\right|^p}{|y|^{N+sp}}dxdy
=C\int_{\R^N}\left(\int_{\R^N}\frac{\left|f(x+y)-f(x)\right|^p}{|y|^{N+sp}}dy\right)dx
\\
= C\int_{\R^N}\int_{\R^N}\frac{\left|f(z)-f(x)\right|^p}{|z-x|^{N+sp}}dzdx\geq C\|u\|^p_{W^{s,p}(\Omega)}.
\end{multline}
\end{proof}

\begin{corollary}
\label{cor:BMO-interpolation for fractional Sobolev functions, the case s=1}
Let $1<p< \infty$ and $s\in (1/p,1)$. If $u \in W^{1,sp}(\R^N) \cap BMO(\R^N)$, then $u \in W^{s,p}(\R^N)$ and
\begin{equation}
\label{eq:interpolation for fractional Sobolev functions in BMO}
\|u\|_{W^{s,p}(\R^N)}
\leq C\, \|\nabla u\|^{s}_{L^{sp}(\R^N)} \|u\|^{1-s}_{BMO(\R^N)},
\end{equation}
where $C=C(s,p,N)$ is a constant depending only on $s$, $p$, and $N$.
\end{corollary}

\begin{proof}
The result follows from Theorem \ref{thm:BMO-interpolation for fractional Sobolev functions, the case s=1 in extension domain,intro} by choosing $\Omega=\R^N$.
\end{proof}

\section{$BMO$-interpolation in weak Besov spaces $(B^s_{p,q})_w$}

\begin{definition}[$k$-differences]
Let \( f:\R^N\to \R \) be a function and \( k \in \mathbb{N} \). The \( k \)-th order difference of \( f \) with step size \( h \in \mathbb{R}^N \) is defined as:
\begin{equation}
\label{eq:finite difference}
\Delta^k_h f(x):= \sum_{j=0}^k \binom{k}{j}(-1)^{k-j} f(x + jh),
\end{equation}
where \( \binom{k}{j} = \frac{k!}{j!(k-j)!} \) is the binomial coefficient.
\end{definition}

\begin{definition}[Modulus of continuity]
\label{Mdcnt}
Let $0< p\leq \infty$.
Given $f\in L^p(\R^N)$, the strong modulus of continuity of $f$ of order $k$ in $L^p$ is defined to be a function of $t\in (0,\infty)$ defined by:
\begin{equation}\label{eq:definition of k-modulus of continuity}
\Omega_{k}(f,t)_{L^p}:=\sup_{|h|\leq t}\|\Delta^k_h f\|_{L^p(\R^N)}\\
\end{equation}
and 
given $f\in L^p_w(\R^N)$, the weak modulus of continuity of $f$ of order $k$ in $L^p_w$ is defined to be a function of $t\in (0,\infty)$ defined by:
\begin{equation}\label{eq:definition of k-modulus of continuity weak}
\Omega_{k}(f,t)_{L^p_w}:=\sup_{|h|\leq t}\|\Delta^k_h f\|_{L^p_w(\R^N)}.
\end{equation}
\end{definition}

Note that by Chebychev's inequality, for $f\in L^p(\R^N)$ we have
\begin{equation}\label{eq:definition of k-modulus of continuity111}
\Omega_{k}(f,t)_{L^p_w}\leq\Omega_{k}(f,t)_{L^p}.
\end{equation}
Notice that
\begin{multline}
\|\Delta_h^k f\|_{BMO(\mathbb{R}^N)}
=\sup_{C}\fint_C\fint_C|\Delta_h^k f(x)-\Delta_h^k f(z)|dxdz
\\
=\sup_{C}\fint_C\fint_C\left|\sum_{j=0}^k \binom{k}{j}(-1)^{k-j} f(x + jh)-\sum_{j=0}^k \binom{k}{j}(-1)^{k-j} f(z + jh)\right|dxdz
\\
\leq \sum_{j=0}^k \binom{k}{j}\sup_{C}\fint_C\fint_C\left|f(x + jh)-f(z + jh)\right|dxdz=
 2^k \|f\|_{BMO(\mathbb{R}^N)}.
\end{multline}
Therefore, in case $f\in BMO(\R^N)$, and $0<\lambda<p<\infty$, it follows from Theorem \ref{thm:the minimum theorem for Lorentz quasi-norms} 
\begin{multline}
\label{eq:modulus of continuity interpolation}
\Omega_{k}(f,t)_{L^p}=\sup_{|h|\leq t}\|\Delta^k_h f\|_{L^p(\R^N)}
\leq C(p,N)\sup_{|h|\leq t}\left(\|\Delta^k_h f\|^{1-\frac{\lambda}{p}}_{BMO(\R^N)}\|\Delta^k_h f\|^{\frac{\lambda}{p}}_{L^\lambda_w(\R^N)}\right)
\\
\leq  C(p,N)2^{k\left(1-\frac{\lambda}{p}\right)} \|f\|^{1-\frac{\lambda}{p}}_{BMO(\R^N)} \left(\Omega_{k}(f,t)_{L^\lambda_w}\right)^{\frac{\lambda}{p}}.
\end{multline}

\begin{definition}[Besov space]
\label{def:definition of general Besov space1}
For $s\in (0,\infty)$, let $k\in \N$ be minimal such that $s<k$. Let $0<q\leq \infty$, $0< p\leq\infty$. We say that a function $f\in L^p(\R^N)$ belongs to $B^s_{p,q}(\R^N)$ if and only if
\begin{equation}
\label{eq:def of Besov quasi-norm}
\|f\|_{B^s_{p,q}(\R^N)}:=
\begin{cases}
\left(\int_0^\infty\left(\frac{\Omega_k(f,t)_{L^p}}{t^s}\right)^q\frac{dt}{t}\right)^{\frac{1}{q}}<\infty;&\quad if\,\, q<\infty
\\
\sup\limits_{0<t<\infty}\frac{\Omega_k(f,t)_{L^p}}{t^s}<\infty;&\quad if\,\, q=\infty
\end{cases}.
\end{equation}
Moreover,  we say that a function $f\in L^p_w(\R^N)$ belongs to $(B^s_{p,q})_w(\R^N)$ if and only if
\begin{equation}
\|f\|_{(B^s_{p,q})_w(\R^N)}:=
\begin{cases}
\left(\int_0^\infty\left(\frac{\Omega_k(f,t)_{L^p_w}}{t^s}\right)^q\frac{dt}{t}\right)^{\frac{1}{q}}<\infty;&\quad if\,\, q<\infty
\\
\sup\limits_{0<t<\infty}\frac{\Omega_k(f,t)_{L^p_w}}{t^s}<\infty;&\quad if\,\, q=\infty
\end{cases}.
\end{equation}
Note that $f\in B^s_{p,q}(\R^N)$ implies $f\in (B^s_{p,q})_w(\R^N)$.
\end{definition}

The following theorem proves Theorem \ref{$BMO$-interpolation for Besov functions,intro}.
\begin{theorem}
Let $0<p<\infty$, $0< q\leq \infty$, $s\in (0,\infty)$ and $0<\lambda<p$. There exists constant $C=C(N,s,p,q,\lambda)$ such that for every $f\in BMO(\R^N)\cap (B^{s}_{\lambda,q\frac{\lambda}{p}})_w(\R^N)$ 
\begin{equation}
\label{eq:interpolations in Besov Bspq1}
\begin{cases}
\|f\|_{B^{s\frac{\lambda}{p}}_{p,q}(\R^N)}
\leq C \|f\|^{1-\frac{\lambda}{p}}_{BMO(\R^N)}\|f\|^{\frac{\lambda}{p}}_{(B^{s}_{\lambda,q\frac{\lambda}{p}})_w(\R^N)}\leq C \|f\|^{1-\frac{\lambda}{p}}_{BMO(\R^N)}\|f\|^{\frac{\lambda}{p}}_{B^{s}_{\lambda,q\frac{\lambda}{p}}(\R^N)}\quad\quad  &q<\infty
\\
\|f\|_{B^{s\frac{\lambda}{p}}_{p,\infty}(\R^N)}
\leq C \|f\|^{1-\frac{\lambda}{p}}_{BMO(\R^N)}\|f\|^{\frac{\lambda}{p}}_{(B^{s}_{\lambda,\infty})_w(\R^N)}\leq C \|f\|^{1-\frac{\lambda}{p}}_{BMO(\R^N)}\|f\|^{\frac{\lambda}{p}}_{B^{s}_{\lambda,\infty}(\R^N)}\quad\quad  &q=\infty\,.
\end{cases}
\end{equation}
\end{theorem}
\begin{proof} 
Let $k\in \N$ be a minimal number such that $s<k$ and $k'\in \N$ be a minimal number such that $s\frac{\lambda}{p}<k'\in\N$, as in Definition \ref{def:definition of general Besov space1}. Recall that,  it is well known (see Theorem \ref{thm:equivalence of Besov quasi-norms}) that
since $s\frac{\lambda}{p}\leq s<k$, the quasi-norms of the space $B^{s\frac{\lambda}{p}}_{p,q}(\mathbb{R}^N)$ obtained via the modulus of continuity of degrees $k$ and $k'$ are equivalent, meaning that:
If $q < \infty$,
\begin{equation}
\label{eq:equivalence for Besov norms in case of finite q1}
\left( \displaystyle\int_0^\infty  \left( \frac{\Omega_{k'}(f,t)_{L^p}}{t^{s\frac{\lambda}{p}}} \right)^{q} \frac{dt}{t} \right)^{\! 1/q} \sim
\left( \displaystyle\int_0^\infty  \left( \frac{\Omega_k(f,t)_{L^p}}{t^{s\frac{\lambda}{p}}} \right)^{q} \frac{dt}{t} \right)^{\! 1/q}.
\end{equation}
If $q = \infty$,
\begin{equation}
\label{eq:equivalence for Besov norms in case of infinite q2}
\sup_{t\in (0,\infty)} \frac{\Omega_{k'}(f,t)_{L^p}}{t^{s\frac{\lambda}{p}}} \sim
\sup_{t\in (0,\infty)} \frac{\Omega_k(f,t)_{L^p}}{t^{s\frac{\lambda}{p}}}.
\end{equation} 
Thus, in case $q<\infty$, by \eqref{eq:modulus of continuity interpolation}
\begin{multline}
\label{eq:estimate for f(+h)-f1jhjjkhjhj}
\left( \displaystyle\int_0^\infty \left( \frac{\Omega_{k'}(f,t)_{L^p}}{t^{s\frac{\lambda}{p}}} \right)^{q}\frac{dt}{t} \right)^{\! 1/q} \leq C \left(\int\limits_0^\infty\left(\frac{\Omega_k(f,t)_{L^p}}{t^{s\frac{\lambda}{p}}}\right)^q\frac{dt}{t}\right)^{\frac{1}{q}}
\\ 
\leq C \|f\|^{1-\frac{\lambda}{p}}_{BMO(\R^N)}\left[\left(\int\limits_0^\infty\left(\frac{\Omega_k(f,t)_{L^\lambda_w}}{t^{s}}\right)^{\frac{q\lambda}{p}}\frac{dt}{t}\right)^{\frac{p}{q\lambda}}\right]^{\frac{\lambda}{p}}
=C \|f\|^{1-\frac{\lambda}{p}}_{BMO(\R^N)}\|f\|^{\frac{\lambda}{p}}_{(B^{s}_{\lambda,q\frac{\lambda}{p}})_w(\R^N)}.
\end{multline} 
In the case $q=\infty$
\begin{multline}\label{eq:estimate for f(+h)-f1jhjjgfujjkjk}
\sup_{t\in (0,\infty)} \frac{\Omega_{k'}(f,t)_{L^p}}{t^{s\frac{\lambda}{p}}}\leq C\sup\limits_{t\in (0,\infty)}\frac{\Omega_k(f,t)_{L^p}}{t^{s\frac{\lambda}{p}}}
\\
\leq C\|f\|^{1-\frac{\lambda}{p}}_{BMO(\R^N)}\left[\sup\limits_{t\in (0,\infty)}\frac{\Omega_{k}(f,t)_{L^\lambda_w}}{t^{s}}\right]^{\frac{\lambda}{p}}=C \|f\|^{1-\frac{\lambda}{p}}_{BMO(\R^N)}\|f\|^{\frac{\lambda}{p}}_{(B^{s}_{\lambda,\infty})_w(\R^N)}.
\end{multline}
Note, however, that the case $q = \infty$ with $s < 1$ has already been established; 
see Theorem~\ref{thm:BMO-interpolation for Besov functions}.
\end{proof}

For $m\in\mathbb{N}$ and $1\leq p<\infty$, we denote the homogeneous Sobolev seminorm by
\[
\|u\|_{{W}^{m,p}(\mathbb{R}^N)}
:=
\left(
\sum_{|\alpha|=m}
\|D^\alpha u\|_{L^p(\mathbb{R}^N)}^p
\right)^{1/p}.
\]
\begin{corollary}
Let $m\in\mathbb{N}$, $0<r<1$, and $1\leq p<\infty$ such that 
\begin{equation}
\frac{2}{p}\leq r<\frac{1}{m}.
\end{equation}
Then there exists a constant $C=C(N,m,p,r)>0$ such that for every $u\in BMO(\R^N)\cap{W}^{m,rp}(\mathbb{R}^N)$
\begin{equation}
\|u\|_{{W}^{mr,p}(\mathbb{R}^N)}
\leq C\|u\|^{1-r}_{BMO(\mathbb{R}^N)}
\|u\|_{{W}^{m,rp}(\mathbb{R}^N)}^{r}.
\end{equation}
\end{corollary}

\begin{proof}

Let $m\in\N$ and $1<p<\infty$. It is known that
\begin{equation}
\|u\|_{{B}^m_{p,\max(p,2)}(\mathbb{R}^N)}\leq C\|u\|_{{W}^{m,p}(\mathbb{R}^N)
}.
\end{equation}
For proof see \cite[Theorem 2.26]{Sawano2018} together with \cite[5.2.3]{Triebel1983}. 
Choosing \(s=m\), \(\lambda=rp\), and \(q=p\), we obtain from \eqref{eq:interpolations in Besov Bspq1}
\begin{multline}
\|u\|_{{B}^{mr}_{p,p}(\mathbb{R}^N)}
\leq C \|u\|^{1-r}_{BMO(\mathbb{R}^N)}
\|u\|_{\left(B^{m}_{rp,rp}\right)_w(\mathbb{R}^N)}^{r}
\leq C \|u\|^{1-r}_{BMO(\mathbb{R}^N)}
\|u\|_{{B}^{m}_{rp,rp}(\mathbb{R}^N)}^{r}
\\
= C\|u\|^{1-r}_{BMO(\mathbb{R}^N)}
\|u\|_{{B}^{m}_{rp,\max(rp,2)}(\mathbb{R}^N)}^{r}
\leq C\|u\|^{1-r}_{BMO(\mathbb{R}^N)}
\|u\|_{{W}^{m,rp}(\mathbb{R}^N)}^{r},
\end{multline}
provided that \(rp\geq 2\). Since \(0<mr<1\), then ${B}^{mr}_{p,p}(\mathbb{R}^N)={W}^{mr,p}(\mathbb{R}^N)$ with equivalent semi-norms; see \cite[Proposition 14.40]{Leoni2017}. Therefore,
\begin{equation}
\label{eq:inter-Sobolev-high deri}
\|u\|_{{W}^{mr,p}(\mathbb{R}^N)}
\leq C\|u\|^{1-r}_{BMO(\mathbb{R}^N)}
\|u\|_{{W}^{m,rp}(\mathbb{R}^N)}^{r}.
\end{equation}
\end{proof}

\section{$BMO$-interpolation in weak Besov spaces $B^{s,p}_w$}
In this section, we establish the $BMO$-interpolation theorem for Besov functions in $B^{s,p}_w$; see Theorem \ref{thm:BMO-interpolation for Besov functions}. We prove a sharp result concerning the limiting Besov constant for $VMO$ functions; see Theorem \ref{thm: Besov constants vanish for VMO} and Remark \ref{rem: the limit of Besov constant in case q=p is not necessarily zero}. Additionally, we prove an interpolation result for weak functions of bounded variation; see Corollary \ref{cor:BMO-interpolation for BV}.
\\
\subsection{$BMO$-interpolation for Besov functions $B^{s,p}_w$}
\begin{definition}[Definition of the space $B^{s,q}$]
\label{def:Besov functions}
Let $0< q < \infty$ and $s\in (0,1]$. For a measurable function $u: \R^N \to \R$, we define the quasi-semi-norm
\begin{equation}
\|u\|_{B^{s,q}(\R^N)} :=\sup_{h \in \R^N \setminus \{0\}} \frac{\|u(\cdot+h)-u\|_{L^q(\R^N)}}{|h|^s}
\\
=\sup_{h \in \R^N \setminus \{0\}} \left(\int_{\R^N} \frac{|u(x+h)-u(x)|^q}{|h|^{sq}} \, dx\right)^{\frac{1}{q}}.
\end{equation}
We say that $u \in B^{s,q}(\R^N)$ if and only if $u \in L^q(\R^N)$ and $\|u\|_{B^{s,q}(\R^N)} < \infty$. 
\end{definition}

\begin{remark}
\label{rem:Bsq coninsides with different function spaces}
The space $B^{s,q}(\mathbb{R}^N)$ equals the Besov space $B^{s}_{q,\infty}(\mathbb{R}^N)$ for $s \in (0,1)$. If $s = 1$, then for $q \in (1, \infty)$, the space $B^{1,q}(\mathbb{R}^N)$ coincides with the Sobolev space $W^{1,q}(\mathbb{R}^N)$ 
, and for $q = 1$, it coincides with the space of functions of bounded variation, i.e., $B^{1,1}(\mathbb{R}^N) = BV(\mathbb{R}^N)$. 
\end{remark}

\begin{remark}
If we define the space $B^{s,q}(\mathbb{R}^N)$ allowing the differentiability degree $s$ to be greater than $1$, we would find that $B^{s,q}(\mathbb{R}^N)$ contains only functions that are equal to zero almost everywhere. More precisely, assume $q \in [1, \infty)$ and $s \in (1, \infty)$, and
\begin{equation}
\sup_{h \in \mathbb{R}^N \setminus \{0\}} \int_{\mathbb{R}^N} \frac{|u(x+h)-u(x)|^q}{|h|^{sq}} \, dx < \infty.
\end{equation}
Let us denote $A_h := \int_{\mathbb{R}^N} \frac{|u(x+h)-u(x)|^q}{|h|^{sq}} \, dx$. Since the family $\{A_h\}_{h \in \mathbb{R}^N \setminus \{0\}}$ is bounded and the term $|h|^{(s-1)q}$ converges to zero as $h \to 0$, we get
\begin{equation}
\label{eq:Sobolev norm equals zero}
0 = \lim_{h \to 0} A_h |h|^{(s-1)q} = \lim_{h \to 0} \int_{\mathbb{R}^N} \frac{|u(x+h)-u(x)|^q}{|h|^q} \, dx.
\end{equation}
It is known that the last term in \eqref{eq:Sobolev norm equals zero} is equivalent to $\|\nabla u\|^q_{L^q(\mathbb{R}^N)}$. Hence, we conclude that $u=0$ almost everywhere in $\mathbb{R}^N$.
\end{remark}

We give now a weak version of Definition \ref{def:Besov functions}:
\begin{definition}[Definition of the weak space $B_w^{s,q}$]
\label{def:Besov functions weak}
Let $0< q < \infty$ and $s\in (0,1]$. For a measurable function $u: \R^N \to \R$, we define the quasi-semi-norm 
\begin{multline}
[u]_{B^{s,q}_w(\R^N)}:=\sup_{h \in \R^N \setminus \{0\}} \frac{[u(\cdot+h)-u]_{L^q_w(\R^N)}}{|h|^s}
\\
=\sup_{h \in \R^N \setminus \{0\}} \frac{1}{|h|^s}\sup\limits_{0<s<\infty}\Big[s\mathcal{L}^N\big(\Set{x\in
\R^N}[|u(x+h)-u(x)|^q>s]\big)\Big]^{\frac{1}{q}}.
\end{multline}
We say that $u \in B^{s,q}_w(\R^N)$ if and only if $u \in L^q_w(\R^N)$ and $[u]_{B^{s,q}_w(\R^N)} < \infty$.
\end{definition}
\begin{remark}
\label{rem:Bsq coninsides with different function spaces111}
Notice again that the space $B^{s,q}(\mathbb{R}^N)$ equals the Besov space $B^{s}_{q,\infty}(\mathbb{R}^N)$ for $s \in (0,1)$. 
\end{remark}

\begin{theorem}[$BMO$-interpolation for weak Besov functions in $B^{s,q}_w$]
\label{thm:BMO-interpolation for Besov functions}
Let $0<p<q<\infty$ and $s\in (0,1]$. If $u\in BMO(\R^N)$, then we get
\begin{enumerate}
\item \begin{equation}
\label{eq:interpolation for Besov functions}
\|u\|_{B^{s\frac{p}{q},q}(\R^N)}\leq C[u]^{\frac{p}{q}}_{B^{s,p}_w(\R^N)}
\|u\|_{BMO(\R^N)}^{1-\frac{p}{q}}\leq C\|u\|^{\frac{p}{q}}_{B^{s,p}(\R^N)}
\|u\|_{BMO(\R^N)}^{1-\frac{p}{q}}.
\end{equation}
Here $C$ is a constant depending only on $p,q,N$.

\item For every $h\in\R^N$
\begin{multline}
\label{eq:upper directional interpolation for Besov functions with weight}
\limsup_{t\to 0}\frac{\|u(\cdot+th) - u\|_{L^q(\R^N)}}{|th|^{s\frac{p}{q}}}
\leq C\, \limsup_{t\to 0}\left(\frac{[u(\cdot+th)-u]_{L^p_w(\R^N)}}{|th|^s}\right)^{\frac{p}{q}}\|u\|_{BMO(\R^N)}^{1-\frac{p}{q}}
\\
\leq C\,\limsup_{t\to 0}\left(\frac{\|u(\cdot+th)-u\|_{L^p(\R^N)}}{|th|^s}\right)^{\frac{p}{q}}\|u\|_{BMO(\R^N)}^{1-\frac{p}{q}},
\end{multline}
and 
\begin{multline}
\label{eq:lower directional interpolation for Besov functions with weight}
\liminf_{t\to 0}\frac{\|u(\cdot+th) - u\|_{L^q(\R^N)}}{|th|^{s\frac{p}{q}}}
\leq C\, \liminf_{t\to 0}\left(\frac{[u(\cdot+th)-u]_{L^p_w(\R^N)}}{|th|^s}\right)^{\frac{p}{q}}\|u\|_{BMO(\R^N)}^{1-\frac{p}{q}}
\\
\leq C\,\liminf_{t\to 0}\left(\frac{\|u(\cdot+th)-u\|_{L^p(\R^N)}}{|th|^s}\right)^{\frac{p}{q}}\|u\|_{BMO(\R^N)}^{1-\frac{p}{q}},
\end{multline}
\end{enumerate}
where $C=C(p,q,N)$.
\end{theorem}

\begin{proof}
By Proposition \ref{prop:interpolation for translations in BMO - weak case}, we get for every $h \in \R^N$
\begin{equation}
\label{eq:interpolation for vh(1)2}
\int_{\R^N} |u(x+h) - u(x)|^q\,
dx \leq C(p,q,N) [u(\cdot+h)-u]^p_{L^p_w(\R^N)}
\|u\|_{BMO(\R^N)}^{q-p}.
\end{equation}
Multiplying both sides of \eqref{eq:interpolation for vh(1)2}
by
$|h|^{-sp}$, for $h \in \R^N \setminus \{0\}$, we get
\begin{equation}
\label{eq:interpolation for vh(2)}
\int_{\R^N} \frac{|u(x+h) - u(x)|^q}{|h|^{s p}}\, dx
\leq C(p,q,N) \left(\frac{ [u(\cdot+h)-u]^p_{L^p_w(\R^N)}}{|h|^{sp}}\right) \|u\|_{BMO(\R^N)}^{q-p}.
\end{equation}
Taking the supremum over $h \in \R^N \setminus \{0\}$ on both sides of \eqref{eq:interpolation for vh(2)}, we have
\begin{equation}
\label{eq:interpolation for vh(3)}
\sup_{h \in \R^N \setminus \{0\}} \int_{\R^N} \frac{|u(x+h) - u(x)|^q}{|h|^{s p}}\, dx
\leq C\, \left(\sup_{h \in \R^N \setminus \{0\}}\frac{ [u(\cdot+h)-u]^p_{L^p_w(\R^N)}}{|h|^{sp}}\right) \|u\|_{BMO(\R^N)}^{q-p},
\end{equation}
where $C=C(p,q,N)$. Raising inequality \eqref{eq:interpolation for vh(3)} to the power $\frac{1}{q}$, we obtain inequality \eqref{eq:interpolation for Besov functions}. We obtain \eqref{eq:upper directional interpolation for Besov functions with weight} and \eqref{eq:lower directional interpolation for Besov functions with weight} by replacing $h$ with $th$ in \eqref{eq:interpolation for vh(2)} and taking the upper limit (lower limit) as $t \to 0$.
\end{proof}

\begin{corollary}
\label{cor:BMO-interpolation for Besov functions in case p=1/s}
Let $1 \leq p < q < \infty$. If $u\in BMO(\R^N)$, then 
\begin{equation}
\label{eq:interpolation for Besov functions in BMO in case p=1/s}
\|u\|_{B^{1/q,q}(\R^N)}
\leq C[u]_{B^{1/p,p}_w(\R^N)}^{\frac{p}{q}} \|u\|_{BMO(\R^N)}^{1-\frac{p}{q}}\leq C\|u\|_{B^{1/p,p}(\R^N)}^{\frac{p}{q}} \|u\|_{BMO(\R^N)}^{1-\frac{p}{q}},
\end{equation}
where $C=C(p,q,N)$ is a constant depending only on $p$, $q$, and $N$.
\end{corollary}

\begin{proof}
The corollary follows from the first part of Theorem \ref{thm:BMO-interpolation for Besov functions} by choosing $s = \frac{1}{p}$.
\end{proof}

\subsection{Besov Functions $B^{s,p}_w$ with Vanishing Mean Oscillation(VMO)}
In this subsection, we prove for $u\in B_w^{s,p}\cap VMO$ that the limit of the expression $\int_{\R^N}\frac{|u(x+h)-u(x)|^q}{|h|^{sp}} \,dx$ as $h\to 0$ equals $0$ whenever $q>p$; see Theorem \ref{thm: Besov constants vanish for VMO}. In Remark \ref{rem: the limit of Besov constant in case q=p is not necessarily zero}, we provide an example showing that this result is optimal in the sense that it is no longer true for $p=q$.

Observe that if $u\in BMO(\mathbb{R}^N)\cap L^\beta_w(\mathbb{R}^N)$ for some $0<\beta<\infty$, then Theorem~\ref{thm:the minimum theorem for Lorentz quasi-norms} implies that $u\in L^p(\mathbb{R}^N)$ for every $p>\beta$. Consequently, if $\eta\in L^1(\mathbb{R}^N)$, then $u*\eta\in L^p(\mathbb{R}^N)$ for every $p>\max\{\beta,1\}$ because 
$\|u*\eta\|_{L^p(\R^N)}\leq \|\eta\|_{L^1(\R^N)}\|u\|_{L^p(\R^N)}$.
\begin{lemma}[Approximation of $VMO$-functions by mollification]
\label{lem:convergence of mollification in BMO}
Let $0<\beta<\infty$. Assume $u\in VMO(\R^N)\cap L_w^\beta(\R^N)$. Let $\eta\in L^1(\R^N)$ be such that $\int_{\R^N}\eta(z)dz=1$, and for every $\e\in (0,\infty)$ we define the mollification
\begin{equation}
u_\e(x):=\int_{\R^N}u(y)\eta_\e(x-y)dy,\quad \eta_\e(y):=\frac{1}{\e^N}\eta\left(\frac{y}{\e}\right).
\end{equation}
Then,
\begin{equation}
\label{eq:convergence of mollification of VMO function in BMO norm}
\lim_{\e\to 0^+}\|u_\e-u\|_{BMO(\R^N)}=0.
\end{equation}
\end{lemma}

\begin{proof}
Let $\e,R\in(0,\infty)$. By definition of the $BMO$-semi-norm
\begin{multline}
\label{eq:BMO norm for mollified u}
\|u_\e-u\|_{BMO(\R^N)}=\sup_{C}\fint_C\fint_C|u_\e(x)-u(x)-[u_\e(z)-u(z)]|dxdz
\\
\leq \sup_{C,\diam(C)\leq R}\fint_C\fint_C|u_\e(x)-u(x)-[u_\e(z)-u(z)]|dxdz
\\
+\sup_{C,\diam(C)>R}\fint_C\fint_C|u_\e(x)-u(x)-[u_\e(z)-u(z)]|dxdz.
\end{multline}
By the triangle inequality we have
\begin{equation}
\label{eq:on sup where diameter is bigger than R}
\sup_{C,\diam(C)>R}\fint_C\fint_C|u_\e(x)-u(x)-[u_\e(z)-u(z)]|dxdz
\leq 2\sup_{C,\diam(C)>R}\fint_C|u_\e(x)-u(x)|dx.
\end{equation}
By the triangle inequality we get
\begin{multline}
\label{eq:estimate for u-epsilon minus u}
\sup_{C,\diam(C)\leq R}\fint_C\fint_C|u_\e(x)-u(x)-[u_\e(z)-u(z)]|dxdz
\\
=\sup_{C,\diam(C)\leq R}\fint_C\fint_C|u_\e(x)-u_\e(z)-[u(x)-u(z)]|dxdz
\\
\leq \sup_{C,\diam(C)\leq R}\fint_C\fint_C|u_\e(x)-u_\e(z)|dxdz+\sup_{C,\diam(C)\leq R}\fint_C\fint_C|u(x)-u(z)|dxdz.
\end{multline}
By the definition of $u_\e$, Fubini's theorem and the change of variable formula, we obtain for a cube $C_0$ with diameter less or equal to $R$
\begin{multline}
\label{eq:estimate for double integral with convolution}
\fint_{C_0}\fint_{C_0}|u_\e(x)-u_\e(z)|dxdz
=\fint_{C_0}\fint_{C_0}\left|\int_{\R^N}(u(x-y)-u(z-y))\eta_\e(y)dy\right|dxdz
\\
\leq \int_{\R^N}|\eta_\e(y)|\left(\fint_{C_0}\fint_{C_0}\big|u(x-y)-u(z-y)\big|dxdz\right)dy
\\
=\int_{\R^N}|\eta_\e(y)|\left(\fint_{C_0-y}\fint_{C_0-y}\big|u(x)-u(z)\big|dxdz\right)dy
\\
\leq \left(\int_{\R^N}|\eta_\e(y)|dy\right)
\left(\sup_{y\in\R^N}\fint_{C_0-y}\fint_{C_0-y}\big|u(x)-u(z)\big|dxdz\right)
\\
\leq \left(\int_{\R^N}|\eta(y)|dy\right)\sup_{C,\diam(C)\leq R}\fint_C\fint_C|u(x)-u(z)|dxdz.
\end{multline}
In the last inequality of \eqref{eq:estimate for double integral with convolution} we used that the translation of a cube, $C_0-y$, is still a cube with the same diameter.
Therefore, by \eqref{eq:estimate for u-epsilon minus u} and \eqref{eq:estimate for double integral with convolution}
\begin{multline}
\label{eq:estimate for double integral with convolution(1)}
\sup_{C,\diam(C)\leq R}\fint_C\fint_C|u_\e(x)-u(x)-[u_\e(z)-u(z)]|dxdz
\\
\leq \left(1+\|\eta\|_{L^1(\R^N)}\right)\sup_{C,\diam(C)\leq R}\fint_C\fint_C|u(x)-u(z)|dxdz.
\end{multline}
Therefore, by \eqref{eq:BMO norm for mollified u}, \eqref{eq:on sup where diameter is bigger than R} and \eqref{eq:estimate for double integral with convolution(1)} we get the estimate
\begin{multline}
\label{eq:estimate for u-epsilon minus u, final}
\|u_\e-u\|_{BMO(\R^N)}\leq  2\sup_{C,\diam(C)>R}\fint_C|u_\e(x)-u(x)|dx
\\
+\left(1+\|\eta\|_{L^1(\R^N)}\right)\sup_{C,\diam(C)\leq R}\fint_C\fint_C|u(x)-u(z)|dxdz.
\end{multline}
Recall that Lebesgue measure of a cube of whose edge length equals $r$ is $r^N$ and its diameter is $r\sqrt{N}$. Therefore, if $C$ is a cube,  then $\mathcal{L}^N(C)=N^{-N/2}\diam(C)^N$.

It it enough to assume that $u\in L^\beta(\R^N)$ for $\beta\geq 1$, because by Theorem \ref{thm:interpolation for BMO1222,intro}, from $u\in L^\beta_w(\R^N)\cap VMO(\R^N)\subset L^\beta_w(\R^N)\cap BMO(\R^N)$, we get that $u\in L^r(\R^N)$ for every $r\in (\beta,\infty)$.
Thus, since $\beta\geq 1$, we get
\begin{multline}
\label{eq:on sup where diameter is bigger than R(1)}
\sup_{C,\diam(C)>R}\fint_C|u_\e(x)-u(x)|dx
\leq \sup_{C,\diam(C)>R}\left(\fint_C|u_\e(x)-u(x)|^\beta dx\right)^{1/\beta}
\\
\leq \left(\frac{N^{N/2}}{R^N}\right)^{1/\beta}\left(\int_{\R^N}|u_\e(x)-u(x)|^\beta dx\right)^{1/\beta}.
\end{multline}
Since $u_\e$ converges to $u$ in $L^\beta(\R^N)$ 
, we obtain by  \eqref{eq:on sup where diameter is bigger than R(1)}
\begin{equation}
\label{eq:converges to zero on cubes with bound from below}
\lim_{\e\to 0^+}\left(\sup_{C,\diam(C)>R}\fint_C|u_\e(x)-u(x)|dx\right)=0.
\end{equation}
Taking the upper limit as $\e\to 0^+$ in \eqref{eq:estimate for u-epsilon minus u, final} and using \eqref{eq:converges to zero on cubes with bound from below} we get
\begin{equation}
\label{eq:bound on the limsup of the BMO norm of u epsilon minus u}
\limsup_{\e\to 0^+}\|u_\e-u\|_{BMO(\R^N)}\leq \left(1+\|\eta\|_{L^1(\R^N)}\right) \sup_{C,\diam(C)\leq R}\fint_C\fint_C|u(x)-u(z)|dxdz.
\end{equation}
Since $u\in VMO(\R^N)$, we get \eqref{eq:convergence of mollification of VMO function in BMO norm} by taking the limit as $R\to 0^+$ in \eqref{eq:bound on the limsup of the BMO norm of u epsilon minus u}.
\end{proof}

\begin{remark}
\label{rem:convolution is uniformly continuous}
Recall that, for $1\leq p<\infty$, $u\in L^p(\R^N)$, $\e\in(0,\infty)$, and $\eta\in C^1_c(\R^N)$, the convolution $u_\e$ has bounded derivative: by H{\"o}lder inequality, for $q\in [1,\infty]$ such that $\frac{1}{p}+\frac{1}{q}=1$, we have
\begin{equation}
Du_\e(x)=\int_{\R^N}u(y)D\eta_\e(x-y)dy\,\,\Longrightarrow\,\, |Du_\e(x)|\leq\|u\|_{L^p(\R^N)}\|D\eta_\e\|_{L^q(\R^N)}.
\end{equation}
Therefore, $u_\e$ is Lipschitz continuous function in $\R^N$ and particularly uniformly continuous in $\R^N$. In the next lemma we will make use of this property of the convolution.
\end{remark}

\begin{lemma}[$BMO$-Continuity of translations for $VMO$-functions]
\label{lem:continuity of translations in BMO norm for UC functions}
Let $u\in VMO(\R^N)\cap L^\beta_w(\R^N)$ for some $0<\beta<\infty$. Then
\begin{equation}
\label{eq:continuity of translations in BMO norm for UC functions}
\lim_{h\to 0}\|u_h-u\|_{BMO(\R^N)}=0,\quad u_h(x):=u(x+h).
\end{equation}
\end{lemma}

\begin{proof}
We first prove that \eqref{eq:continuity of translations in BMO norm for UC functions} holds for every function $u$ which is uniformly continuous in $\R^N$. 
By definition of the $BMO$-semi-norm
\begin{equation}
\label{eq:BMO norm for u minus translation }
\|u_h-u\|_{BMO(\R^N)}=\sup_{C}\fint_C\fint_C|u(x+h)-u(x)-[u(z+h)-u(z)]|dxdz,
\end{equation}
where the supremum is taken over all  cubes $C\subset \R^N$. Since $u$ is uniformly continuous in $\R^N$, then for arbitrary $\e\in (0,\infty)$ there exists $\delta$ such that, for every $x,z\in\R^N$, if $|x-z|<\delta$, then $|u(x)-u(z)|<\e$. Therefore, we get by \eqref{eq:BMO norm for u minus translation } for every $h$ such that $|h|<\delta$
\begin{equation}
\label{eq:BMO norm for u minus translation (1)}
\|u_h-u\|_{BMO(\R^N)}\leq 2\sup_{C}\fint_C|u(x+h)-u(x)|dx\leq 2\e.
\end{equation}
Therefore,
\begin{equation}
\label{eq:BMO norm for u minus translation (2)}
\limsup_{h\to 0}\|u_h-u\|_{BMO(\R^N)}\leq 2\e.
\end{equation}
It completes the proof of the lemma in case the function $u$ is uniformly continuous. We want to use Remark \ref{rem:convolution is uniformly continuous}, so we assume that $u\in L^\beta(\R^N)$ with $\beta\geq 1$ without loss of generality, because by Theorem \ref{thm:interpolation for BMO1222,intro} , we get from $u\in BMO(\R^N)\cap L^\beta_w(\R^N)$ that $u\in L^r(\R^N)$ for every $r> \beta$. By Lemma \ref{lem:convergence of mollification in BMO} and Remark \ref{rem:convolution is uniformly continuous}, for arbitrary number $\delta\in (0,\infty)$ there exists a function $g$ which is uniformly continuous and $\|u-g\|_{BMO(\R^N)}<\delta$. Therefore,
\begin{equation}
\label{eq:triangle inequality for u-uh}
\|u_h-u\|_{BMO(\R^N)}\leq \|u_h-g_h\|_{BMO(\R^N)}+\|g_h-g\|_{BMO(\R^N)}+\|u-g\|_{BMO(\R^N)}.
\end{equation}
Note that $\|u_h-g_h\|_{BMO(\R^N)}=\|u-g\|_{BMO(\R^N)}$:
\begin{multline}
\|u_h-g_h\|_{BMO(\R^N)}=\sup_{C}\fint_C\fint_C|u(x+h)-g(x+h)-(u(y+h)-g(y+h))|dxdy
\\
=\sup_{C}\fint_{C+h}\fint_{C+h}|u(x)-g(x)-(u(y)-g(y))|dxdy=\|u-g\|_{BMO(\R^N)}.
\end{multline}
Since $g$ is uniformly continuous, we get from \eqref{eq:triangle inequality for u-uh} and the previous case
\begin{equation}
\limsup_{h\to 0}\|u_h-u\|_{BMO(\R^N)}\leq 2 \delta.
\end{equation}
Since $\delta$ can be arbitrarily small, we obtain \eqref{eq:continuity of translations in BMO norm for UC functions}.
\end{proof}

The following theorem contains the proof of Theorems \ref{thm: Besov constants vanish for VMO,intro} and \ref{thm: Besov constants vanish for VMO,intro111}.
\begin{theorem}
\label{thm: Besov constants vanish for VMO}
Let $0<p<q<\infty$ and $s\in(0,1]$. Let $u\in B^{s,p}_w(\R^N)\cap VMO(\R^N)$. Then
\begin{equation}
\label{eq: Besov constants vanish for VMO }
\lim_{h\to 0}\int_{\R^N}\frac{|u(x+h)-u(x)|^q}{|h|^{sp}}dx=0.
\end{equation}
\end{theorem}

\begin{proof}
Let $h\in \R^N\setminus\{0\}$ and let us define $v(x):=u(x+h)-u(x)$. Since $u\in L^p_w(\R^N)\cap BMO(\R^N)$, then, by Proposition \ref{prop:interpolation for translations in BMO - weak case}, $u,v\in L^q(\R^N)$ and
\begin{equation}
\label{eq:BMO inequality for v}
\int_{\R^N} |u(x+h) - u(x)|^q \, dx
\leq C\|u_h-u\|_{BMO(\R^N)}^{q-p}[u(\cdot+h) - u]_{L^p_w(\R^N)}^{p},
\end{equation}
where $C=C(p,q,N)$ is a constant depending only on $p,q,N$, and $u_h(x):=u(x+h)$. Multiplying both sides of \eqref{eq:BMO inequality for v} by $|h|^{-sp}$ we obtain
\begin{multline}
\label{eq:BMO inequality for uh-u(1)}
\int_{\R^N}\frac{|u(x+h)-u(x)|^q}{|h|^{sp}} dx
\leq C\left(\frac{[u(\cdot+h) - u]_{L^p_w(\R^N)}^{p}}{|h|^{sp}}\right)\|u_h-u\|^{q-p}_{BMO(\R^N)}
\\
\leq C[u]^p_{B^{s,p}_w(\R^N)}\|u_h-u\|^{q-p}_{BMO(\R^N)}.
\end{multline}
Since $u\in L^p_w(\R^N)\cap VMO(\R^N)$, we get by and Lemma \ref{lem:continuity of translations in BMO norm for UC functions}
\begin{equation}
\lim_{h\to 0}\|u_h-u\|_{BMO(\R^N)}=0.
\end{equation}
Therefore, we get \eqref{eq: Besov constants vanish for VMO } by taking the limit as $h\to 0$ in \eqref{eq:BMO inequality for uh-u(1)}.
\end{proof}

\begin{remark}
\label{rem: the limit of Besov constant in case q=p is not necessarily zero}
The result of Theorem \ref{thm: Besov constants vanish for VMO} is not true in general when $q=p$.

In \cite{PAjumps} (see Remark 3.5 in \cite{PAjumps}), it was proven that there exists a function $f\in C^{0,\frac{1}{q}}_c(\R^N)\subset B^{1/q,q}(\R^N)$, $q\in (1,\infty)$, such that
\begin{equation}
\label{eq:Besov constant is bounded from below by positive number}
\limsup_{\e \to 0^+}\int_{\R^N}\left(\fint_{B_\e(x)}\frac{|f(x)-f(y)|^q}{|x-y|}dy\right)dx>0.
\end{equation}
Here $C^{0,\frac{1}{q}}_c(\R^N)$ is the space of H{\"o}lder continuous functions with compact support and exponent $1/q$.
Note that since $f\in C^{0,\frac{1}{q}}(\R^N)$, then $f\in VMO(\R^N)$. Indeed, since $f$ is continuous with compact support it is bounded, hence it lies in the space $BMO(\R^N)$. Since there exists a constant $A$ such that $|f(x)-f(z)|\leq A|x-z|^{1/q}$, for every $x,z\in\R^N$, then we get
\begin{multline}
\lim_{R\to 0^+}\left(\sup_{C,\diam(C)\leq R}\fint_C\fint_C|f(x)-f(z)|dxdz\right)
\leq A\lim_{R\to 0^+}\left(\sup_{C,\diam(C)\leq R}\fint_C\fint_C |x-z|^{1/q}dxdz\right)
\\
\leq A\lim_{R\to 0^+}R^{1/q}=0.
\end{multline}
It proves that $f\in VMO(\R^N)$.
Note also that
\begin{multline}
\label{eq:Besov constant is bounded from below by positive number(1)}
\int_{\R^N}\left(\fint_{B_\e(x)}\frac{|f(x)-f(y)|^q}{|x-y|}dy\right)dx=\int_{\R^N}\left(\fint_{B_\e(0)}\frac{|f(x)-f(x+y)|^q}{|y|}dy\right)dx
\\
=\fint_{B_\e(0)}\left(\int_{\R^N}\frac{|f(x)-f(x+y)|^q}{|y|}dx\right)dy
\leq \sup_{|y|\leq \e}\left(\int_{\R^N}\frac{|f(x)-f(x+y)|^q}{|y|}dx\right).
\end{multline}
Therefore, by \eqref{eq:Besov constant is bounded from below by positive number} and \eqref{eq:Besov constant is bounded from below by positive number(1)}
\begin{equation}
\limsup_{y\to 0}\int_{\R^N}\frac{|f(x+y)-f(x)|^q}{|y|}dx=\lim_{\e \to 0^+}\left(\sup_{|y|\leq \e}\int_{\R^N}\frac{|f(x+y)-f(x)|^q}{|y|}dx\right)>0.
\end{equation}
It shows that the result of Theorem \ref{thm: Besov constants vanish for VMO} is not true in case $1<p=q<\infty$ and $s=\frac{1}{q}$.
\end{remark}

\subsection{$BMO$-interpolation in the space of weak functions with bounded variation $BV_w$}

The proof of Theorem \ref{cor:BMO-interpolation for BV,intro} is contained in the following corollary.
\begin{corollary}[$BMO$-interpolation for $BV$-functions]
\label{cor:BMO-interpolation for BV}
Let $1 < q < \infty$. If $u\in BMO(\R^N)$, then 
\begin{equation}
\label{eq:interpolation for BV functions in BMO}
\|u\|_{B^{1/q,q}(\R^N)}
\leq C [u]^{\frac{1}{q}}_{BV_w(\R^N)} \|u\|^{1-\frac{1}{q}}_{BMO(\R^N)}\leq  C\|u\|^{\frac{1}{q}}_{BV(\R^N)} \|u\|^{1-\frac{1}{q}}_{BMO(\R^N)},
\end{equation}
where $C=C(q,N)$ is a constant depending only on $q$ and $N$.
\end{corollary}

\begin{proof}
Choosing $s = p = 1$ in the first part  of Theorem \ref{thm:BMO-interpolation for Besov functions}, we get \eqref{eq:interpolation for BV functions in BMO}.
\end{proof}

\section{Appendix}
\label{sec:appendix}  
    
\begin{proposition}
\label{prop:BMO-function is integrable of any order on a cube}
Let $C\subset \R^N$ be a cube. Assume that $u\in BMO(C)$ and let $p\in (0,\infty),\gamma\in (0,\infty]$. Then, $u\in L^{p,\gamma}(C)$. 
\end{proposition}

\begin{proof}
By the additivity of the integral, we obtain for $\gamma < \infty$ and every number $\delta \in (0,\infty)$ 
\begin{multline}
\label{eq:estimate for Lp norm in a cube for positive p}
\|u-u_C\|^\gamma_{L^{p,\gamma}(C)} = \int_0^\infty \mathcal{L}^N \left(\Set{x \in C}[|u(x)-u_C|^p > t]\right)^{\frac{\gamma}{p}}t^{\frac{\gamma}{p}-1} dt
\\
= \int_0^\delta \mathcal{L}^N \left(\Set{x \in C}[|u(x)-u_C|^p > t]\right)^{\frac{\gamma}{p}}t^{\frac{\gamma}{p}-1} dt
+ \int_\delta^\infty \mathcal{L}^N \left(\Set{x \in C}[|u(x)-u_C|^p > t]\right)^{\frac{\gamma}{p}}t^{\frac{\gamma}{p}-1} dt
\\
\leq \delta^{\frac{\gamma}{p}} \frac{p}{\gamma}\mathcal{L}^N(C)^{\frac{\gamma}{p}}
+ \int_\delta^\infty \mathcal{L}^N \left(\Set{x \in C}[|u(x)-u_C| > t^{1/p}]\right)^{\frac{\gamma}{p}}t^{\frac{\gamma}{p}-1} dt.
\end{multline}
Let $a, A, \alpha$ be as in Theorem \ref{thm:BMO inequality}. Assume that $\|u\|_{BMO(C)} > 0$; otherwise, if $\|u\|_{BMO(C)}=0$, then the function $u$ is constant on $C$, a set of finite Lebesgue measure, and hence $u \in L^{p,\gamma}(C)$. We choose $\delta := \left(a\textbf{k}\right)^p$, where $\textbf{k} := \|u\|_{BMO(C)}$. Then, for every $t \in (\delta,\infty)$, we have $\frac{t^{1/p}}{\textbf{k}} \geq a$. Therefore, from \eqref{eq:estimate for Lp norm in a cube for positive p} and Theorem \ref{thm:BMO inequality}, we obtain
\begin{multline}
\label{eq:estimate for Lp norm in a cube for positive p(1)}
\|u-u_C\|^\gamma_{L^{p,\gamma}(C)} \leq \delta^{\frac{\gamma}{p}} \frac{p}{\gamma}\mathcal{L}^N(C)^{\frac{\gamma}{p}}
+ \int_\delta^\infty \left[\frac{A}{\textbf{k}} e^{-\frac{\alpha}{\textbf{k}} t^{1/p}} \int_C |u(x)-u_C| dx \right]^{\frac{\gamma}{p}}t^{\frac{\gamma}{p}-1}  dt
\\
= \delta^{\frac{\gamma}{p}} \frac{p}{\gamma}\mathcal{L}^N(C)^{\frac{\gamma}{p}}
+ \left[\frac{A}{\textbf{k}} \int_C |u(x)-u_C| dx\right]^{\frac{\gamma}{p}} \int_{\delta}^\infty e^{-\frac{\alpha}{\textbf{k}} \frac{\gamma}{p}t^{1/p}}t^{\frac{\gamma}{p}-1}  dt < \infty.
\end{multline}
From \eqref{eq:estimate for Lp norm in a cube for positive p(1)}, we see that $u - u_C \in L^{p,\gamma}(C)$. Since $C$ is a bounded set, it follows that $u = (u - u_C) + u_C \in L^{p,\gamma}(C)$. The case $\gamma=\infty$ is proved in a similar manner.
\end{proof}

\subsection{An alternative proof for $BMO$-interpolation in Lorentz spaces via real interpolation and sharp maximal function}
\begin{definition}[$K$-functional for quasi-normed spaces]
\label{subsec: An alternative proof for $BMO$-interpolation in Lorentz spaces via real interpolation and sharp maximal function}
Let $(X_0,\|\cdot\|_{X_0})$ and $(X_1,\|\cdot\|_{X_1})$ be quasi-normed spaces continuously embedded into a common topological Hausdorff vector space. Define the sum space
\[
X_0 + X_1 := \{ f = f_0 + f_1 : f_0 \in X_0,\; f_1 \in X_1 \}.
\]
For $f \in X_0 + X_1$ and $t>0$, the $K$-functional of $f$ with respect to the couple $(X_0,X_1)$ is defined by
\[
K(t,f; X_0, X_1)
:= \inf \left\{ \|f_0\|_{X_0} + t\,\|f_1\|_{X_1} : f = f_0 + f_1,\; f_0 \in X_0,\; f_1 \in X_1 \right\}.
\]
\end{definition}

\begin{definition}[Real interpolation space]
\label{def:Real interpolation space}
Let $X_0$ and $X_1$ be quasi-normed spaces continuously embedded in a common topological Hausdorff vector space. Let $0<\theta<1$ and $0<\gamma \le \infty$. For $f \in X_0 + X_1$, define
\[
\|f\|_{(X_0,X_1)_{\theta,\gamma}} :=
\begin{cases}
\left( \displaystyle \int_0^\infty \bigl( t^{-\theta} K(t,f;X_0,X_1) \bigr)^\gamma \, \frac{dt}{t} \right)^{1/\gamma}, & 0<\gamma<\infty,\\[1ex]
\displaystyle \sup_{0<t<\infty} t^{-\theta} K(t,f;X_0,X_1), & \gamma = \infty,
\end{cases}
\]
and
\[
(X_0,X_1)_{\theta,\gamma} := \{ f \in X_0 + X_1 : \|f\|_{(X_0,X_1)_{\theta,\gamma}} < \infty \}.
\]
\end{definition}
The spaces $(X_0,X_1)_{\theta,\gamma}$ are intermediate spaces between the intersection space $X_0 \cap X_1$, endowed with the quasi-norm $\|f\|_{X_0} + \|f\|_{X_1}$, and the sum space $X_0 + X_1$, endowed with the quasi-norm $\inf_{f=f_0+f_1}\left(\|f_0\|_{X_0} + \|f_1\|_{X_1}\right)$. More precisely,
\begin{equation}
X_0 \cap X_1 \hookrightarrow (X_0,X_1)_{\theta,\gamma} \hookrightarrow X_0 + X_1.
\end{equation}
The order of the couple matters; however, one has the symmetry $(X_0,X_1)_{\theta,\gamma} = (X_1,X_0)_{1-\theta,\gamma}$.
\begin{definition}[Sharp maximal function]
Let $f:\R^N\to \R$ be a Lebesgue measurable function. For \(0 < \alpha \leq \frac{1}{2}\), define
\begin{equation}
\label{eq:local-sharp-max}
M_{\alpha} f(x) := \sup\limits_{Q \owns x}\inf_{c\in \R} \inf\Bigl\{ A \geq  0 : \mathcal{L}^N\left(\{ y \in Q : |f(y) - c| > A \}\right) < \alpha \mathcal{L}^N(Q) \Bigr\},
\end{equation}
where $Q$ is an arbitrary finite cube in $\R^N$ with sides parallel to the coordinate axes. Note, in \eqref{eq:local-sharp-max} $x$ is fixed and the supremum runs over cubes $Q$ containing $x$.
\end{definition}

The operator $M_\alpha$ is denoted in \cite{JawerthTorchinsky1985} by $M^{\#}_{0,\alpha}$. 

What is its meaning? Fix a point $x\in\mathbb{R}^N$. Let $Q$ be a cube containing $x$. Let us first fix a reference center $c\in\mathbb{R}$ around which we examine the function $f$. Assume that $0<\alpha\leq \tfrac12$ is a small number.

Consider the quantity
\begin{equation}
\label{eq:oscillation dependent of center c}
\inf\Bigl\{ A \ge 0 : \mathcal{L}^N\big(\{ y \in Q : |f(y) - c| > A \}\big) < \alpha \mathcal{L}^N(Q) \Bigr\}.
\end{equation}
Denote this infimum by $A_0$. Then, for every $\varepsilon>0$, we have
\[
\mathcal{L}^N\big(\{ y \in Q : |f(y) - c| > A_0+\varepsilon \}\big) < \alpha \mathcal{L}^N(Q).
\]
This means that, outside a subset of $Q$ of relative measure less than $\alpha$, the function $f$ satisfies
\[
|f(y)-c|\le A_0+\varepsilon.
\]
In other words, on most of the cube $Q$, the function $f$ is concentrated in an interval of length $2(A_0+\varepsilon)$ centered at $c$.

The restriction $\alpha\leq \tfrac12$ is essential. Indeed, let $c_1,c_2\in\mathbb{R}$ be two admissible centers corresponding to some $A$, meaning that
\begin{equation}
\mathcal{L}^N(Q\setminus E_i)<\alpha\mathcal{L}^N(Q), \quad i=1,2,
\end{equation}
where we define
\begin{equation}
E_i := \{y\in Q : |f(y)-c_i|\le A\}, \quad i=1,2.
\end{equation}
Then
\begin{equation}
\mathcal{L}^N(E_i) > (1-\alpha)\mathcal{L}^N(Q) \ge \tfrac12 \mathcal{L}^N(Q), \quad i=1,2.
\end{equation}
Hence $E_1\cap E_2 \neq \emptyset$. For any $y\in E_1\cap E_2$, we have
\[
|f(y)-c_1|\le A, \quad |f(y)-c_2|\le A,
\]
and therefore
\[
|c_1-c_2|\le 2A.
\]
Thus all admissible centers are close to each other, which makes the notion of center stable.

However, the choice of the reference center $c$ is arbitrary and not intrinsic to the function. Therefore, we take the infimum in \eqref{eq:oscillation dependent of center c} over all $c \in \mathbb{R}$. This produces a quantity that is independent of the choice of reference level and invariant under the addition of constants. Roughly speaking, it selects the most suitable center around which $f$ is most concentrated on $Q$.

Note that, so far, the fixed point $x$ has not played any role; we have used only that $f$ is a measurable function on $Q$. In other words, the quantity
\begin{equation}
\label{eq:oscillation dependent on Q}
\inf_{c\in\mathbb{R}} \inf\Bigl\{ A \ge 0 : \mathcal{L}^N\big(\{ y \in Q : |f(y) - c| > A \}\big) < \alpha \mathcal{L}^N(Q) \Bigr\}
\end{equation}
describes the oscillation of a function on a given cube $Q$ for the parameter $\alpha$.

Finally, taking the supremum over all cubes $Q$ containing $x$, we obtain $M_\alpha f(x)$, which measures the largest oscillation of $f$ on cubes containing $x$. Thus, $M_\alpha f(x)$ captures the worst oscillatory behaviour of $f$ around $x$ across all scales. A larger value of $M_\alpha f(x)$ indicates stronger oscillation on some cube containing $x$, whereas a smaller value indicates that $f$ is closer to being nearly constant on every such cube.

\begin{definition}[Distribution function and decreasing rearrangement]
Let $(X,\mu)$ be a measure space and let $f:X\to \mathbb{R}$ be a measurable function.  
The {\it distribution function} of $f$ is defined by
\[
d_f(h):=\mu(\{x\in X:|f(x)|>h\}), \qquad h>0.
\]
The {\it decreasing rearrangement} of $f$ is the function $f^*:(0,\infty)\to[0,\infty]$ defined by
\[
f^*(s):=\inf\{h>0:d_f(h)\le s\}.
\]
\end{definition}

The function $d_f(h)$ measures the size of the level set where $|f|$ exceeds $h$, and it is decreasing in $h$.

The decreasing rearrangement $f^*(s)$ is equivalently characterized by
\[
f^*(s)=\inf\{\alpha\ge 0:\mu(\{|f|>\alpha\})\le s\}.
\]

Thus, $f^*(s)=\alpha$ means that $\alpha$ is the smallest level such that the set where $|f|>\alpha$ has measure at most $s$, i.e.
\[
\mu(\{x\in X:|f(x)|>\alpha\})\le s.
\]

In this sense, $f^*$ provides a “sorted” version of $|f|$, arranged in decreasing order of magnitude with respect to the measure of the underlying space.

\begin{theorem}[Corollary 3.3 in \cite{JawerthTorchinsky1985}]
\label{thm:equivalence between K of Lorentz and BMO and M}
Let $0 < p < \infty$. Then for $t>0$ and 
$f \in L^{p}_w(\mathbb{R}^N) + BMO(\mathbb{R}^N)$, we have that
\begin{equation}
\label{eq:K equi to M}
K\!\left(t,f; L^{p}_w, BMO \right) \approx 
\sup_{0 < s < t^p} \; s^{1/p} \bigl( M_{\alpha}f \bigr)^{*}(s)
\end{equation}
provided that \(0 < \alpha \leq  \alpha_0 = \alpha_0(N)\) is sufficiently small. The constants of equivalence do not depend on $t$ or $f$.
\end{theorem}

\begin{theorem}[Theorem 4.2 in \cite{Holmstedt1970}]
\label{thm:equi bet K and truncated norm}
Let $0<p<\infty$ and $f\in L^{p}_w(\R^N)+L^\infty(\R^N)$. Then, for every $t>0$
\begin{equation}
K\!\left(t,f; L^{p}_w, L^\infty \right)\approx \sup_{0 < s < t^p} \; s^{1/p}f^{*}(s).
\end{equation}
The constants of equivalence do not depend on $t$ or $f$.
\end{theorem}

\begin{theorem}[Lemma 4.2 in \cite{Holmstedt1970}]
\label{thm:equi bet quasi Lorentz norms}
For $0<p<q<\infty$, $0<\gamma\leq \infty$, and $f\in L^p_w(\R^N)+L^\infty(\R^N)$
\begin{equation}
\|f\|_{(L^p_w,L^\infty)_{\theta,\gamma}}\approx \|f\|_{L^{q,\gamma}},\quad \theta:=1-\frac{p}{q}.
\end{equation}
\end{theorem}
In Lemma 4.2 of \cite{Holmstedt1970}, Theorem \ref{thm:equi bet quasi Lorentz norms} is stated for $\gamma<\infty$, but the result remains valid also for $\gamma=\infty$. This follows from Theorem \ref{thm:equi bet K and truncated norm}. We now prove it.

\begin{proposition}
For $0<p<q<\infty$ and $f\in L^p_w(\R^N)+L^\infty(\R^N)$, we have
\begin{equation}
\|f\|_{(L^p_w,L^\infty)_{\theta,\infty}}\approx \|f\|_{L^{q,\infty}},
\qquad \theta:=1-\frac{p}{q}.
\end{equation}
\end{proposition}

\begin{proof}
Using Theorem \ref{thm:equi bet K and truncated norm}, we get
\begin{multline}
\|f\|_{(L^p_w,L^\infty)_{\theta,\infty}}
:=\sup_{0<t<\infty}t^{-\theta}K\left(t,f; L^{p}_w, L^\infty \right)\approx \sup_{0<t<\infty}\sup_{0 < s < t^p} \;t^{-\theta} s^{1/p}f^{*}(s)
=\sup_{0<t<\infty}\sup_{0 < \xi < t} \;t^{-\theta} \xi f^{*}(\xi^p)
\\
=\sup_{0<\xi<\infty}\sup_{\xi<t<\infty} \;t^{-\theta} \xi f^{*}(\xi^p)
=\sup_{0<\xi<\infty}\;\xi^{1-\theta} f^{*}(\xi^p)
=\sup_{0<\xi<\infty}\;\left(\xi^{p}\right)^{1/q} f^{*}(\xi^p)
=\sup_{0<t<\infty}\;t^{1/q} f^{*}(t)
\approx \|f\|_{L^{q,\infty}}.
\end{multline}
\end{proof}
By Theorem \ref{thm:equivalence between K of Lorentz and BMO and M} and Theorem \ref{thm:equi bet K and truncated norm}
\begin{equation}
\label{eq:K equi to M1}
K\!\left(t,u; L^{p}_w, BMO \right) \approx 
\sup_{0 < s < t^p} \; s^{1/p} \bigl( M_{\alpha} u \bigr)^{*}(s)\approx K\!\left(t,M_{\alpha}u; L^{p}_w, L^\infty \right)
\end{equation}
provided $u\in L^{p}_w(\R^N)+BMO(\R^N)$ and $M_{\alpha}u\in L^{p}_w(\R^N)+L^\infty(\R^N)$.
Therefore, from Definition \ref{def:Real interpolation space}, and \eqref{eq:K equi to M1}, we obtain 
\begin{equation}
\label{eq:part of the equiv}
\|u\|_{(L^p_w,BMO)_{\theta,\gamma}}\approx \|M_{\alpha}u\|_{(L^p_w,L^\infty)_{\theta,\gamma}},\quad \forall \theta\in (0,1),\,\ p\in (0,\infty),\,\gamma\in (0,\infty].
\end{equation}
Using \eqref{eq:part of the equiv} and Theorem \ref{thm:equi bet quasi Lorentz norms}, we get for $0<p<q<\infty$ and $0<\gamma\leq \infty$
\begin{equation}
\label{eq:full of the equiv}
\|u\|_{(L^p_w,BMO)_{\theta,\gamma}}\approx \|M_{\alpha}u\|_{(L^p_w,L^\infty)_{\theta,\gamma}}\approx \|M_{\alpha}u\|_{L^{q,\gamma}},\quad \theta:=1-\frac{p}{q}.
\end{equation}

\begin{proposition}[Corollary 2.5 in \cite{JawerthTorchinsky1985}]
\label{prop:equivalence of M-BMO and L-infty}
There exists $\alpha_0 = \alpha_0(N)$ such that for $0<\alpha\leq \alpha_0$, there exist constants $c_j = c_j(\alpha,N)$, $j=1,2$, such that for every measurable function $u$ we have
\begin{equation}
c_1 \|u\|_{BMO(\R^N)} \le \|M_{\alpha}u\|_{L^\infty(\R^N)} \le c_2 \|u\|_{BMO(\R^N)}.
\end{equation}
\end{proposition}

Therefore, we get from \eqref{eq:full of the equiv} the following corollary:
\begin{corollary}
\label{cor:full of the equiv2}
Let $0<p<q<\infty$, $0<\gamma\leq \infty$, $u\in L^p_w(\R^N)\cap BMO(\R^N)$. Then, there exists a constant $C$ such that 
\begin{equation}
\label{eq:full of the equiv2}
\|M_{\alpha}u\|_{L^{q,\gamma}}\leq C\|u\|_{(L^p_w,BMO)_{\theta,\gamma}},\quad \theta:=1-\frac{p}{q},
\end{equation}
provided that \(0 < \alpha \leq  \alpha_0 = \alpha_0(N)\) is sufficiently small.
\end{corollary}

\begin{proof}
Since we have \eqref{eq:full of the equiv}, we only need to explain why $u\in L^{p}_w(\R^N)+BMO(\R^N)$ and $M_{\alpha}u\in L^{p}_w(\R^N)+L^\infty(\R^N)$. This it true because
\begin{equation}
u\in L^{p}_w(\R^N)\cap BMO(\R^N)\subset L^{p}_w(\R^N)+BMO(\R^N),
\end{equation}
and by Proposition \ref{prop:equivalence of M-BMO and L-infty}
\begin{equation}
M_{\alpha}u \in L^\infty(\R^N)\subset L^{p}_w(\R^N)+L^\infty(\R^N).
\end{equation}
\end{proof}

\begin{proposition}
\label{prop:interpolation-trivial-general}
Let $X_0,X_1$ be quasi-normed spaces continuously embedded into a common topological Hausdorff vector space. Let $0<\theta<1$, $0<\gamma\leq \infty$, and let $u \in X_0 \cap X_1$. 
Then there exists a constant $C=C(\theta,\gamma)$ such that
\begin{equation}
\label{eq:interpolation-trivial-general}
\|u\|_{(X_0,X_1)_{\theta,\gamma}} \leq C \, \|u\|_{X_0}^{1-\theta} \, \|u\|_{X_1}^{\theta}.
\end{equation}
\end{proposition}

\begin{proof}
Let us denote $K(t):=K(t,u;X_0,X_1)$.
For any $t>0$, consider the trivial decompositions:
$u = u_0 + u_1$, $u_0 = u, \ u_1 = 0$ and $u_0 = 0, u_1 = u$.
Then by definition of the $K$-functional,
\[
K(t) \leq \min\{\|u\|_{X_0}, \ t \, \|u\|_{X_1}\}.
\]

Assume that $\|u\|_{X_0}, \|u\|_{X_1} \in (0,\infty)$; otherwise, if one of them is equal to zero, then $u=0$ in $X_0+X_1$, so $\|u\|_{(X_0,X_1)_{\theta,\gamma}}=0$, therefore, \eqref{eq:interpolation-trivial-general} holds trivially. Let $t_0 := \frac{\|u\|_{X_0}}{\|u\|_{X_1}}$. Split the interpolation integral at $t_0$: for $\gamma<\infty$
\begin{equation}
\int_0^\infty \bigl( t^{-\theta} \, K(t) \bigr)^{\gamma} \, \frac{dt}{t}
=
\int_0^{t_0} \bigl( t^{-\theta} \, K(t) \bigr)^{\gamma} \, \frac{dt}{t}
+
\int_{t_0}^\infty \bigl( t^{-\theta} \, K(t) \bigr)^{\gamma} \, \frac{dt}{t}.
\end{equation}

Since $K(t) \le t \, \|u\|_{X_1}$, then
\begin{equation}
\int_0^{t_0} (t^{-\theta} K(t))^\gamma \frac{dt}{t} 
\le \|u\|_{X_1}^\gamma \int_0^{t_0} t^{(1-\theta)\gamma} \frac{dt}{t} 
= \frac{\|u\|_{X_1}^\gamma}{(1-\theta)\gamma} \, t_0^{(1-\theta)\gamma}
= \frac{1}{(1-\theta)\gamma}\|u\|_{X_0}^{(1-\theta)\gamma} \, \|u\|_{X_1}^{\theta \gamma}.
\end{equation}

Since $K(t) \le \|u\|_{X_0}$, then
\begin{equation}
\int_{t_0}^{\infty} (t^{-\theta} K(t))^\gamma \frac{dt}{t} 
\le \|u\|_{X_0}^\gamma \int_{t_0}^\infty t^{-\theta \gamma} \frac{dt}{t} 
= \frac{\|u\|_{X_0}^\gamma}{\theta \gamma} \, t_0^{-\theta \gamma}
= \frac{1}{\theta\gamma}\|u\|_{X_0}^{(1-\theta)\gamma} \, \|u\|_{X_1}^{\theta \gamma}.
\end{equation}

Combining the above estimates yields \eqref{eq:interpolation-trivial-general}. The case $\gamma=\infty$ is proved in a similar way.
\end{proof}

From \eqref{eq:full of the equiv2} and Proposition \ref{prop:interpolation-trivial-general}, we obtain
\begin{equation}
\label{eq:intepolation bound with lp weak and bmo}
\|M_{\alpha}u\|_{L^{q,\gamma}}\leq C\|u\|_{L^p_w}^{1-\theta}\|u\|^{\theta}_{BMO}.
\end{equation}
For \eqref{eq:intepolation bound with lp weak and bmo}, we regard $\mathrm{BMO}$ as the space of equivalence classes of measurable functions, where two functions $f$ and $g$ are identified if their difference is constant almost everywhere, that is,
\[
f \sim g \quad \Longleftrightarrow \quad f-g = c \ \text{a.e. for some constant } c \in \mathbb{R}.
\]
Thus, elements of $\mathrm{BMO}$ are equivalence classes $[f]$, and the BMO seminorm is well-defined on these classes and becomes a norm on the quotient space $\mathrm{BMO}/\mathbb{R}$.
\begin{theorem}[Corollary 2.6 in \cite{JawerthTorchinsky1985}]
\label{thm:sharp-maximal-constants}
Let $0 < q < \infty$ and $0 < \gamma \leq  \infty$. There exists $\alpha_0 = \alpha_0(N)$ such that for every $0 < \alpha < \alpha_0$, there exist constants $c_1 = c_1(\alpha,N)$ and $c_2 = c_2(\alpha,N)$ with the following property:  

For each measurable function $u$ on $\mathbb{R}^N$, there exists a constant $c_u \in \mathbb{R}$ such that
\[
c_1 \, \|u - c_u\|_{L^{q,\gamma}(\mathbb{R}^N)} \le 
\left\| M_{\alpha} u \right\|_{L^{q,\gamma}(\mathbb{R}^N)} 
\le c_2 \, \|u - c_u\|_{L^{q,\gamma}(\mathbb{R}^N)}.
\]
\end{theorem}

Therefore, by Theorem~\ref{thm:sharp-maximal-constants} and \eqref{eq:intepolation bound with lp weak and bmo}, there exists $c_u \in \mathbb{R}$ with the property
\begin{equation}
\|u - c_u\|_{L^{q,\gamma}} \leq C \|u\|_{L^p_w}^{1-\theta} \|u\|_{BMO}^{\theta}.
\end{equation}

Now, if $u \in L^p_w \cap BMO$, then $\|u - c_u\|_{L^{q,\gamma}} < \infty$. Therefore, $c_u$ must be zero by Proposition \ref{prop:f in lp weak and f+C in lorentz implies C=0}. Hence:

\begin{corollary}
Let $0 < p <q < \infty$ and $0 < \gamma \leq  \infty$. Let $u \in L^p_w(\R^N) \cap BMO(\R^N)$. Then
\begin{equation}
\label{eq:interp result through l.s.m}
\|u\|_{L^{q,\gamma}(\R^N)} \leq C \|u\|_{L^p_w(\R^N)}^{1-\theta} \|u\|_{BMO(\R^N)}^{\theta},\quad \theta:=1-\frac{p}{q}.
\end{equation}

\end{corollary}
Note that the assumption $u \in L^p_w(\mathbb{R}^N) \cap BMO(\mathbb{R}^N)$ is not essential. Indeed, if this condition fails, then at least one of the quantities
$\|u\|_{L^p_w(\mathbb{R}^N)}$ or $\|u\|_{\mathrm{BMO}(\mathbb{R}^N)}$
is infinite, and therefore the right-hand side of \eqref{eq:interp result through l.s.m} is equal to $\infty$.

\begin{proposition}
\label{prop:f in lp weak and f+C in lorentz implies C=0}
Let $(X,\mu)$ be a measure space with $\mu(X)=\infty$. Let $0<p<\infty$, $0<q<\infty$, and $0<\gamma\le\infty$. Suppose that 
$f\in L^{p}_w(X)$ and $f+C \in L^{q,\gamma}(X)$ for some constant $C\in\mathbb{R}$. Then necessarily $C=0$.
\end{proposition}

\begin{proof}
Since $f\in L^{p}_w(X)$, we have
\[
\mu\big(\{x\in X:\ |f(x)|>\lambda\}\big)\leq \frac{\|f\|^p_{L^p_w}}{\lambda^p}<\infty
\quad \text{for all } \lambda>0.
\]
Assume by contradiction that $C\neq 0$, and set $\lambda=\frac{|C|}{2}$. Then
\[
\{x\in X:\ |f(x)|\le \lambda\}\subset \{x\in X:\ |f(x)+C|\ge \lambda\}.
\]
Indeed, if $|f(x)|\le \lambda=\frac{|C|}{2}$, then $|f(x)+C|\ge |C|-|f(x)|\ge \frac{|C|}{2}=\lambda$.
Therefore,
\[
\mu\big(\{|f|\le \lambda\}\big)\le \mu\big(\{|f+C|\ge \lambda\}\big).
\]
Since $\mu\big(\{|f|>\lambda\}\big)<\infty$ and $\mu(X)=\infty$, it follows that $\mu\big(\{|f|\le \lambda\}\big)=\infty$.
Consequently, $\mu\big(\{|f+C|\ge \lambda\}\big)=\infty$.
Now, for every $0<t<\lambda^q$, we have
\[
\{x\in X:\ |f(x)+C|\ge \lambda\}\subset \{x\in X:\ |f(x)+C|^q>t\}.
\]
Therefore,
\begin{multline}
\int_0^\infty \Big(\mu\big(\{x\in X:\ |f(x)+C|^q>t\}\big)\Big)^{\frac{\gamma}{q}} t^{\frac{\gamma}{q}-1}\,dt \\
\ge \int_0^{\lambda^q} \Big(\mu\big(\{x\in X:\ |f(x)+C|^q>t\}\big)\Big)^{\frac{\gamma}{q}} t^{\frac{\gamma}{q}-1}\,dt
\\
\geq \mu\big(\{x\in X:\ |f(x)+C|\geq\lambda\}\big)^{\frac{\gamma}{q}}\frac{q}{\gamma}\lambda^\gamma
= \infty.
\end{multline}
Thus $\|f+C\|_{L^{q,\gamma}(X)}=\infty$, i.e., $f+C\notin L^{q,\gamma}(X)$, a contradiction. Hence $C=0$.
\end{proof}

\subsection{An alternative proof for $BMO$-interpolation in $W^{s,p}_w$ via real interpolation theory and the spaces $\dot {\mathcal B}^{s}_{p}(-sp,\infty)$}
\label{An alternative proof for BMO-interpolation in Wsp via real interpolation theory and the spaces}
Here we use the usual dot notation to distinguish between homogeneous and non-homogeneous spaces: a dot over the symbol indicates the homogeneous case, while the absence of a dot indicates the non-homogeneous case. For example, the homogeneous Besov space is denoted by $\dot{B}^s_{p,q}(\R^N)$, whereas the non-homogeneous Besov space is denoted by $B^s_{p,q}(\R^N)$.

Let $0<s<1$ and $1< p<\infty$. Let $\theta\in (0,1)$ be such that $p_0:=p(1-\theta)\geq 1$. Let $s_0:=\frac{s}{1-\theta}$.
From \cite[Theorem 1.14]{Dominguez2023}, we obtain
\begin{equation}
\label{eq:B with parameters and Besov}
\dot {\mathcal B}^{s}_{p}(-sp,\infty)
=
\left(\dot B^{s_0}_{p_0,p_0},\, \dot B^{0}_{\infty,\infty}\right)_{\theta,\infty}.
\end{equation}
Let $q$ be such that \(1<p<q<\infty\). Then, by \eqref{eq:B with parameters and Besov}, we obtain
\begin{equation}
\label{eq:inter Binfity}
\left( \dot {\mathcal B}^{s}_{p}(-sp,\infty), \dot B^{0}_{\infty,\infty} \right)_{1-\frac{p}{q},\,q}
=
\left(
\left(\dot B^{s_0}_{p_0,p_0},\, \dot B^{0}_{\infty,\infty}\right)_{\theta,\infty},
\dot B^{0}_{\infty,\infty}
\right)_{1-\frac{p}{q},\,q}.
\end{equation}
Recall the following version of the reiteration theorem \cite[Remark 3.2, equation (3.17)]{Holmstedt1970}:
\begin{equation}
\left(
(A_0,A_1)_{\theta,q},\,A_1
\right)_{\lambda,p}
=(A_0,A_1)_{(1-\lambda)\theta+\lambda,\,p}\quad \theta,\lambda\in (0,1),\,p,q\in (0,\infty].
\end{equation}
This reiteration theorem applies to quasi-normed spaces $A_0$ and $A_1$. We regard the homogeneous Besov spaces $\dot B^{s}_{p,q}$ as quotient spaces modulo constants; therefore, they are quasi-normed spaces.

Using this reiteration and  interpolation in Besov spaces \cite[6.4.5 Theorem, equation $(3)$]{BerghLofstrom1976}, we get
\begin{equation}
\label{eq:inter Binfity1}
\left(
\left(\dot B^{s_0}_{p_0,p_0},\, \dot B^{0}_{\infty,\infty}\right)_{\theta,\infty},
\dot B^{0}_{\infty,\infty}
\right)_{1-\frac{p}{q},\,q}
=
\left(\dot B^{s_0}_{p_0,p_0},\, \dot B^{0}_{\infty,\infty}\right)_{\frac{p}{q}\theta+1-\frac{p}{q},\,q}=\dot B^{\frac{ps}{q}}_{q,q}.
\end{equation}
Since we have $\dot B^{\frac{ps}{q}}_{q,q}=\dot{W}^{\frac{ps}{q},\,q}$ \cite[Proposition 14.40]{Leoni2017}, then by \eqref{eq:inter Binfity} and  \eqref{eq:inter Binfity1}, we obtain
\begin{equation}
\label{eq:inter Binfity3}
\left( \dot {\mathcal B}^{s}_{p}(-sp,\infty), \dot B^{0}_{\infty,\infty} \right)_{1-\frac{p}{q},\,q}
=
\dot{W}^{\frac{ps}{q},\,q}.
\end{equation}
In particular, from \eqref{eq:inter Binfity3} and Proposition \ref{prop:interpolation-trivial-general} we obtain
\begin{equation}
\label{eq:weak W and B coincide1}
\|u\|_{\dot{W}^{\frac{ps}{q},q}}
\le C \,
\|u\|_{\dot{\mathcal{B}}^{s}_{p}(-sp,\infty)}^{\frac{p}{q}}
\,
\|u\|_{\dot{{B}}^{0}_{\infty,\infty}}^{1-\frac{p}{q}}.
\end{equation}
Therefore, using \eqref{eq:weak W and B coincide} and \eqref{eq:weak W and B coincide1}, we obtain
\begin{equation}
\label{eq:interpolation for BMO and B0 infty infty}
\|u\|_{\dot{W}^{\frac{ps}{q},q}}
\le C \,
[u]_{W^{s,p}_w}^{\frac{p}{q}}
\,
\|u\|_{\dot{B}^{0}_{\infty,\infty}}^{1-\frac{p}{q}}.
\end{equation}
Recall that $BMO(\mathbb{R}^N)$ continuously embeds into $\dot{{B}}^{0}_{\infty,\infty}(\mathbb{R}^N)$. This inclusion holds because $BMO(\mathbb{R}^N)=\dot{F}^0_{\infty,2}(\mathbb{R}^N)$ \cite[5.2.4]{Triebel1983}, and $\dot{F}^0_{\infty,2}(\mathbb{R}^N)$ continuously embeds into $\dot{{B}}^{0}_{\infty,\infty}(\mathbb{R}^N)$ \cite[Lemma 3]{BenalliaMoussai2019}. Thus, we obtain
\begin{equation}
\label{eq:interpolation for BMO and Sobolev weak, our result}
\|u\|_{\dot{W}^{\frac{ps}{q},q}}
\le C \,
[u]_{W^{s,p}_w}^{\frac{p}{q}}
\,
\|u\|_{BMO}^{1-\frac{p}{q}}.
\end{equation}

{\it We note that \eqref{eq:interpolation for BMO and B0 infty infty} is more general than \eqref{eq:interpolation for BMO and Sobolev weak, our result}. However, it requires the assumption $p>1$ in order to apply \eqref{eq:B with parameters and Besov}, whereas we established \eqref{eq:interpolation for BMO and Sobolev weak, our result} in Section \ref{sec:interpolation of weak Gagliardo spaces} for all $p>0$.}

\subsection{(Local) $bmo$-interpolation in Besov spaces $B^s_{p,q}$} 
\label{sec: bmo(local) interpolation in Besov space}
In this section, we establish a $\mathrm{bmo}$ interpolation result in Besov spaces, using the definition of Besov spaces via the Littlewood--Paley decomposition. We also make use of Triebel--Lizorkin spaces.

The authors do not attempt to derive the $\mathrm{bmo}$ interpolation result from the main tools developed here, namely the John--Nirenberg inequality and its consequences. Instead, following the referee’s suggestion, we indicate an alternative approach based on the general theory of function spaces, in particular as presented in \cite{Triebel1983}. We provide the reader with the necessary definitions and refer to standard results available in the literature.

\begin{definition}[\cite{Triebel1983}, 2.3.1]
\label{def:dyadic decomposition}
Let $\Phi(\mathbb{R}^N)$ be the collection of all systems 
$\varphi = \{\varphi_j(x)\}_{j=0}^\infty \subset \mathcal{S}(\mathbb{R}^N)$ 
such that
\begin{equation}
\label{eq:supports in rings}
\operatorname{supp} \varphi_0 \subset \{ x \in \mathbb{R}^N : |x| \le 2 \},\quad 
\operatorname{supp} \varphi_j \subset 
\{ x \in \mathbb{R}^N : 2^{j-1} \le |x| \le 2^{j+1} \}, 
\quad j = 1,2,3,\dots,
\end{equation}
and for every multi-index $\alpha$ there exists a constant $c_\alpha > 0$ such that
\begin{equation}
2^{j|\alpha|} |D^\alpha \varphi_j(x)| \le c_\alpha
\quad \text{for all } j = 0,1,2,\dots \text{ and all } x \in \mathbb{R}^N,
\end{equation}
and
\begin{equation}
\sum_{j=0}^\infty \varphi_j(x) = 1
\quad \text{for every } x \in \mathbb{R}^N.
\end{equation}
The systems in $\Phi(\mathbb{R}^N)$ are called dyadic resolutions of unity.
\end{definition}

\begin{lemma}[\cite{Triebel1983}, 2.3.1, Remark 1]
\label{lem:diadic system from one function}
There exists $0 \le \varphi \in \mathcal{S}(\mathbb{R}^N)$ such that
\begin{equation}
\label{eq:supports in rings1}
\operatorname{supp}\varphi \subset \{x \in \mathbb{R}^N : \tfrac{1}{2} \le |x| \le 2\},\quad 
\varphi(x) > 0 \quad \text{for } \tfrac{1}{\sqrt{2}} \le |x| \le \sqrt{2},
\end{equation}
and $\varphi_j(x) := \varphi(2^{-j}x)$ lies in  $\Phi(\mathbb{R}^N)$.
\end{lemma}

\begin{proposition}
\label{prop:boundedness in L1 of inverse Fourier transform}
Let $\varphi \in \mathcal{S}(\mathbb{R}^N)$ and define 
\begin{equation}
\varphi_j(\xi) := \varphi(2^{-j}\xi), \qquad  0\leq j\in\N.
\end{equation}
Let $K_j := \mathcal{F}^{-1}(\varphi_j)$ and $K := \mathcal{F}^{-1}(\varphi)$.
Then
\begin{equation}
\label{eq:L1 norm of phij is stable}
\|K_j\|_{L^1(\mathbb{R}^N)}=\|K\|_{L^1(\mathbb{R}^N)}.
\end{equation}
\end{proposition}

\begin{proof}
By the definition of the inverse Fourier transform,
\begin{equation}
K_j(x)
=
\int_{\mathbb{R}^N} e^{ix\cdot \xi}\, \varphi(2^{-j}\xi)\, d\xi.
\end{equation}
We perform the change of variables $\eta = 2^{-j}\xi,\, d\xi = 2^{jN} d\eta$.
This gives
\begin{equation}
K_j(x)
=
2^{jN} \int_{\mathbb{R}^N} e^{i(2^j x)\cdot \eta}\, \varphi(\eta)\, d\eta=2^{jN} K(2^j x),\quad K := \mathcal{F}^{-1}(\varphi).
\end{equation}
Using the change of variables $y = 2^j x$, $dx = 2^{-jN} dy$, we obtain
\begin{equation}
\|K_j\|_{L^1}
=\int_{\mathbb{R}^N} |2^{jN} K(2^j x)|\, dx=
\int_{\mathbb{R}^N} |K(y)|\, dy
=
\|K\|_{L^1}.
\end{equation}
\end{proof}

\begin{definition}
\label{def:definition of Besov space for real s}
Let $s \in \mathbb{R}$, $0 < p \le \infty$.

\textbf{(i) Case $0 < q < \infty$:}
\[
B^s_{p,q}(\mathbb{R}^N)
:=
\left\{
f \in \mathcal{S}'(\mathbb{R}^N)
:
\left(
\sum_{j=0}^{\infty} 2^{jsq}\|\mathcal{F}^{-1} \varphi_j \mathcal{F} f\|_{L^p}^q
\right)^{1/q}
< \infty
\right\}.
\]

\textbf{(ii) Case $q = \infty$:}
\[
B^s_{p,\infty}(\mathbb{R}^N)
:=
\left\{
f \in \mathcal{S}'(\mathbb{R}^N)
:
\sup_{j \ge 0} 2^{js}\|\mathcal{F}^{-1} \varphi_j \mathcal{F} f \|_{L^p}
< \infty
\right\}.
\]
\end{definition}

Recall that if $f\in \mathcal{S}'(\mathbb{R}^N)$ and $\varphi\in \mathcal{S}(\mathbb{R}^N)$, then $\mathcal{F}^{-1}\bigl(\varphi\,\mathcal{F}f\bigr)$
defines a smooth function on $\mathbb{R}^N$. More precisely,
\begin{equation}
\mathcal{F}^{-1}\bigl(\varphi\,\mathcal{F}f\bigr)
=
\mathcal{F}^{-1}\!\left[\mathcal{F}\bigl(\mathcal{F}^{-1}\varphi\bigr)\,\mathcal{F}f\right]
=
\bigl(\mathcal{F}^{-1}\varphi\bigr)*f,
\end{equation}
where the convolution is defined by
\begin{equation}
\bigl((\mathcal{F}^{-1}\varphi)*f\bigr)(x)
:=
\left\langle f,\; (\mathcal{F}^{-1}\varphi)(x-\cdot)\right\rangle,
\qquad x\in \mathbb{R}^N,
\end{equation}
and belongs to $C^\infty(\mathbb{R}^N)$.
\begin{definition}
If $-\infty < s < \infty$, $1 < q \le \infty$ and 
$\varphi = \{\varphi_k(x)\}_{k=0}^{\infty} \in \Phi(\mathbb{R}^N)$, then
\begin{multline}
\label{eq:def of F infty}
F_{\infty,q}^s(\mathbb{R}^N)
:=\Big\{f\,: f \in \mathcal{S}'(\mathbb{R}^N),\ \exists \{f_k\}_{k=0}^{\infty} \subset L^\infty(\mathbb{R}^N) \quad \text{such that}
\\
f = \sum_{k=0}^{\infty} \mathcal{F}^{-1} \varphi_k \mathcal{F} f_k
\ \text{in } \mathcal{S}'(\mathbb{R}^N),
\ \text{and }
\left\| 2^{sk} f_k \mid L^\infty(\mathbb{R}^N, \ell_q) \right\| < \infty
\Big\};
\end{multline}
and 
\begin{equation}
\|f\|_{F_{\infty,q}^s(\mathbb{R}^N)}:=\inf \left\| 2^{sk} f_k \mid L^\infty(\mathbb{R}^N, \ell_q) \right\|,
\end{equation}
where the infimum is taken over all admissible representations of $f$ in the sense of \eqref{eq:def of F infty}. Here

\[
\left\| 2^{sk} f_k \mid L^\infty(\mathbb{R}^N,\ell_q) \right\|
:= \operatorname*{ess\,sup}_{x \in \mathbb{R}^N}
\left(
\sum_{k=0}^\infty |2^{sk} f_k(x)|^q
\right)^{1/q},
\qquad 1 < q < \infty.
\]

\[
\left\| 2^{sk} f_k \mid L^\infty(\mathbb{R}^N,\ell_\infty) \right\|
:= \operatorname*{ess\,sup}_{x \in \mathbb{R}^N}
\sup_{k \ge 0} |2^{sk} f_k(x)|.
\]
\end{definition}

\begin{proposition}
\label{prop:TL sub Besov}
Let $s\in\mathbb{R}$ and $1 < q \le \infty$. Then
\begin{equation}
F^s_{\infty,q}(\mathbb{R}^N)\hookrightarrow B^s_{\infty,\infty}(\mathbb{R}^N),
\end{equation}
i.e., there exists a constant $C>0$ such that
\begin{equation}
\|f\|_{B^s_{\infty,\infty}}
\le C \|f\|_{F^s_{\infty,q}}
\quad \text{for all } f \in F^s_{\infty,q}(\mathbb{R}^N).
\end{equation}
\end{proposition}

\begin{proof}
Let $f \in F^s_{\infty,q}(\mathbb{R}^N)$. By definition, there exists a sequence
$\{f_k\}_{k=0}^\infty \subset L^\infty(\mathbb{R}^N)$ such that
\begin{equation}
f = \sum_{k=0}^\infty \mathcal{F}^{-1}(\varphi_k \mathcal{F}{f_k})
\quad \text{in } \mathcal{S}'(\mathbb{R}^N),
\end{equation}
where $\varphi_k$ are as in Lemma \ref{lem:diadic system from one function},
and
\begin{equation}
\left\| 2^{ks}f_k \mid L^\infty(\mathbb{R}^N,\ell_q) \right\|
:=
\operatorname*{ess\,sup}_{x\in\mathbb{R}^N}
\left(\sum_{k=0}^\infty |2^{ks}f_k(x)|^q\right)^{1/q}
<\infty.
\end{equation}

Set $K_j := \mathcal{F}^{-1}(\varphi_j)$. Then
\begin{equation}
\mathcal{F}^{-1}(\varphi_j \mathcal{F} f) = K_j * f.
\end{equation}
By the definition of $\varphi_j$ as in Definition \ref{def:dyadic decomposition}, there exists a constant $C_0$ (one can take $C_0=3$) such that 
\begin{equation}
\varphi_j \varphi_k \equiv 0 \quad \text{if } |j-k|>C_0.
\end{equation}
Thus, for $j,k\in \mathbb{N}$ such that $|k-j|>C_0$ we have
\begin{equation}
0=\mathcal{F}^{-1}(0)=\mathcal{F}^{-1}(\varphi_j \varphi_k)=\mathcal{F}^{-1}(\mathcal{F}\left(\mathcal{F}^{-1}\varphi_j\right) \mathcal{F}\left(\mathcal{F}^{-1}\varphi_k\right) )=K_j * K_k.
\end{equation}

Using the representation of $f$, we obtain for every $j\in \mathbb{N}$ and $x\in\mathbb{R}^N$
\begin{multline}
\label{eq:representation of convolution Kj and f}
K_j * f(x)
=
K_j * \left(\sum_{k=0}^\infty K_k * f_k\right)(x)
=
\sum_{k=0}^\infty K_j * (K_k * f_k)(x)
\\
=
\sum_{k=0}^\infty (K_j * K_k) * f_k(x)=
\sum_{|k-j|\le C_0} (K_j * K_k) * f_k(x).
\end{multline}
For each $j,k\in \mathbb{N}$, by Young's convolution inequality and Proposition \ref{prop:boundedness in L1 of inverse Fourier transform},
\begin{equation}
\|(K_j * K_k)*f_k\|_{L^\infty}
\le \|K_j * K_k\|_{L^1}\,\|f_k\|_{L^\infty}
\le \|K_j\|_{L^1}\|K_k\|_{L^1}\,\|f_k\|_{L^\infty}
= C\,\|f_k\|_{L^\infty},
\end{equation}
where $C:=\|K\|^2_{L^1}$.
Therefore,
\begin{equation}
\|K_j * f\|_{L^\infty}
\le
C \sum_{|k-j|\le C_0} \|f_k\|_{L^\infty}.
\end{equation}

Multiplying by $2^{js}$, we obtain
\begin{equation}
2^{js}\|K_j * f\|_{L^\infty}
\le
C \sum_{|k-j|\le C_0} 2^{js}\|f_k\|_{L^\infty}.
\end{equation}
Since $|k-j|\le C_0$, we have
\begin{equation}
2^{js}=2^{(j-k)s}2^{ks}\le 2^{|j-k||s|}2^{ks}\leq C_s\,2^{ks},
\end{equation}
where $C_s:=2^{|s|C_0}$. Note that for every $j$ there exists at most $1+2C_0$ integers $k$ such that $|k-j|\leq C_0$. Hence
\begin{equation}
2^{js}\|K_j * f\|_{L^\infty}
\le
C C_s \sum_{|k-j|\le C_0} 2^{ks}\|f_k\|_{L^\infty}
\le
C C_s(1+2C_0)\sup_k 2^{ks}\|f_k\|_{L^\infty}.
\end{equation}

Set $C':=C C_s(1+2C_0)$. Taking the supremum over $j\ge0$, we obtain
\begin{equation}
\|f\|_{B^s_{\infty,\infty}}
=
\sup_{j\ge0} 2^{js}\|\mathcal{F}^{-1}(\varphi_j \mathcal{F} f)\|_{L^\infty}
=
\sup_{j\ge0} 2^{js}\|K_j * f\|_{L^\infty}
\le
C' \sup_k 2^{ks}\|f_k\|_{L^\infty}.
\end{equation}
Finally, since for every $j\in \N\cup\{0\}$
\begin{equation}
 2^{js}\|f_j\|_{L^\infty}
\le
\operatorname*{ess\,sup}_{x\in\mathbb{R}^N}
\left(\sum_{k=0}^\infty |2^{ks}f_k(x)|^q\right)^{1/q}
=
\left\| 2^{ks}f_k \mid L^\infty(\mathbb{R}^N,\ell_q) \right\|,
\end{equation}
it follows that
\begin{equation}
\|f\|_{B^s_{\infty,\infty}}
\le
C' \left\| 2^{ks}f_k \mid L^\infty(\mathbb{R}^N,\ell_q) \right\|.
\end{equation}
Taking the infimum over all admissible representations of $f$ yields
\begin{equation}
\|f\|_{B^s_{\infty,\infty}}
\le
C' \|f\|_{F^s_{\infty,q}}.
\end{equation}
This completes the proof.
\end{proof}

Recall the Besov--Lorentz spaces, denoted by $B^s_{q}L^{p,r}(\mathbb{R}^N)$, where $s\in\mathbb{R}$ and $p,q,r\in(0,\infty]$ ($L^{\infty,r}(\R^N)=\{0\}$ for $r<\infty$). They are defined analogously to Besov spaces, as in Definition \ref{def:definition of Besov space for real s}, with the $L^p$ quasi-norm replaced by the Lorentz quasi-norm $L^{p,r}$. Note that
\begin{equation}
\label{eq:Besov Lorentz coinside with Lebesgue}
B^s_{q}L^{p,p}(\mathbb{R}^N)=B^s_{p,q}(\mathbb{R}^N).
\end{equation}

Let $p_0,p_1\in [1,\infty]$, $q_0,q_1\in (0,\infty]$, $s_0,s_1\in \R$, $0<\theta<1$, $0<r\leq \infty$. By \cite[Page 106, Theorem 6]{Peetre1976}
\begin{equation}
\label{eq:inclusion in Besov-Lorentz space}
\left(B^{s_0}_{p_0,q_0}, B^{s_1}_{p_1,q_1}\right)_{\theta,r}
\hookrightarrow
B^{s}_{\max\{q,r\}} L^{p,r},
\end{equation}
where
\[
\frac{1}{p} = \frac{1-\theta}{p_0} + \frac{\theta}{p_1}, \qquad
\frac{1}{q} = \frac{1-\theta}{q_0} + \frac{\theta}{q_1}, \qquad
s = (1 - \theta)s_0 + \theta s_1.
\]

Note that, for quasi-normed spaces $X_0,X_1,X_1'$, if $X_1\hookrightarrow X_1'$, then the pair $(X_0,X_1)_{\theta,\gamma}$ is continuously embedded into $(X_0,X_1')_{\theta,\gamma}$, where $0<\theta<1$, $0<\gamma \le \infty$. 

It follows from the inequality $K(t,f,X_0,X_1')\leq C K(t,f,X_0,X_1)$ for $t>0$ and $f\in X_0+X_1$, where $C$ does not depend on $t$ and $f$. This inequality holds for $f\in X_0+X_1$:
\begin{multline}
K(t,f; X_0, X'_1)
:= \inf\limits_{f=f_0+f_1}\left\{ \|f_0\|_{X_0} + t\,\|f_1\|_{X'_1} \right\}
\leq \inf\limits_{f=f_0+f_1}\left\{ \|f_0\|_{X_0} + tC\,\|f_1\|_{X_1} \right\}
\\
\leq \max\{1,C\}\inf\limits_{f=f_0+f_1}\left\{ \|f_0\|_{X_0} + t\,\|f_1\|_{X_1} \right\}=\max\{1,C\} K(t,f; X_0, X_1).
\end{multline} 

\begin{definition}
The space $\mathrm{bmo}(\mathbb{R}^N)$ consists of all $f \in L^1_{\mathrm{loc}}(\mathbb{R}^N)$ such that
\[
\|f\|_{\mathrm{bmo}}
:=
\sup_{\mathcal{L}^N(Q)\le 1}
\frac{1}{\mathcal{L}^N(Q)}
\int_Q |f(x)-f_Q|dx
+
\sup_{\mathcal{L}^N(Q)>1}
\frac{1}{\mathcal{L}^N(Q)}
\int_Q |f(x)|dx
< \infty,
\]
where $Q\subset\mathbb{R}^N$ ranges over cubes with sides parallel to the axes, and
\[
f_Q := \frac{1}{\mathcal{L}^N(Q)}\int_Q f(x)dx.
\]
\end{definition} 
By \cite[2.5.8, Theorem 2]{{Triebel1983}} and Proposition \ref{prop:TL sub Besov}, we obtain $bmo(\mathbb{R}^N)=F^0_{\infty,2}(\R^N)\hookrightarrow  B^{0}_{\infty,\infty}(\R^N)$. 

Let $s\in\R$, $1\leq w<p<\infty$, $r,q\in (0,\infty]$. 

Therefore, by \eqref{eq:inclusion in Besov-Lorentz space}, we get
\[
\left(B^{s}_{w,\frac{qw}{p}}, bmo\right)_{1-\frac{w}{p},r}
\hookrightarrow
\left(B^{s}_{w,\frac{qw}{p}}, B^{0}_{\infty,\infty}\right)_{1-\frac{w}{p},r}
\hookrightarrow
B^{s\frac{w}{p}}_{\max\{q,r\}} L^{p,r}.
\]
Recall that for $0<p\leq \infty$, $0<r_1\leq r_2\leq \infty$, we have $L^{p,r_1}\hookrightarrow L^{p,r_2}$ \cite[Theorem 6.3]{CastilloRafeiro2016}. In particular, $B^s_qL^{p,r_1}\hookrightarrow B^s_q L^{p,r_2}$, for $s\in \R$, $q\in (0,\infty]$. Setting $r = \min\{p, q\}$ one has by \eqref{eq:Besov Lorentz coinside with Lebesgue}
\[
B^{s\frac{w}{p}}_{\max\{q,\min\{p,q\}\}} L^{p,\min\{p,q\}}
=
B^{s\frac{w}{p}}_{q} L^{p,\min\{p,q\}}
\hookrightarrow
B^{s\frac{w}{p}}_{q} L^{p,p}=
B^{s\frac{w}{p}}_{p,q}.
\]
Hence
\[
\left(B^{s}_{w,w\frac{q}{p}}, bmo\right)_{1-\frac{w}{p},\min\{p, q\}}
\hookrightarrow
B^{s\frac{w}{p}}_{p,q}.
\]
Using Proposition \ref{prop:interpolation-trivial-general}, we get 
\begin{equation}
\|f\|_{B^{s\frac{w}{p}}_{p,q}}\leq C
\|f\|_{\left(B^{s}_{w,w\frac{q}{p}}, bmo\right)_{1-\frac{w}{p},\min\{p, q\}}}
\leq C
 \|f\|^{\frac{w}{p}}_{B^{s}_{w,w\frac{q}{p}}} \|f\|^{1-\frac{w}{p}}_{bmo}.
\end{equation}
We summarize this result in the following theorem:

\begin{theorem}
Let $s\in \R$, $q\in (0,\infty]$, $1\le w\leq p<\infty$. Then, there exists a constant $C$ such that 
for every $\mathcal{L}^N$-measurable function $f:\R^N\to \R$
\begin{equation}
\|f\|_{B^{s\frac{w}{p}}_{p,q}}\leq C
 \|f\|^{\frac{w}{p}}_{B^{s}_{w,w\frac{q}{p}}} \|f\|^{1-\frac{w}{p}}_{bmo}.
\end{equation}
\end{theorem}
The restriction $1 \le w$ comes solely from the use of the inclusion \eqref{eq:inclusion in Besov-Lorentz space}.

\subsection{Equivalence of Besov quasi-norms}

\begin{proposition}
\label{prop:the bigger the degree of differences the smaller it is in Lp norm}
Let $q\in (0,\infty)$. Given $u \in L^q(\mathbb{R}^N)$, $n = 1, 2, \ldots$ and $h \in \mathbb{R}^N$, consider
\begin{equation}
\label{eq:definition of Delta nhu}
\Delta_h^n u(x) := \sum_{j=0}^n (-1)^{n-j} \binom{n}{j} u(x + jh).
\end{equation}
Then, for every $n,m\in \N$ and $h\in\R^N$, we get
\begin{equation}
\label{eq:final,prop}
\|\Delta_h^{n+m} u\|_{L^q(\R^N)}
\leq 2^m \|\Delta_h^{n} u\|_{L^q(\R^N)},\quad\quad q\in [1,\infty);
\end{equation}
and
\begin{equation}
\label{eq:final1,prop}
\|\Delta_h^{n+m} u\|_{L^q(\R^N)}
\leq 2^{\frac{m}{q}} \|\Delta_h^{n} u\|_{L^q(\R^N)}\quad\quad q\in (0,1).
\end{equation}
In particular, for every $n,k\in \N$, $k\geq n$, we have for every $t\in (0,\infty)$
\begin{equation}
\Omega_k(f,t)_{L^q}\leq\max\left(2^{k-n},2^{\frac{k-n}{q}}\right) \Omega_n(f,t)_{L^q},\quad\quad \Omega_n(f,t)_{L^q} := \sup_{|h| \leq t} \|\Delta_h^n f\|_{L^q(\mathbb{R}^N)}.
\end{equation}
\end{proposition}

\begin{proof}
In case $q\in [1,\infty)$, we get by Minkowski inequality for the $L^q$ norm and change of variable:
\begin{equation}
\label{eq:ineq01}
\left\| \Delta_h^n u(\cdot + h) - \Delta_h^n u \right\|_{L^q(\R^N)}
\leq \left\| \Delta_h^n u(\cdot + h)\right\|_{L^q(\R^N)}+\left\|\Delta_h^n u \right\|_{L^q(\R^N)}=2\left\|\Delta_h^n u \right\|_{L^q(\R^N)}.
\end{equation}

In case $q\in (0,1)$, we get by the properties of the $L^q$ quasi-norm:
\begin{multline}
\label{eq:ineq011}
\left\| \Delta_h^n u(\cdot + h) - \Delta_h^n u \right\|_{L^q(\R^N)}
\leq 2^{\frac{1}{q}-1}\left(\left\| \Delta_h^n u(\cdot + h)\right\|_{L^q(\R^N)}+\left\|\Delta_h^n u \right\|_{L^q(\R^N)}\right)=2^{\frac{1}{q}}\left\|\Delta_h^n u \right\|_{L^q(\R^N)}.
\end{multline}

We now prove the identity:
\begin{equation}
\label{eq:identity for differences}
\Delta_h^n u(x + h) - \Delta_h^n u(x)=\Delta_h^{n+1} u(x)
\end{equation}
for almost every $x\in \R^N$.

Using the definition of $\Delta_h^n u$, we obtain for almost every $x\in \R^N$
\begin{multline}
\Delta_h^n u(x + h) - \Delta_h^n u(x)
= \sum_{j=0}^n (-1)^{n-j} \binom{n}{j} u((x+h) + jh)
   - \sum_{j=0}^n (-1)^{n-j} \binom{n}{j} u(x + jh) 
\\   
= \sum_{j=1}^{n+1} (-1)^{(n+1)-j} \binom{n}{j-1} u(x + jh)
   + \sum_{j=0}^n (-1)^{n+1-j} \binom{n}{j} u(x + jh) 
\\
= \binom{n}{n} u\big( x + (n+1)h \big)
   + \sum_{j=1}^n (-1)^{(n+1)-j} \left[ \binom{n}{j-1} + \binom{n}{j} \right] u(x + jh)
   + (-1)^{n+1} \binom{n}{0} u(x) \\
= \binom{n+1}{n+1} u\big( x + (n+1)h \big)
   + \sum_{j=1}^n (-1)^{(n+1)-j} \binom{n+1}{j} u(x + jh)
   + (-1)^{n+1} \binom{n+1}{0} u(x) 
\\
= \sum_{j=0}^{n+1} (-1)^{(n+1)-j} \binom{n+1}{j} u(x + jh)
= \Delta_h^{n+1} u(x).
\end{multline}

This shows that for $q\in [1,\infty)$ we get from \eqref{eq:ineq01} and \eqref{eq:identity for differences}
\begin{equation}
\label{eq:recurrence2}
 \left\|\Delta_h^{n+1} u\right\|_{L^q(\R^N)}
\leq 2 \left\|\Delta_h^{n} u\right\|_{L^q(\R^N)},
\end{equation}
and for $q\in (0,1)$ we get from \eqref{eq:ineq011} and \eqref{eq:identity for differences}
\begin{equation}
\label{eq:recurrence1}
 \left\|\Delta_h^{n+1} u\right\|_{L^q(\R^N)}
\leq 2^{\frac{1}{q}} \left\|\Delta_h^{n} u\right\|_{L^q(\R^N)}.
\end{equation}
Applying \eqref{eq:recurrence2} and \eqref{eq:recurrence1} $m$ times finally gives \eqref{eq:final,prop} and \eqref{eq:final1,prop}, respectively.
This completes the proof.
\end{proof}

{\textbf{Marchaud inequality}:
\begin{theorem}
\label{thm:Marchaud inequality non hom}
Let $p\in (0,\infty)$, $f\in L^p(\R^N)$, $n,k\in \N$, $n\leq k$. Then, for every $0<\mu\leq \min(1,p)$ and $t\in(0,\infty)$
\begin{equation}\label{ijjihjhj}
\Omega_{n}(f,t)_{L^p} \leq C\, t^{n} \left[ 
\|f\|_{L^p(\R^N)} +
\left( \int_{t}^{\infty} \left( s^{-n} \Omega_{k}(f,s)_{L^p} \right)^{\mu} \frac{ds}{s} \right)^{1/\mu} \right].
\end{equation}
The constant $C$ depends on $p,k,\mu$ only.
\end{theorem}
One can find a proof of Theorem \ref{thm:Marchaud inequality non hom} in \cite{DeVorePopov1988,Ditzian1988}.

\begin{corollary}[Homogenous Marchaud inequality]
\label{thm:Marchaud inequality}
Let $p\in (0,\infty)$, $f\in L^p(\R^N)$, $n,k\in \N\setminus\{0\}$, $n\leq k$. Then, for every $0<\mu\leq \min(1,p)$ and $t\in(0,\infty)$
\begin{equation}\label{ijjihjhj111}
\Omega_{n}(f,t)_{L^p} \leq C\, t^{n} \left[ 
\left( \int_{t}^{\infty} \left( s^{-n} \Omega_{k}(f,s)_{L^p} \right)^{\mu} \frac{ds}{s} \right)^{1/\mu} \right].
\end{equation}
The constant $C$ depends on $p,k,\mu$ only.
\end{corollary}
\begin{proof} 
The idea of the proof was taken from \cite{Kolomoitsev-Tikhonov}. Given $p\in (0,\infty)$, $f\in L^p(\R^N)$, $n\leq k\in \N$, $0<\mu\leq \min(1,p)$ and $t\in(0,\infty)$, for every $\varepsilon>0$ define $f_\e\in L^p(\R^N)$ by defining: $g_\e(x):=f(\e x)$, $f_\e(x):=\e^{\frac{N}{p}}g_\e(x)=\e^{\frac{N}{p}}f(\e x)$ and define $t_\e:=t/\e$. Then, by \eqref{ijjihjhj}, applied to $f_\e$ in the place of $f$ and $t_\e$ instead of $t$ we have
\begin{equation}\label{ijjihjhjjjkjk}
\Omega_{n}(f_\e,t_\e)_{L^p} \leq C\, t^{n}_\e \left[ 
\|f_\e\|_{L^p(\R^N)} +
\left( \int_{t_\e}^{\infty} \left( s^{-n} \Omega_{k}(f_\e,s)_{L^p} \right)^{\mu} \frac{ds}{s} \right)^{1/\mu} \right].
\end{equation}
Since $g_\e(x)=f(\e x)$, then $(\Delta^j_{(h/\e)} g_\e)(x)=(\Delta^j_{h} f)(\e x)$ for every $h\in\R^N$, $j\in\N$ and almost every $x\in\R^N$. Therefore, by change of variable $\e^{\frac{N}{p}}\|\Delta^j_{(h/\e)} g_\e\|_{L^p(\R^N)}=\|\Delta^j_{h} f\|_{L^p(\R^N)}$. Thus, for every $j\in \N$
\begin{multline}\label{ijjihjhjjjkjkklklk}
\Omega_{j}(f_\e,t_\e)_{L^p}:=\sup_{|h|\leq t_\e}\|\Delta^j_h f_\e\|_{L^p(\R^N)}=\sup_{|h|\leq t/\e}\|\Delta^j_h f_\e\|_{L^p(\R^N)}=\sup_{|h'|\leq t}\|\Delta^j_{(h'/\e)} f_\e\|_{L^p(\R^N)}\\
=\sup_{|h'|\leq t}\e^{\frac{N}{p}}\|\Delta^j_{(h'/\e)} g_\e\|_{L^p(\R^N)}=
\sup_{|h'|\leq t}\|\Delta^j_{h'} f\|_{L^p(\R^N)}=\Omega_{j}(f,t)_{L^p}.
\end{multline}
Similarly (with $(\e s)$ instead of $t$):
\begin{equation}\label{ijjihjhjjjkjkklklkjjk}
\Omega_{j}(f_\e,s)_{L^p}=\Omega_{j}(f,\e s)_{L^p}.
\end{equation}
Observe that
\begin{equation}\label{ijjihjhjjjkjkdffdjkjk}
\|f_\e\|^p_{L^p(\R^N)}=\e^N\|g_\e\|^p_{L^p(\R^N)}=\|f\|^p_{L^p(\R^N)},
\end{equation}
Thus, we rewrite
\eqref{ijjihjhjjjkjk} as 
\begin{equation}\label{ijjihjhjjjkjkjlkjjjjkjjbn}
\Omega_{n}(f,t)_{L^p} \leq C\, t^{n} \left[ 
\frac{1}{\e^n}\|f\|_{L^p(\R^N)} +
\frac{1}{\e^n}\left( \int_{t/\e}^{\infty} \left( s^{-n} \Omega_{k}(f,\e s)_{L^p} \right)^{\mu} \frac{ds}{s} \right)^{1/\mu} \right].
\end{equation}
Thus, changing variables of integration $\tau=\e s$ in \eqref{ijjihjhjjjkjkjlkjjjjkjjbn} gives
\begin{equation}\label{ijjihjhjjjkjkjlkjjjjkjjbnllkk}
\Omega_{n}(f,t)_{L^p} \leq C\, t^{n} \left[ 
\frac{1}{\e^n}\|f\|_{L^p(\R^N)} +
\left( \int_{t}^{\infty} \left( \tau^{-n} \Omega_{k}(f,\tau)_{L^p} \right)^{\mu} \frac{d\tau}{\tau} \right)^{1/\mu} \right].
\end{equation}
Finally, letting $\e\to \infty$ in \er{ijjihjhjjjkjkjlkjjjjkjjbnllkk} implies \eqref{ijjihjhj111}.
\end{proof}

\begin{theorem}
\label{thm:equivalence of Besov quasi-norms}
Let $\alpha \in (0,\infty)$, $p,q \in (0,\infty]$, $n,k \in \mathbb{N}$ with $n,k > \alpha$, and let $f \in L^p(\mathbb{R}^N)$. Then, the following quasi-norms are equivalent:

\medskip
If $q < \infty$,
\begin{equation}
\label{eq:equivalence for Besov norms in case of finite q}
\left( \displaystyle\int_0^\infty  \left( \frac{\Omega_n(f,t)_{L^p}}{t^{\alpha}} \right)^{q} \, \frac{dt}{t} \right)^{\! 1/q} \sim
\left( \displaystyle\int_0^\infty \left( \frac{\Omega_k(f,t)_{L^p}}{t^{\alpha}} \right)^{q} \, \frac{dt}{t}  \right)^{\! 1/q}.
\end{equation}

If $q = \infty$,
\begin{equation}
\label{eq:equivalence for Besov norms in case of infinite q}
\sup_{t\in (0,\infty)} \frac{\Omega_n(f,t)_{L^p}}{t^{\alpha}} \sim
\sup_{t\in (0,\infty)} \frac{\Omega_k(f,t)_{L^p}}{t^{\alpha}}.
\end{equation}
\end{theorem}

Recall the Hardy inequality, see, for example, Lemma 1 in \cite{Mironescu2018}.

{\textbf{Hardy's inequality}: If \( 1\leq q<\infty \), \( 0<r<\infty \), and \( h \) is a non-negative $\mathcal{L}^1$-measurable function defined on \( (0, \infty) \), then:
\begin{equation}
 \int_0^\infty t^{r - 1}\left( \int_t^\infty h(z) \,  \, dz \right)^q \, dt  \leq \left(\frac{q}{r}\right)^q \int_0^\infty t^{r+q-1} \left(h(t)\right)^q dt .
\end{equation}}

\begin{proof}[Proof of Theorem \ref{thm:equivalence of Besov quasi-norms}] 
Notice that, by the definition of $\Omega_{n}(f,t)_{L^p}$ and the change of variable formula, we have
\begin{multline}
\label{eq:estimate of modulus of continuity by Lp norm}
\begin{cases}
\Omega_{n}(f,t)_{L^p}\leq 2^n\|f\|_{L^p(\R^N)},\quad &p\geq 1,
\\
\Omega_{n}(f,t)_{L^p}\leq \left(\sum_{j=0}^n \binom{n}{j}^p\right)^{1/p}\|f\|_{L^p(\R^N)},\quad & p<1,
\end{cases}
\quad C_0:=\max\left(2^n,\left(\sum_{j=0}^n \binom{n}{j}^p\right)^{1/p}\right).
\end{multline}
We prove first of all the equivalence \eqref{eq:equivalence for Besov norms in case of finite q}, so we assume that $q<\infty$. 
Assume that $k>n$. Let $0<\mu\leq\min(1,p)$. From Corollary \ref{thm:Marchaud inequality}, we have for every $t\in (0,\infty)$
\begin{equation}
\label{eq:Marchaud inequality1}
\Omega_n(f,t)_{L^p}\leq
Ct^{n}\left(\int_t^\infty\left(s^{-n}\Omega_{k}(f,s)_{L^p}\right)^{\mu}
\frac{ds}{s}\right)^{\frac{1}{\mu}}\,.
\end{equation}
Therefore, \eqref{eq:Marchaud inequality1}, we get
\begin{multline}
\label{eq:appyling Hardy}
\int_0^\infty\left(\frac{\Omega_n(f,t)_{L^p}}{t^{\alpha}}\right)^q\frac{dt}{t}\leq 
C\int_0^\infty\left(t^{n-\alpha}\left(\int_t^\infty\left(s^{-n}\Omega_{k}(f,s)_{L^p}\right)^{\mu}
\frac{ds}{s}\right)^{\frac{1}{\mu}}\right)^q\frac{dt}{t}
\\
=C\int_0^\infty t^{(n-\alpha)q-1}\left(\int_t^\infty\frac{1}{s}\left(s^{-n}\Omega_{k}(f,s)_{L^p}\right)^{\mu}
ds\right)^{\frac{q}{\mu}}dt.
\end{multline}
We choose $\mu\leq q$ and apply Hardy's inequality and $n>\alpha$ to obtain
\begin{multline}
\label{eq:second decomposition}
\int_0^\infty t^{(n-\alpha)q-1}\left(\int_t^\infty\left(s^{-n}\Omega_{k}(f,s)_{L^p}\right)^{\mu}
\frac{ds}{s}\right)^{\frac{q}{\mu}}dt
\\
\leq  \left(\frac{\frac{q}{\mu}}{(n-\alpha)q}\right)^{\frac{q}{\mu}} \int_0^\infty t^{(n-\alpha)q+\frac{q}{\mu}-1} \left[\frac{1}{t}\left(t^{-n}\Omega_{k}(f,t)_{L^p}\right)^{\mu}
\right]^{\frac{q}{\mu}}dt
\\
=\left(\frac{\frac{q}{\mu}}{(n-\alpha)q}\right)^{\frac{q}{\mu}} \int_0^\infty \left(\frac{\Omega_{k}(f,t)_{L^p}}{t^{\alpha}}\right)^{q}\frac{dt}{t}.
\end{multline}
Thus, by \eqref{eq:appyling Hardy} and \eqref{eq:second decomposition}, we obtain
\begin{equation}
\label{eq:first decomposition1}
\int_0^\infty\left(\frac{\Omega_n(f,t)_{L^p}}{t^{\alpha}}\right)^q\frac{dt}{t}
\leq
C_1\left(\frac{\frac{q}{\mu}}{(n-\alpha)q}\right)^{\frac{q}{\mu}} \int_0^\infty \left(\frac{\Omega_{k}(f,t)_{L^p}}{t^{\alpha}}\right)^{q}\frac{dt}{t}.
\end{equation}
Therefore, by \eqref{eq:first decomposition1}
\begin{equation}
\label{eq:first decomposition3}
\left(\int_0^\infty\left(\frac{\Omega_n(f,t)_{L^p}}{t^{\alpha}}\right)^q\frac{dt}{t}\right)^{\frac{1}{q}}
\leq C_1^{1/q}\left(\frac{\frac{q}{\mu}}{(n-\alpha)q}\right)^{\frac{1}{\mu}} \left(\int_0^\infty \left(\frac{\Omega_{k}(f,t)_{L^p}}{t^{\alpha}}\right)^{q}\frac{dt}{t}\right)^{\frac{1}{q}}\,.
\end{equation}
This proves that the left-hand side of~\eqref{eq:equivalence for Besov norms in case of finite q} is bounded above by a constant times the right-hand side of~\eqref{eq:equivalence for Besov norms in case of finite q}. The reverse inequality follows from Proposition~\ref{prop:the bigger the degree of differences the smaller it is in Lp norm}.

We now treat the equivalence in \eqref{eq:equivalence for Besov norms in case of infinite q}. Let $0<\mu\leq\min(1,p)$. For every $t\in (0,\infty)$, we obtain by $n>\alpha$ and \eqref{eq:Marchaud inequality1}
\begin{multline}
\label{eq:estimate for Omegan}
\frac{\Omega_n(f,t)_{L^p}}{t^\alpha}\leq
Ct^{n-\alpha}\left(\int_t^\infty\frac{1}{s}\left(\frac{1}{s^{n}}\Omega_{k}(f,s)_{L^p}\right)^{\mu}
ds\right)^{\frac{1}{\mu}}
\\
= Ct^{n-\alpha}\left(\int_t^\infty\frac{1}{s}\left(\frac{s^\alpha}{s^{n}}\frac{\Omega_{k}(f,s)_{L^p}}{s^\alpha}\right)^{\mu}
ds\right)^{\frac{1}{\mu}}
\leq Ct^{n-\alpha}\sup_{t<s<\infty}\frac{\Omega_{k}(f,s)_{L^p}}{s^\alpha}\left(\int_t^\infty\frac{1}{s^{1+(n-\alpha)\mu}}
ds\right)^{\frac{1}{\mu}}
\\
= Ct^{n-\alpha}\sup_{t<s<\infty}\frac{\Omega_{k}(f,s)_{L^p}}{s^\alpha}\left(\frac{1}{(n-\alpha)\mu\, t^{(n-\alpha)\mu}}
\right)^{\frac{1}{\mu}}
=\frac{C}{\left((n-\alpha)\mu\right)^{\frac{1}{\mu}}}
\sup_{t<s<\infty}\frac{\Omega_{k}(f,s)_{L^p}}{s^\alpha}.
\end{multline}
Therefore, from \eqref{eq:estimate for Omegan},
\begin{equation}
\label{eq:estimate for Omegan1}
\sup_{t\in (0,\infty)}\frac{\Omega_n(f,t)_{L^p}}{t^\alpha}\leq \frac{C}{\left((n-\alpha)\mu\right)^{\frac{1}{\mu}}}
\sup_{s\in (0,\infty)}\frac{\Omega_{k}(f,s)_{L^p}}{s^\alpha}.
\end{equation}
This proves that the left-hand side of \eqref{eq:equivalence for Besov norms in case of infinite q} is bounded above by a constant times the right-hand side of \eqref{eq:equivalence for Besov norms in case of infinite q}. The reverse estimate follows from Proposition \ref{prop:the bigger the degree of differences the smaller it is in Lp norm}. This completes the proof of the theorem.
\end{proof}

\section*{Acknowledgment}
The first author was supported by the Israel Science Foundation (grant No.\ 569/21).
\section*{Acknowledgment}
The authors thank the anonymous referee for helpful comments and for suggesting alternative approaches to some of the results in this work.

\end{document}